%% file: li-li-liu-liao.tex
\documentclass[preprint,12pt,3p]{elsarticle}
\usepackage{amssymb}
\usepackage{amsmath}
\usepackage{algorithm,algpseudocode}
\usepackage{txfonts}
\usepackage{afterpage}
\usepackage{lipsum}
\usepackage{graphicx}
\usepackage{epstopdf}
\usepackage{amsfonts}
\usepackage{amsopn}
\usepackage{color}
\usepackage{tabularray}
\usepackage{hyperref}

\newcommand{\bc}{\boldsymbol{c}}
\newcommand{\bz}{\boldsymbol{z}}
\newcommand{\bs}{\boldsymbol{s}}
\newcommand{\bt}{\boldsymbol{t}}

\newcommand{\ba}{\boldsymbol{\alpha}}
\newcommand{\bw}{\boldsymbol{\omega}}

\newcommand{\pE}{\mathbf{E}}
\newcommand{\pV}{\mathbf{V}}
\newcommand{\pL}{\mathcal{L}}
\newcommand{\pB}{\mathfrak{b}}
\newcommand{\ph}{\mathfrak{h}}
\newcommand{\pD}{\partial D}
\newcommand{\pjoi}{\partial_{j} D_{i}}
\newcommand{\mw}{\mathbf{W}}
\newcommand{\mN}{\mathbf{N}}

\newcommand{\mR}{\mathbb{R}}
\newcommand{\dsR}{\mathbb{R}}

\newcommand{\linter}{\Lambda}
\newcommand{\linteri}{\Lambda_{i}}
\newcommand{\lr}{\color{black}}
\newcommand{\Non}{N_{\rm {on}}}
\newcommand{\Noff}{N_{\rm {off}}}
\renewcommand{\div}{\nabla \cdot}

\newdefinition{definition}{Definition}
\newtheorem{theorem}{Theorem}
\newtheorem{lemma}[theorem]{Lemma}
\newproof{proof}{Proof}

\makeatletter
\def\ps@pprintTitle{%
   \let\@oddhead\@empty
   \let\@evenhead\@empty
   \let\@oddfoot\@empty
   \let\@evenfoot\@oddfoot
}

\journal{}

\begin{document}

\begin{frontmatter}

\title{A conditional normalizing flow for domain decomposed uncertainty quantification} 

\author[1]{Sen Li\texorpdfstring{\fnref{equal}}{equal}} 
\ead{lisen@shanghaitech.edu.cn}
\author[1,2]{Ke Li\texorpdfstring{\fnref{equal}}{equal}}
\ead{like1@shanghaitech.edu.cn}
\author[1]{Yu Liu}
\ead{liuyu@shanghaitech.edu.cn}
\author[1]{Qifeng Liao\texorpdfstring{\corref{cor1}}{cor1}}
\ead{liaoqf@shanghaitech.edu.cn}

\affiliation[1]{
    organization={School of Information Science and Technology, ShanghaiTech University}, 
    city={Shanghai 201210},
    country={China}
}
\affiliation[2]{
    organization={Shanghai Alliance Investment Ltd}, 
    city={Shanghai 201210},
    country={China}
}

\cortext[cor1]{Corresponding author.}
\fntext[equal]{These authors contributed equally to this work.}

\begin{abstract}
In this paper we present a conditional KRnet (cKRnet) based domain decomposed uncertainty quantification (CKR-DDUQ) approach to propagate uncertainties across different physical domains in models governed by partial differential equations (PDEs) with random inputs. This approach is based on the domain decomposed uncertainty quantification (DDUQ) method presented in [Q. Liao and K. Willcox, SIAM J. Sci. Comput., 37 (2015), pp. A103--A133], which suffers a bottleneck of density estimation for local joint input distributions in practice. In this work, we reformulate the required joint distributions as conditional distributions, and propose a new conditional normalizing flow model, called cKRnet, to efficiently estimate the conditional probability density functions. We present the general framework of CKR-DDUQ, conduct its convergence analysis, validate its accuracy and demonstrate its efficiency with numerical experiments. 
\end{abstract}


\end{frontmatter}

\section{Introduction}\label{section_intro}
\input{sections/section1}

\section{Problem setup and domain decomposed uncertainty quantification}\label{section_pre}
\input{sections/section2}

\section{Conditional KRnet and conditional density estimation}\label{section_ckrnet}
\input{sections/section3}

\section{Conditional KRnet for DDUQ (CKR-DDUQ)}\label{section_method}
\input{sections/section4}

\section{Numerical study}\label{section_experiments}
\input{sections/section5}

\section{Conclusion}\label{section_conclude}
\input{sections/section6}

\section*{Acknowledgments}
This work is supported by the National Natural Science Foundation of China (No. 12071291) and the Model Reduction Theory and Algorithms for Complex Systems Foundation of Institute of Mathematical Sciences, ShanghaiTech University (No. 2024X0303-902-01).

\bibliographystyle{elsarticle-num}
\bibliography{ckr_dduq_ref}
\end{document}

%% file: sections/section1.tex
During the last two decades there has been a rapid development in numerical methods for uncertainty quantification (UQ). This explosion in interest has been driven by the need to increase the reliability of numerical simulations for complex engineering problems, e.g., partial differential equations (PDEs) with random inputs, inverse problems and optimization under uncertainty. In many UQ related areas, such as multidisciplinary optimization \cite{diedrich2006multidisciplinary}, subsurface hydrology \cite{tar13}, porous media flow modelling \cite{cheung2020deep} and microscopic simulations \cite{lei2010direct}, an essential problem is uncertainty propagation across different physical domains. For this purpose, decomposition based UQ methods have been actively developed. These include, for example, hybrid particle and continuum models \cite{aletar02}, parallel domain decomposition methods for stochastic Galerkin \cite{cai07}, random domain decomposition with sparse grid collocation \cite{lin10}, coupling dimension reduction for polynomial chaos expansion \cite{arngha11a,ghaspa03},
stochastic model reduction  based on low-rank separated representation \cite{doo14}, local polynomial chaos expansion \cite{xiu15,xiu2010numerical}, importance sampling based decomposition \cite{amaral2014decomposition,liao2015domain}, and propagation methods of uncertainty across heterogeneous domains \cite{HeyrimCho15}. In \cite{DongkunZhang18}, stochastic domain decomposition via moment minimization is proposed, where the moment minimizing interface condition is introduced  to match the stochastic local solutions. A domain decomposition model reduction method for stochastic convection-diffusion equations is studied in  \cite{mu2019domain}, and a network UQ method for component-based systems is introduced in  \cite{carlberg2019network}. In addition,  stochastic domain decomposition based on variable-separation is proposed in \cite{chen2024stochastic}, and domain decomposition for Bayesian inversion is studied in \cite{xu2024domain}. 

In this work, we focus on the importance sampling based decomposition for uncertainty analysis \cite{amaral2014decomposition,liao2015domain}. In \cite{amaral2014decomposition}, a distributed uncertainty analysis method for feed-forward multicomponent systems is introduced, which conducts local uncertainty analyses for individual local components and synthesizes  them into a global system uncertainty analysis through importance sampling. Based on this idea, a domain-decomposed uncertainty quantification approach (DDUQ) for systems governed by PDEs is proposed in \cite{liao2015domain}. In DDUQ, non-overlapping domain decomposition is applied to decompose global PDEs into local problems posed on physical subdomains, and uncertainty analyses are performed for local problems with corresponding proposal probability density functions (PDFs). The target probability density function of the inputs for each local problem is estimated  through a density estimation technique, where samples are generated through domain decomposition iterations with coupling surrogates. Then, importance sampling is utilized to assemble global uncertainty analysis results.

Density estimation is a crucial step in DDUQ. Over the last few decades, various density estimation methods have been developed. The classical density estimation methods can be broadly categorized into parametric and non-parametric approaches. In parametric methods, the given data are assumed to follow specific distributions, such as Gaussian distributions and mixed Gaussian distributions \cite{redner1984mixture, bishop2006pattern}, and parameters in these distributions are typically determined through maximum likelihood estimation and expectation maximization \cite{DeGroot2012probability}. This type of approaches is effective when the specified distributions can fit the underlying data well; otherwise, their performance degrades. For the non-parametric methods, they offer greater flexibility to fit the underlying data than parametric methods, as they are not restricted to a predetermined form of distributions. Techniques such as histogram estimation, kernel density estimation (KDE) \cite{hansen2008uniform}, and k-nearest neighbor methods \cite{peterson2009k} are popular in this category. Although achieving great success, challenges still exist for these classic methods when the underlying data are high-dimensional. Recently, new deep learning-based techniques have shown great potential for high-dimensional density estimation, e.g., the autoregressive models \cite{van2016pixel, papamakarios2017masked}, variational autoencoders (VAEs) \cite{kingma2013auto}, generative adversarial networks (GANs) \cite{goodfellow2014generative}, normalizing flow models \cite{dinh2014nice, dinh2016density, kingma2018glow, tang2020deep} and diffusion models \cite{ho2020denoising}. In this work, we are concentrated on KRnet \cite{tang2020deep}, a normalizing flow model that incorporates the Knothe-Rosenblatt (KR) rearrangement \cite{carlier2010knothe}, and it is applied to solve high-dimensional Fokker-Planck equations in  \cite{tang2022adaptive}. The KRnet is constructed as an invertible block-triangular mapping between the target distribution and a prior distribution, and the estimated PDF can then be written explicitly by the change of variables.  

The aim of this work is to systematically improve DDUQ \cite{liao2015domain} through utilizing deep learning methods for both density estimation and coupling surrogate modelling. In the offline stage of DDUQ, each subdomain is assigned a local proposal input PDF, and Monte Carlo simulations are performed to generate samples of local outputs and coupling functions. Based on the samples, surrogates of the coupling functions, which are referred to as the coupling surrogates, are built using deep neural networks (DNNs) in this work. In the online stage, a large number of samples are generated from the given PDF of global system input parameters, and domain decomposition iterations are conducted with the coupling surrogates to compute target local input samples. While the importance sampling procedure is applied to re-weight the local output samples, the weights are computed based on conditional PDFs of the local input parameters in this work, rather than the joint PDFs in the original DDUQ. Estimating conditional PDFs is inspired by the specific structure of the target local input samples. The local input parameters are divided into system input parameters and interface parameters, where interface parameters are dependent on system input parameters (details are shown in Section \ref{ddqu-alg}). Furthermore, the PDFs of the system input parameters are given, and focusing on the conditional PDFs can avoid unnecessary complexity to estimate the joint PDFs. However, efficient estimation for conditional distributions is another challenge. To address this issue, we generalize KRnet to a conditional version, which is referred to as conditional KRnet (cKRnet) in the following. To summarize, the main contributions of this work are two-fold: first we propose cKRnet to estimate general conditional PDFs; second, we modify the importance sampling procedure in DDUQ through replacing the joint PDFs of local inputs by conditional PDFs, and efficiently apply cKRnet to estimate the conditional PDFs for the local inputs. This new version of DDUQ is referred to as cKRnet based domain decomposed uncertainty quantification (CKR-DDUQ) approach in the following.  

The rest of the paper is organized as follows. Section~\ref{section_pre} presents the problem setup, reviews the DDUQ approach following \cite{liao2015domain}, where the  weights in the importance sampling procedure are reformulated. Then we propose the conditional KRnet (cKRnet) for conditional density estimation in Section~\ref{section_ckrnet}.  Section~\ref{section_method} presents our CKR-DDUQ approach and its convergence analysis. In Section~\ref{section_experiments}, we demonstrate the effectiveness of CKR-DDUQ with numerical experiments, where both diffusion and Stokes problems are studied. Finally Section~\ref{section_conclude} concludes the paper.

%% file: sections/section2.tex
Let $D\in\dsR^{d}$ $(d=1,2,3)$ denote a physical domain which is bounded, connected and with a polygonal boundary $\partial D$, and let $x\in \dsR^{d}$ denote a physical variable. Letting $\xi$ denote a vector collecting a finite number of random variables, the support of $\xi$ is denoted by $\Gamma$ and its probability density function (PDF) is denoted by $\pi_{\xi}(\xi)$. Our physics of problems considered are governed by a PDE over the domain $D$ and boundary conditions on the boundary $\partial D$, which are stated as: find a function $u(x,\xi)$ mapping $D\times \Gamma$ to  $\dsR$, such that
\begin{eqnarray}
	\pL\left(x,\xi;u\left(x,\xi\right)\right)=f(x,\xi)  \qquad
	&\forall& \left(x,\xi\right) \in D\times \Gamma,
	\label{spdexi1}\\
	\pB\left(x,\xi;u\left(x,\xi\right)\right)=g(x,\xi) \qquad
	&\forall& \left(x,\xi\right)\in \partial D\times \Gamma,
	\label{spdexi2}
\end{eqnarray}
where $\pL$ and $\pB$ denote a partial differential operator and a boundary condition operator respectively, $f$ is the source function and $g$ specifies the boundary conditions. The operators $\pL$ and $\pB$, and the functions $f$ and $g$ can depend on random variables. 

In the rest of this section, we review the domain decomposed uncertainty quantification (DDUQ) approach \cite{liao2015domain}, where the main challenge is to estimate conditional density functions in an importance sampling procedure.   

\subsection{Settings for non-overlapping domain decomposition}
\label{dd}
We assume that the physical domain $D$ is decomposed into a finite number, $M$ of non-overlapping subdomains (components), i.e., $D=\cup^{M}_{i=1}D_i$, and the intersection of two subdomains is a connected interface with positive $(d-1)$-dimensional measure or an empty set. For each subdomain $D_i$, $\partial D_i$ denotes the set of its boundaries, and $\linteri$ denotes the set of its neighboring subdomain indices, i.e.,  $\linteri:=\{j \,|\, j\in\{1,\dots,M\}, j \neq i \textrm{ and } \partial D_j \cap \partial D_i \neq \emptyset\}$. Moreover, the boundaries of each subdomain $\partial D_i$ can be split into two parts: exterior boundaries $\partial_{\rm ex} D_i:=\partial D_i \cap \partial D$, and interface boundaries $\partial_{\rm in} D_i:=\cup_{j\in \linteri}\{\partial_{j} D_i\}$ where $\partial_{j} D_i:=\partial D_i \cap \partial D_j$, so that $\partial D_i=\partial_{\rm ex} D_i\cup\partial_{\rm in} D_{i}$. For simplicity, we label an interface $\partial_j D_i\in \partial_{\rm in}D_i$ with the index pair $(i,j\,)$. Grouping all interface indices associated with all subdomains $\{D_i\}_{i=1}^{M}$ together, we define $\linter:=\{(i,j\,)\,|\, i\in \{1,\dots,M\} \textrm{ and } j\in \linteri\}$, and assume that two different interfaces do not connect, i.e., $\pjoi\cap \partial_k D_i=\emptyset$ for $j\neq k$. 

We next assume that the input random vector $\xi$ can be decomposed into $\xi^{\top}=[\xi_1^{\top},\dots,\xi_M^{\top}]$, where $\xi_i$ is a vector collecting $\mathfrak{N}_i$ (a positive integer) random variables and is associated with the subdomain $D_i$. For each $\xi_i$, its support is denoted by $\Gamma_i$ and its local joint PDF is denoted by $\pi_{\xi_i}(\xi_i)$. The local outputs are denoted by $y_i\left(u(x,\xi)|_{D_i}\right)$, $i=1,\dots,M$, which depends only on the solution in one local subdomain and can be any quantity such as integral quantities and nodal point values.

For a given realization of the random vector $\xi$, \eqref{spdexi1}--\eqref{spdexi2} is a deterministic PDE problem. We herein review the non-overlapping domain decomposition method for deterministic PDEs, following the presentation in \cite{quavalbook,liao2015domain}. Let $g_i(x,\xi_i)$ denote boundary functions on exterior boundaries $\partial_{\rm ex}D_i$ for $i=1,\dots,M$, and let $\{g_{i,j}(x,\tau_{j,i})\}_{(i,j\,)\in\Lambda}$ denote interface functions on the interfaces $\{\pjoi\}_{(i,j\,)\in\Lambda}$. Each $g_{i,j}(x,\tau_{j,i})$ is dependent on a $\mN_{j,i}$-dimensional vector parameter $\tau_{j,i}$. For each subdomain $D_i$, we call $\tau_{j,i}$ an interface input parameter, and $\tau_{i,j}$ is called an interface output parameter. Grouping all the interface input parameters for $D_i$ together, we define $\tau_i\in \dsR^{\mN_i}$ with $\mN_i=\sum_{j\in \linteri}\mN_{j,i}$ as an entire interface input parameter, i.e.,
\begin{eqnarray}
	\tau_i:=\otimes_{j\in\linteri} \tau_{j,i},  \label{lam_in}
\end{eqnarray}
where $\otimes$ combines vectors as follows
\begin{eqnarray}
	\tau_{j_1,i}\otimes\tau_{j_2,i}:=\left\{\begin{array}{cl}{\lr \left[\tau^{\top}_{j_1,i}\, ,\tau^{\top}_{j_2,i}\right]^{\top}} & \textrm{if $j_1< j_2$}\\
		  {\lr \left[\tau^{\top}_{j_2,i}\, , \tau^{\top}_{j_1,i}\right]^{\top}} & \textrm{if $j_1> j_2$}\\
		  {\tau_{j_1,i}} & \textrm{if $j_1= j_2$}
		  \end{array} \right..
\end{eqnarray}

For each subdomain $D_i$ (for $i=1,\dots,M$), decomposed local operators and functions are defined as $\{\pL_i:=\pL|_{D_i}\}^M_{i=1}$, $\{\pB_i:=\pB|_{ D_i}\}^M_{i=1}$, $\{f_i:=f|_{D_i}\}^M_{i=1}$, $\{g_i:=g|_{D_i}\}^M_{i=1}$, which are the global operators and functions restricted to the subdomains. Given an initial guess $\tau_{i,j}^0$ (e.g., $\tau_{i,j}^0=0$) for each interface parameter, at domain decomposition iteration step $k$ ($k=0,1,\ldots$), we solve $M$ local problems: find $u(x,\xi_i,\tau^k_i): D_i\to \dsR$ for $i=1,\dots,M$, where $\tau^k_i=\otimes_{j\in \linteri} \tau^k_{j,i}$ is the interface input parameter on the $k$-th iteration, such that
\begin{eqnarray}
	\pL_i\left(x,\xi_i;u\left(x,\xi_i,\tau^k_i \right)\right)=f_i\left(x,\xi_i\right)\qquad & \textrm{in}& D_i,\label{dd1}\\
	\pB_i\left(x,\xi_i;u\left(x,\xi_i,\tau^k_i\right)\right)=g_i\left(x,\xi_i\right)\qquad  &\textrm{on}& \partial_{\rm ex} D_i,\label{dd2}\\
	\pB_{i,j}\left(x,\xi_i;u\left(x,\xi_i,\tau^k_i\right)\right)=g_{i,j}\left(x,\tau^k_{j,i}\right)\qquad  &\textrm{on}& \partial_j D_i, \textrm{ for all }
	j\in \linteri, \label{dd3}
\end{eqnarray}
where $\pB_{i,j}$ is an appropriate boundary operator posed on the interface $\pjoi$. 

For the next iteration, the interface data of each current local solution (e.g., the Dirichlet, the Neumann and the Robin type boundary condition values) need to be computed, and they are denoted by $\ph_{i,j}(u(x,\xi_i,\tau^k_i))$, where each $\ph_{i,j}$ is a functional to take the interface data of $u(x,\xi_i,\tau^k_i)$ (for $i=1,\dots,M$ and $j \in \Lambda_i$). 
For simplicity, we define 
\begin{eqnarray}
	h_{i,j}\left(\xi_i,\tau^k_i\right):=\ph_{i,j}\left(u\left(x,\xi_i,\tau^k_i\right)\right), \label{he}
\end{eqnarray}
where $h_{i,j}$ is referred to as a coupling function. After getting $h_{i,j}(\xi_i,\tau^k_i)$ at the $k$-th iteration step, the interface parameters are updated through
\begin{eqnarray}
	\tau^{k+1}_{i,j}=\theta_{i,j} h_{i,j}\left(\xi_i,\tau^k_i\right)+(1-\theta_{i,j})\tau^{k}_{i,j}, \label{ddupdate}
\end{eqnarray}
where $\theta_{i,j}$ is a non-negative acceleration parameter \cite{toselli2006domain}. 

In this work, the deterministic domain decomposition iteration is assumed to be convergent; that is, the following DD-convergence condition \cite{quavalbook} is satisfied. 
\begin{definition}
\label{dd_convergence}
{\textrm{(DD-convergence)}}. A domain decomposition method is convergent if
\begin{eqnarray}
\lim _{k \rightarrow \infty}\left\|u\left(x, \xi_i, \tau_i^k\right)-\left.u(x, \xi)\right|_{D_i}\right\|_{D_i} \rightarrow 0, \quad i=1, \ldots, M,
\end{eqnarray}
\end{definition}
where $\|v(x)\|_{D_i}:=\int_{D_i} v^2(x)\, {\rm{d}}x$ for any arbitrary function $v$ defined in $D_i$. 

The quantity $\tau^{\infty}_i:=\otimes_{j\in \linteri} \tau^{\infty}_{j,i}$ (where $\tau^{\infty}_{j,i}:=\lim_{k\to \infty} \tau^k_{j,i}$) is called a target interface input parameter for subdomain $D_i$, and $u(x,\xi_i,\tau^{\infty}_i)$ is a local stationary solution, which is consistent with the global solution of \eqref{spdexi1}--\eqref{spdexi2}, i.e., $u(x,\xi_i,\tau^{\infty}_i)=u(x,\xi)|_{D_i}$. For the stochastic problem, i.e., considering the inputs $\xi^{\top}=[\xi_1^{\top},\dots,\xi_M^{\top}]$ as a random vector, the joint PDF of $\xi_i$ and $\tau_i^{\infty}$ is called the target input PDF for the subdomain $D_i$ and is denoted by $\pi_{\xi_i,\tau_i}(\xi_i,\tau_i)$. 

\subsection{The domain-decomposed uncertainty quantification (DDUQ) approach} \label{ddqu-alg}
Our DDUQ approach proposed in \cite{liao2015domain} proceeds mainly through an offline stage and an online stage. In the offline stage, for each subdomain $D_i$, $i=1,\dots,M$, an initial guess for each (unknown) target input PDF $\pi_{\xi_i,\tau_i}(\xi_i,\tau_i)$ is constructed, which is called a proposal PDF and is denoted by 
\begin{eqnarray}
p_{\xi_i,\tau_i}(\xi_i,\tau_i)=\pi_{\xi_i}(\xi_i)p_{\tau_i}(\tau_i), \label{proposal_PDF}
\end{eqnarray}
where $\xi_i$ and $\tau_i$ are set to be independent in the proposal PDF and $p_{\tau_i}(\tau_i)$ is a proposal PDF for $\tau_i$. After that, Monte Carlo simulation with the proposal PDF is conducted for each subdomain $D_i$. That is, a large number $\Noff$ of samples $\{(\xi^{(s)}_i , \tau^{(s)}_i)\}^{\Noff}_{s=1}$ are drawn from  $p_{\xi_i,\tau_i}(\xi_i,\tau_i)$, and local solutions $\{u(x,\xi^{(s)}_i,\tau^{(s)}_i)\}^{\Noff}_{s=1}$ are obtained through solving the following local problem using simulators (e.g., finite element methods \cite{elman2014finite}), 
\begin{eqnarray}
	\pL_i\left(x,\xi^{(s)}_i;u\left(x,\xi^{(s)}_i,\tau^{(s)}_i \right)\right)=f_i\left(x,\xi^{(s)}_i\right)
	\qquad & \textrm{in}& D_i,\label{lp1}\\
	\pB_i\left(x,\xi^{(s)}_i;u\left(x,\xi^{(s)}_i,\tau^{(s)}_i\right)\right)=g_i\left(x,\xi^{(s)}_i\right)
	\qquad  &\textrm{on}& \partial_{\rm ex} D_i,\label{lp2}\\
	\pB_{i,j}\left(x,\xi^{(s)}_i;u\left(x,\xi^{(s)}_i,\tau^{(s)}_i\right)\right)=g_{i,j}\left(x,\tau^{(s)}_{j,i}\right)
	\qquad  &\textrm{on}& \partial_j D_i, \textrm{ for all } {\lr j\in \linteri},\label{lp3}
\end{eqnarray}
where $\tau^{(s)}_i=\otimes_{j\in \linteri} \tau^{(s)}_{j,i}$. With the local solution samples $\{u(x,\xi^{(s)}_i,\tau^{(s)}_i)\}^{\Noff}_{s=1}$, the outputs of interest $\{y_i(u(x,\xi^{(s)}_i,\tau^{(s)}_i))\}^{\Noff}_{s=1}$ are  evaluated, and surrogates for the coupling functions $\{h_{i,j}\}_{(i,j)\in\Lambda}$ (introduced in Section \ref{dd}) are constructed. These surrogates for $\{h_{i,j}\}_{(i,j)\in\Lambda}$ are denoted by $\{\tilde{h}_{i,j}\}_{(i,j)\in\Lambda}$ and are called coupling surrogates in the following.

In the online stage, a large number $\Non$ of samples $\{\xi^{(s)}\}_{s=1}^{\Non}$ are drawn from the joint distribution of the system input $\pi_{\xi}(\xi)$, and the corresponding target interface parameters $\{\tau_{i,j}^{\infty}(\xi^{(s)})\}_{s=1}^{\Non}$ are computed through domain decomposition iteration using the coupling surrogates $\tilde{h}_{i,j}$ for $(i,j)\in \Lambda$. With the samples $\{(\xi^{(s)}_i,\tau^{\infty}_i(\xi^{(s)}))\}^{\Non}_{s=1}$, one can estimate each target input PDF $\pi_{\xi_i,\tau_i}(\xi_i,\tau_i)$ (introduced in Section \ref{dd}) through density estimation techniques, and the estimated target input PDF is denoted by $\hat\pi_{\xi_i,\tau_i}(\xi_i,\tau_i)$. At the end of the online stage, the offline samples $\{y_i(u(x,\xi^{(s)}_i,\tau^{(s)}_i))\}^{\Noff}_{s=1}$ are re-weighted with the following weights so as to satisfy the domain decomposition coupling conditions,
    \begin{eqnarray} w^{(s)}_i=\frac{\hat{\pi}_{\xi_i,\tau_i}\left(\xi^{(s)}_i,\tau^{(s)}_i\right)}{p_{\xi_i,\tau_i}\left(\xi^{(s)}_i,\tau^{(s)}_i\right)},\quad s=1,\dots,\Noff,\quad i=1,\dots,M.
    \label{dduq_weights}
\end{eqnarray}

The density estimation step is a major challenge in DDUQ, and density estimation for high-dimensional random vectors in general is an open problem in machine learning. Let $\pi_{\tau_i\mid\xi_i}(\tau_i\mid\xi_i)$ denote the target conditional PDF of $\tau_i$ given $\xi_i$ for subdomain $D_i$, i.e., the conditional PDF of $\tau_i^{\infty}$ (introduced in Section \ref{dd}) given $\xi_i$,
and let $\hat{\pi}_{\tau_i\mid\xi_i}(\tau_i\mid\xi_i)$ denote an approximation of $\pi_{\tau_i\mid\xi_i}(\tau_i\mid\xi_i)$ (details of the conditional density estimation is discussed in the next section). Then, the estimated target input PDF can be expressed as $\hat{\pi}_{\xi_i,\tau_i}(\xi_i,\tau_i)=\hat{\pi}_{\tau_i \mid \xi_i}(\tau_i\mid\xi_i)\pi_{\xi_i}(\xi_i)$, and \eqref{dduq_weights} can be reformulated as 
\begin{eqnarray}
	w_i^{(s)}
    &=&\frac{\hat{\pi}_{\tau_i\mid\xi_i}\left(\tau_i^{(s)}\mid\xi_i^{(s)}\right)\pi_{\xi_i}\left(\xi_i^{(s)}\right)}{\pi_{\xi_i}\left(\xi_i^{(s)}\right)p_{\tau_i}\left(\tau_i^{(s)}\right)}\nonumber\\
    &=&\frac{\hat{\pi}_{\tau_i\mid\xi_i}\left(\tau_i^{(s)}\mid\xi_i^{(s)}\right)}{p_{\tau_i}\left(\tau_i^{(s)}\right)},\quad s=1,\ldots,\Noff, \quad i=1,\dots,M.
    \label{conditional_weights}
\end{eqnarray}
From \eqref{conditional_weights}, the original density estimation problem for the joint distribution $\pi_{\xi_i,\tau_i}(\xi_i,\tau_i)$ is replaced by density estimation for the conditional distribution $\pi_{\tau_i\mid\xi_i}(\tau_i\mid\xi_i)$. For this purpose, a new flow-based model is proposed to estimate conditional PDFs in the next section.

%% file: sections/section3.tex
Let $(\ba,\bc)$ with $\ba\in \dsR^{|\ba|}$ and $\bc\in \dsR^{|\bc|}$ be a pair of random variables associated with a given data set, and $p(\ba\mid\bc)$ denote the conditional probability density function of $\ba$ given $\bc$, where $|\boldsymbol{v}|$ denotes the dimensionality of any arbitrary vector $\boldsymbol{v}$. The goal of conditional flow models is to find an approximation of $p(\ba\mid\bc)$ from available data. Let $\bz \in \dsR^{|\ba|}$ represent a random vector associated with a PDF $p_{\bz}(\bz)$, where $p_{\bz}(\bz)$ is a prior distribution (e.g., Gaussian distributions). {\color{black} The conditional flow-based generative modeling is to seek a distinct invertible mapping $\bz=f_{\bc}(\ba)$ for each realization of the conditional variable $\bc$. By the change of variables, the conditional PDF of $\ba = f_{\bc}^{-1}(\bz)$ given $\bc$ can be computed through
\begin{eqnarray}
p_{\ba\mid\bc}(\ba\mid\bc)=p_{\bz}(f_{\bc}(\ba))\left|\det \nabla_{\ba}f_{\bc}(\ba)\right|.\label{chang_of_var} 
\end{eqnarray}
In flow-based generative models, the invertible mapping $f_{\bc}(\cdot)$ is constructed by stacking a sequence of simple invertible layers. Each of these layers is a shallow neural network, and thus the overall mapping is a deep neural network. The mapping $f_{\bc}(\cdot)$ can be expressed in a composite form:
\begin{eqnarray}
    \qquad \quad \bz=f_{\bc}(\ba)=f_{\bc,[R]}\circ \cdots \circ f_{\bc,[1]}(\ba) \quad \text{and} \quad \ba=f_{\bc}^{-1}(\bz)=f^{-1}_{\bc,[R]}\circ \cdots \circ f^{-1}_{\bc,[1]}(\bz),\label{inver_mapping}
\end{eqnarray}
where $f_{\bc,[r]}(\cdot)$ (for $r=1,\dots,R$) is the bijection given $\bc$ at stage $r$. In this way, the determinant of the Jacobian in equation \eqref{chang_of_var} can be computed by the chain rule:
\begin{eqnarray}
    \left|\det \nabla_{\ba}f_{\bc}(\ba)\right|=\prod_{r=1}^{R}\left|\det  \nabla_{\ba_{[r-1]}}f_{\bc,[r]}(\ba)\right|,\label{total_Jacobian}
\end{eqnarray}
where $\ba_{[r-1]}$ is the input of $f_{\bc,[r]}(\cdot)$. Once the prior distribution $p_{\bz}(\bz)$ and each bijection $f_{\bc,[r]}(\cdot)$ are specified, equations \eqref{chang_of_var} and \eqref{total_Jacobian} provide an explicit conditional density model of $\ba$ given $\bc$. In the rest of this section, the construction of the mapping $f_{\bc}(\cdot)$ in \eqref{inver_mapping} is presented, which is modified from our flow model KRnet \cite{tang2020deep,tang2022adaptive}, and this modified  KRnet is referred to as the conditional KRnet (cKRnet) in the following.}

\subsection{Conditional affine coupling layers} \label{sec:affine} 
The design of each bijection $f_{\bc,[r]}(\cdot)$ in \eqref{inver_mapping} necessitates careful consideration to ensure a well-conditioned inverse and a tractable Jacobian determinant. In real-valued non-volume preserving (real NVP) transformations \cite{dinh2016density}, a general efficient architecture to construct bijections is proposed, which is called the affine coupling layers. In our recent work \cite{tang2020deep,tang2022adaptive}, modifications for the affine coupling layers are developed, of which the determinant of the corresponding Jacobian is bounded. For the specific setting in this work which focuses on conditional distributions, related works include conditional affine coupling layers for point cloud generation \cite{pumarola2020c} and image super-resolution \cite{lugmayr2020srflow, liang2021hierarchical}, temporal normalizing flows \cite{both2019temporal,feng2022solving}, and a systematic extension of KRnet for time dependent problems \cite{he2024adaptive}. 

In this work, conditional affine coupling layers are proposed as follows. Let $\ba$ be the input of a conditional affine coupling layer and $[\ba_{1}, \ba_{2}]^{\top}$ be a partition of $\ba$ with $\ba_{1} \in \dsR^{|\ba_1|}$, $\ba_{2}\in \dsR^{|\ba_2|}$ and $|\ba|=|\ba_1|+|\ba_2|$. Given the conditional variable $\bc$, the output of the conditional affine coupling layer $\hat{\ba}=[\hat{\ba}_{1}, \hat{\ba}_{2}] ^{\top}$ is defined as
\begin{eqnarray}
    \begin{aligned}
        \hat{\ba}_{1}&=\ba_{1},\\ 
        \hat{\ba}_{2}&=\ba_{2} \odot\left(1+\gamma \tanh \left(\bs_{\bc}\left(\ba_{1}\right)\right)\right)+e^{\boldsymbol{\beta}} \odot \tanh \left(\bt_{\bc}\left(\ba_{1}\right)\right),
    \end{aligned}\label{conditional_affine_coupling_layer}
\end{eqnarray}
where $0< \gamma <1$ is a hyperparameter, $\boldsymbol{\beta}$ is a trainable parameter, and $\odot$ represents the Hadamard product (or the element-wise product). In \eqref{conditional_affine_coupling_layer}, different from our original KRnet, $\bs_{\bc}$ and $\bt_{\bc}$ are two functions mapping $\dsR^{|\ba_1|+|\bc|}$ to $\dsR^{|\ba_2|}$ rather than $\dsR^{|\ba_1|}$ to $\dsR^{|\ba_2|}$, which are used to generate different scaling and translation for different realization of $\bc$ respectively. In practice, $\bs_{\bc}$ and $\bt_{\bc}$ are implemented through neural networks typically. The Jacobian of $\hat{\ba}$ with respect to $\ba$ given $\bc$ is
\begin{eqnarray*}
    \nabla_{\ba}\hat{\ba}=\left[\begin{array}{rl}
    \mathbb{I} & \mathbf{0} \\
    \frac{\partial \hat{\ba}_{2}}{\partial \ba_{1}} & \operatorname{diag}\left(\mathbf{1}+ \gamma \tanh \left(\bs_{\bc}\left(\ba_{1}\right)\right)\right)
    \end{array}\right].
\end{eqnarray*}
It can be seen that the Jacobian of $\hat{\ba}$ is a lower triangular matrix, and its determinant can be easily computed. 

According to the above definition of conditional affine coupling layers \eqref{conditional_affine_coupling_layer}, only $\hat{\ba}_{2}$ is updated while $\hat{\ba}_{1}$ remains unchanged. This issue can be addressed by composing conditional affine coupling layers in an alternating pattern. Specifically, the next conditional affine coupling layer is defined as
\begin{eqnarray}
    \begin{aligned}
        \hat{\ba}_{1}&=\ba_{1} \odot\left(1+\gamma \tanh \left(\bs_{\bc}\left(\ba_{2}\right)\right)\right)+e^{\boldsymbol{\beta}} \odot \tanh \left(\bt_{\bc}\left(\ba_{2}\right)\right),\\
        \hat{\ba}_{2}&=\ba_{2}.    \end{aligned}\label{conditional_affine_coupling_layer2}
\end{eqnarray}
Then, \eqref{conditional_affine_coupling_layer} and \eqref{conditional_affine_coupling_layer2} alternate in each bijection $f_{\bc,[r]}(\cdot)$. 

\subsection{Conditional KRnet (cKRnet) structure}
The basic idea of cKRnet is to define the structure of $f_{\bc}(\ba)$ (in \eqref{inver_mapping}) in terms of the Knothe-Rosenblatt rearrangement (Cf. \cite{carlier2010knothe, tang2020deep} for more details). Specifically, let $\ba=[\ba^{(1)},\ldots,\ba^{(R)}]^{\top}$ be a partition of $\ba$, where $\ba^{(r)}=[\ba^{(r)}_{1},\ldots,\ba^{(r)}_{m}]^{\top}$ with $2\leq R \leq |\ba|$, $1\leq m < |\ba|$, and $\sum_{r=1}^{R} |\ba^{(r)}|=|\ba|$. The overall structure of cKRnet is defined as a triangular structure:
\begin{eqnarray}
    \bz=f_{\mathsf{cKR}}\left(\ba,\bc\right)=\left[\begin{array}{l}
    f_{\bc,[1]}\left(\ba^{(1)}, \ldots, \ba^{(R)}\right) \\
    \vdots \\
    f_{\bc,[R-1]}\left(\ba^{(1)}, \ba^{(2)}\right) \\
    f_{[R]}\left(\ba^{(1)}\right)
    \end{array}\right],\label{cKR_form}
\end{eqnarray}
where each $f_{\bc,[r]}(\cdot)$ (for $r=1,\ldots,R-1$) is an invertible mapping depending on $\bc$ and the last mapping $f_{[R]}(\cdot)$ is a nonlinear invertible layer $L_{N}$ (defined in \cite{tang2022adaptive}) independent of $\bc$ and the detailed workflow of cKRnet can be represented as 
\begin{eqnarray}
    \bz=f_{\mathsf{cKR}}(\ba, \bc)=L_{N}\circ f_{\bc,[R-1]}\circ \cdots \circ f_{\bc,[1]}(\ba),\label{outer_loop}
\end{eqnarray}
where $f_{\bc,[r]}$ (for $r=1,\ldots,R-1$) is defined as
\begin{eqnarray}
    f_{\bc,[r]}=L_{S} \circ f_{\bc,[r,L]} \circ \cdots \circ f_{\bc,[r,1]}(\ba).\label{inner_loop}
\end{eqnarray}
Here, $f_{\bc,[r,l]}(\cdot)$ is a combination of a conditional affine coupling \eqref{conditional_affine_coupling_layer} (or \eqref{conditional_affine_coupling_layer2}) layer and a scale-bias layer for $l=1,\dots,L$, and $L_{S}$ indicates the squeezing layer (the scale-bias layer and the squeezing layer are defined in \cite{tang2020deep,tang2022adaptive}). 

\subsection{cKRnet for conditional density estimation}\label{cKR_for_cds}  
Let $p_{\ba,\bc}(\ba,\bc)$ be a joint probability density function (PDF) associated with a pair of random variables $(\ba,\bc)$ and $p_{\bc}(\bc)$ is the marginal PDF associated with $\bc$. The conditional PDF of $\ba$ given $\bc$ is denoted as $p_{\ba\mid\bc}(\ba\mid\bc)=p_{\ba,\bc}(\ba,\bc)/p_{\bc}(\bc)$. Given a set of samples $\mathcal{S}=\{(\ba^{(s)},\bc^{(s)})\}_{s=1}^{ N_t}$ drawn from $p_{\ba,\bc}(\ba,\bc)$, the goal of conditional density estimation is to find an approximation of $p_{\ba\mid\bc}(\ba\mid\bc)$. With a prior distribution $p_{\bz}(\bz)$, using our cKRnet $f_{\mathsf{cKR}}(\ba,\bc)$ as the mapping $f_{\bc}(\ba)$ in \eqref{chang_of_var} gives an approximation of $p_{\ba\mid\bc}(\ba\mid\bc)$,  i.e.,
\begin{eqnarray}
    p_{\ba\mid\bc}(\ba\mid\bc)\approx p_{\mathsf{cKR}}\left(\ba\mid\bc;\Theta\right):=p_{\bz}\left(f_{\mathsf{cKR}}\left(\ba,\bc;\Theta\right)\right)|\det \nabla_{\ba}f_{\mathsf{cKR}}\left(\ba,\bc;\Theta\right)|,\label{estimated_PDF}
\end{eqnarray}
where $\Theta$ represents the parameters of the cKRnet that needs to be trained. To determine the optimal parameters, we minimize the Kullback-Leibler (KL) divergence between $p_{\ba\mid\bc}(\ba\mid\bc)$ and $p_{\mathsf{cKR}}(\ba\mid\bc;\Theta)$, which is defined as
\begin{eqnarray*}
    \begin{aligned}
    \mathcal{D}_{K L}(p_{\ba\mid\bc}(\ba\mid\bc) \| p_{\mathsf{cKR}}(\ba\mid\bc;\Theta))
    &=\mathbb{E}_{\bc \sim p_{\bc}(\bc)}\left[\mathbb{E}_{\ba \sim p_{\ba\mid\bc}(\ba\mid\bc)}  \left[\log\frac{p_{\ba\mid\bc}(\ba\mid\bc)}{p_{\mathsf{cKR}}(\ba\mid\bc;\Theta)} \right] \right]\\
    &=\mathbb{E}_{(\ba,\bc) \sim p_{(\ba,\bc)}(\ba,\bc)}\left[\log p_{\ba\mid\bc}(\ba\mid\bc)-\log p_{\mathsf{cKR}}(\ba\mid\bc;\Theta)\right]\\
    &=\mathbb{E}_{(\ba,\bc) \sim p_{(\ba,\bc)}(\ba,\bc)}\log p_{\ba\mid\bc}(\ba\mid\bc)-\mathbb{E}_{(\ba,\bc) \sim p_{(\ba,\bc)}(\ba,\bc)}\log p_{\mathsf{cKR}}(\ba\mid\bc;\Theta). \label{KL_divergence}
    \end{aligned}
\end{eqnarray*}
Since the first expectation above is independent of $\Theta$, minimizing the KL divergence is equivalent to minimizing the following loss functional
\begin{eqnarray}
\mathcal{L}\left(p_{\mathsf{cKR}}\left(\ba\mid\bc;\Theta\right)\right)= - \mathbb{E}_{\left(\ba,\bc\right) \sim p_{\left(\ba,\bc\right)}\left(\ba,\bc\right)}\log p_{\mathsf{cKR}}\left(\ba\mid\bc;\Theta\right).
\label{loss_function1}
\end{eqnarray}
Then we use the given data set $\mathcal{S}$ to approximate the loss functional, i.e.,
\begin{eqnarray}
\begin{aligned}
\hat{\mathcal{L}}\left(p_{\mathsf{cKR}}\left(\ba\mid\bc;\Theta\right)\right)&= - \frac{1}{N_t}\sum_{s=1}^{N_t}\log p_{\mathsf{cKR}}\left(\ba^{(s)}\mid\bc^{(s)};\Theta\right),
\label{loss_function2}
\end{aligned}
\end{eqnarray}
based on which the optimal parameter $\Theta^{*}$ is chosen as:
\begin{eqnarray}
\Theta^{*}=\arg \min_{\Theta}\hat{\mathcal{L}}\left(p_{\mathsf{cKR}}\left(\ba\mid\bc;\Theta\right)\right),
\label{optimal_parameter}
\end{eqnarray}
which can be solved through the stochastic gradient descent methods \cite{bottou2018optimization}.

\subsection{Performance of cKRnet}\label{cKR_per} 
Here, we compare the performance of cKRnet and the standard KRnet (introduced in \cite{tang2020deep,tang2022adaptive}) for density estimation, and the following $16$-dimensional Gaussian mixture distribution is considered,
\begin{equation}
    p_{\mathsf{ref}}(\boldsymbol{\omega})=\beta_1p_1(\boldsymbol{\omega})+\beta_2p_2(\boldsymbol{\omega})+\beta_3p_3(\boldsymbol{\omega}),\label{ref_pdf}
\end{equation}
where $\boldsymbol{\omega}\in \mR^{16}$ and each $p_m(\boldsymbol{\omega})$ (for $m=1,2,3$) is  a probability density function of the normal distribution with mean $\mu_m \in \mR^{16}$ and covariance $\Sigma_m \in \mR^{16\times 16}$. The mean vectors $\mu_m$ and the parameters $\beta_m$ (for $m=1,2,3$) are set as 
\begin{eqnarray*}
\begin{aligned}
    &\mu_1=[-1,-1,-0.3,-0.3,-0.4,-0.4,-1.6,-1.6,-0.8,-0.8,-0.5,-0.5,-0.5,-0.3,-1,-1]^{\top},\\
    &\mu_3=[32,32,32.6,32.6,32.8,32.8,33.2,33.2,32.4,32.4,30.8,30.8,31,31,31.6,31.6]^{\top},\\
    &\mu_2=\frac{\mu_1+\mu_3}{2},\,\beta_1=0.3,\,\beta_2=0.4,\,\beta_3=0.3.
\end{aligned}
\end{eqnarray*}
The covariance matrices $\Sigma_m$ (for $m=1,2,3$) are randomly constructed. Specifically, three matrices $\tilde{\Sigma}_m \in \mR^{16\times 16}$ (for $m=1,2,3$) are constructed, of which the entries are generated through the uniform distribution with range [0,1], and then each covariance matrix is obtained as $\Sigma_m=\tilde{\Sigma}_m\tilde{\Sigma}_m^{\top}$ for $m=1,2,3$.

In this test, we split the random variable $\boldsymbol{\omega}$ into two parts ($\ba$ and $\bc$), i.e., $\boldsymbol{\omega}^{\top}=[\ba^{\top},\bc^{\top}]$. Denoting the dimensionality of $\bc$ by $|\bc|$, three cases of the decomposition are considered, which are associated with $|\bc|=8$, $|\bc|=12$ and $|\bc|=14$ respectively. Let $\mathcal{S}=\{\boldsymbol{\omega}^{(s)}\}_{s=1}^{ N_t}=\{(\ba^{(s)},\bc^{(s)})\}_{s=1}^{ N_t}$ denote a training dataset with size $N_t$, which is generated through drawing samples from $p_{\mathsf{ref}}(\boldsymbol{\omega})$, and our test here is to estimate the PDF using the training dataset. Different sizes of the training dataset are considered in the following, which are $N_t=2\times 10^5$, $N_t=4\times 10^5$ and $N_t=6\times 10^5$. Then, KRnet and cKRnet are applied to estimate $p_{\mathsf{ref}}(\bw)$. For KRnet (details are shown in \cite{tang2020deep}), the estimated PDF using the data set $\mathcal{S}$ is denoted by $p_{\mathsf{KR}}(\bw;\Theta)$.
For cKRnet, the marginal PDF $p_{\bc}(\bc)$ is considered to be known and the conditional PDF $p_{\ba\mid\bc}(\ba\mid\bc)$ is estimated (see Section \ref{cKR_for_cds}). Denoting $p_{\mathsf{cKR}}(\ba\mid\bc;\Theta)$ as an approximation of $p_{\ba\mid\bc}(\ba\mid\bc)$, the joint PDF estimated using cKRnet is obtained as $p_{\mathsf{cKR}}(\bw;\Theta)=p_{\mathsf{cKR}}(\ba\mid\bc;\Theta)p(\bc)$.

\begin{figure}[!ht]
    \centerline{
    \begin{tabular}{ccc}
    \includegraphics[width=0.32\textwidth]{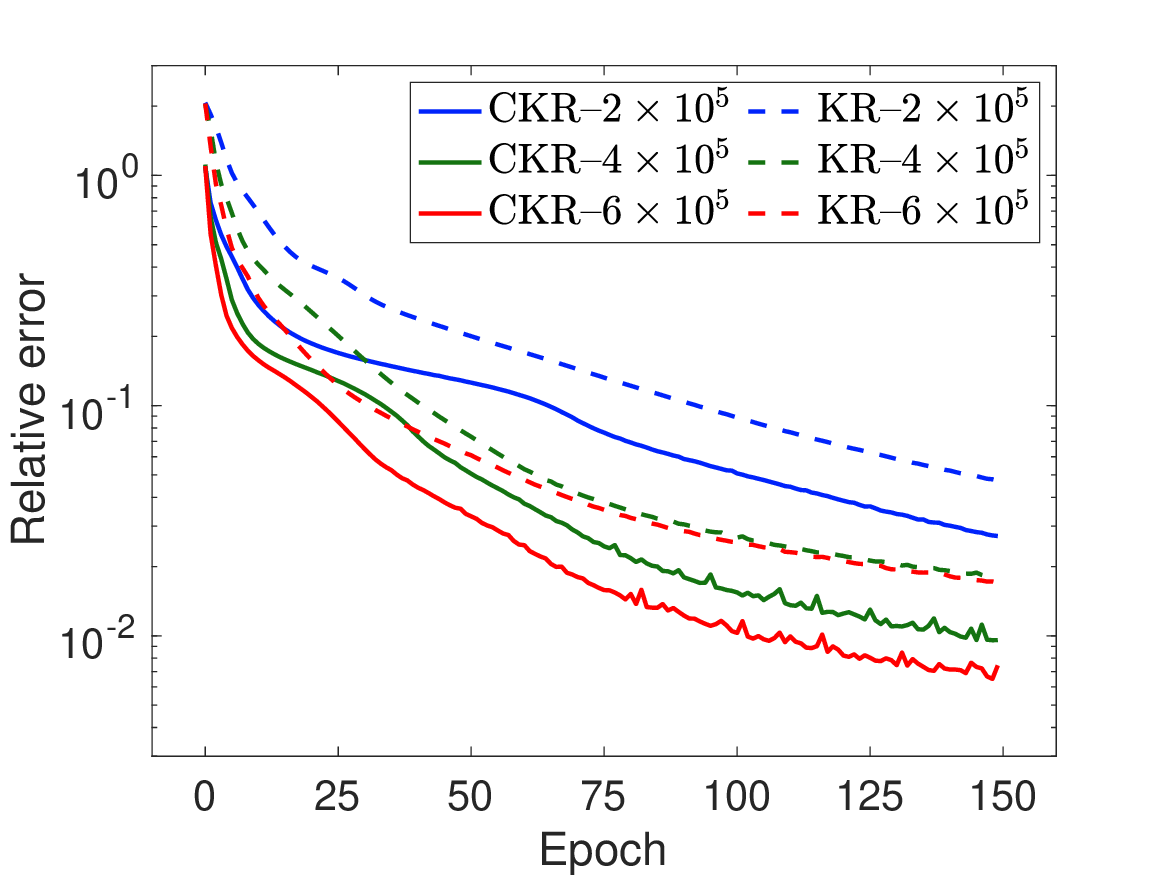}&
    \includegraphics[width=0.32\textwidth]{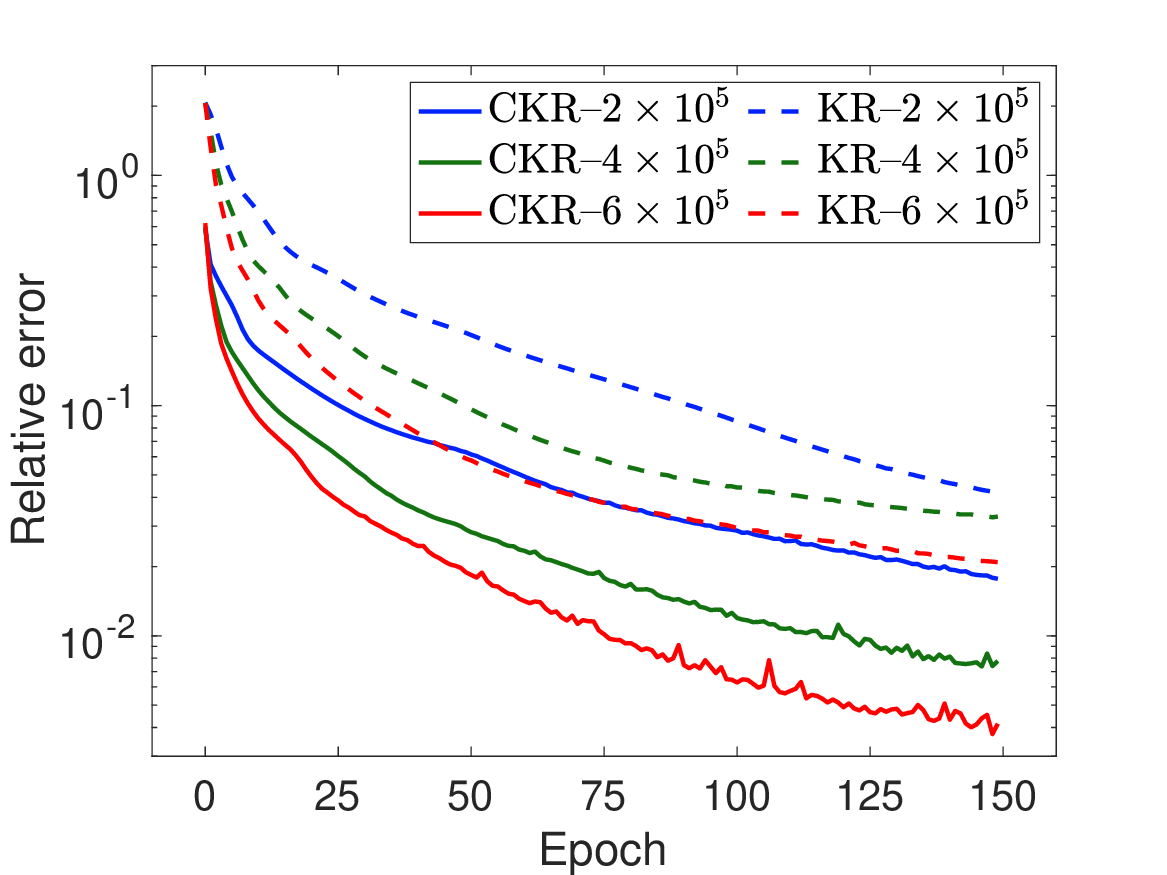}&
    \includegraphics[width=0.32\textwidth]{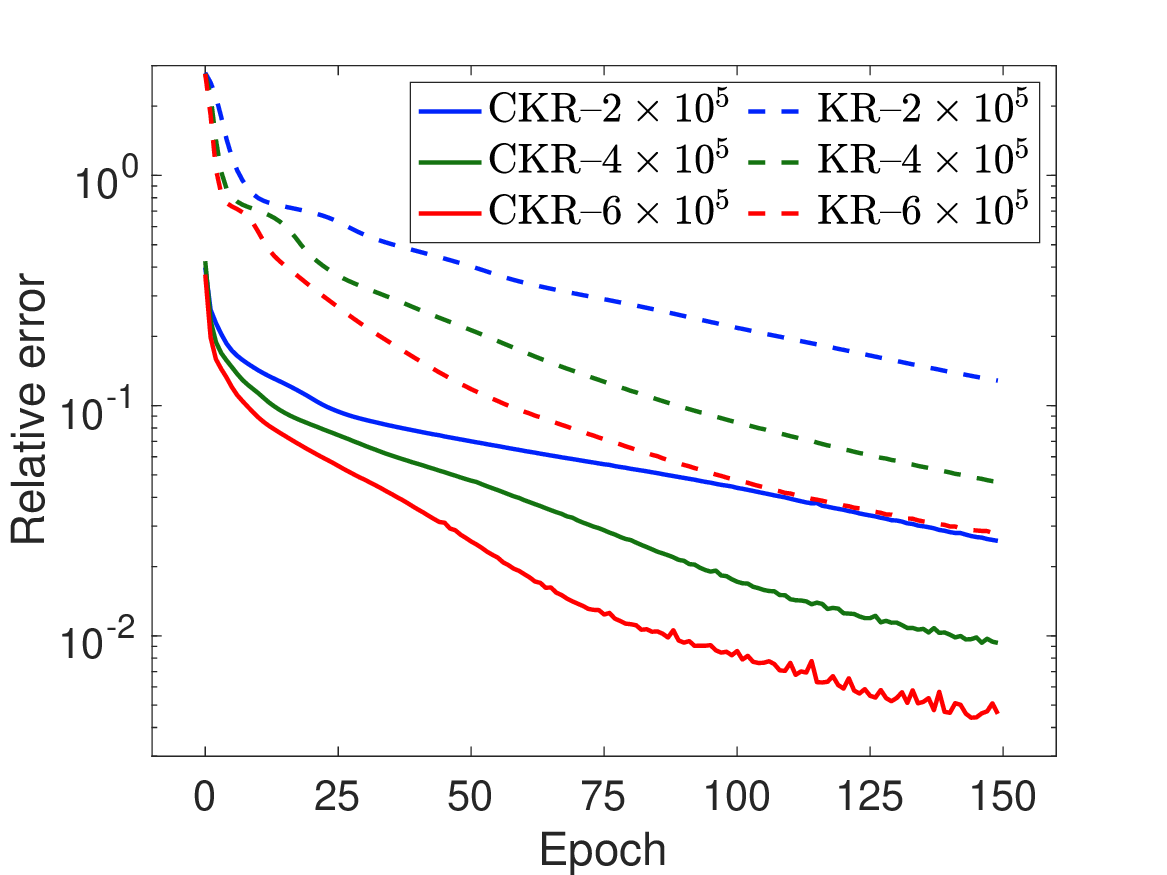}\\
    (a) $|\bc|=8$&
    (b) $|\bc|=12$& 
    (c) $|\bc|=14$\\
    \end{tabular}}
    \caption{\lr
    Relative errors of cKRnet (CKR) and KRnet (KR) for the Gaussian mixture distribution test problem.}
    \label{fig_cKR_vs_KR}
\end{figure}

For each case of $|\bc|$, cKRnet contain $R=16-|\bc|$ invertible mappings (i.e., $f_{\bc,[r]}$ defined in \eqref{inner_loop} or the nonlinear invertible layer $L_N$), and each invertible mapping except $L_N$ is configured with $L=4$ conditional affine coupling layers. In each conditional affine coupling layer, $\bs_{\bc}$ and $\bt_{\bc}$ are implemented through a fully connected network containing two hidden layers with thirty two neurons in each layer. Additionally, the rectified linear unit function (ReLU) function is employed as the activation function. For fair comparison, the setting of KRnet is the same as that for cKRnet as above. Then, the Adam optimizer \cite{kingma2014adam} with a learning rate of 0.001 is used to train cKRnet and KRnet. To assess the accuracy of estimated PDFs, the following relative error for $p_{\mathsf{cKR}}(\bw;\Theta)$ is defined 
\begin{eqnarray}
\delta(p_{\mathsf{ref}}(\bw),p_{\mathsf{cKR}}(\bw;\Theta))=\frac{\mathcal{D}_{K L}(p_{\mathsf{ref}}(\bw)\|p_{\mathsf{cKR}}(\bw;\Theta))}{H(p_{\mathsf{ref}}(\bw))}, \label{relative_error}
\end{eqnarray}
where $H(p_{\mathsf{ref}}(\bw))$ is the entropy of $p_{\mathsf{ref}}(\bw)$; similarly, the error $\delta(p_{\mathsf{ref}}(\bw),p_{\mathsf{KR}}(\bw;\Theta))$ is defined through replacing $p_{\mathsf{cKR}}(\bw;\Theta)$ in \eqref{relative_error} by $p_{\mathsf{KR}}(\bw;\Theta)$. The KL divergence and the entropy in \eqref{relative_error} are approximated using the Monte Carlo method with a large validation dataset generated from $p_{\mathsf{ref}}(\bw)$, and the size of the validation dataset is set to $10^7$ in this test. During the training process, the relative errors $\delta(p_{\mathsf{ref}}(\bw),p_{\mathsf{cKR}}(\bw;\Theta))$ and $\delta(p_{\mathsf{ref}}(\bw),p_{\mathsf{KR}}(\bw;\Theta))$ are computed for each epoch on the validation dataset. Figure \ref{fig_cKR_vs_KR} shows the relative errors of cKRnet and KRnet, considering varying sample sizes and values of $|\bc|$. It can be seen that the relative errors of cKRnet are smaller than those of KRnet for all the cases. The relative errors decrease as the number of training samples increases from $2\times 10^5$ to $6\times 10^5$ for both cKRnet and KRnet. Furthermore, the larger the value of $|\bc|$, the greater the discrepancies between the relative errors of cKRnet and the relative errors of KRnet.

%% file: sections/section4.tex
Our goal is to perform the propagation of uncertainty from the inputs $\xi$ to the outputs of interest $\{y_i(u(x,\xi)|_{D_i})\}_{i=1}^M$. In the following, Section~\ref{main_alg} presents the details of our proposed CKR-DDUQ algorithm and Section~\ref{analysis} gives its  convergence analysis.

\subsection{CKR-DDUQ Algorithm}\label{main_alg}
In this section, our cKRnet based domain decomposed uncertainty quantification (CKR-DDUQ) approach is presented, which involves deep neural networks (DNNs) for coupling surrogates modeling and cKRnet for density estimation. The CKR-DDUQ approach is composed of two stages: the offline stage and the online stage. In the offline stage, the Monte Carlo simulation is conducted for each subdomain to generate samples, and based on the samples, the surrogates of coupling functions are constructed using DNNs. The online stage conducts the domain decomposition iteration with the coupling surrogates and trains cKRnets to re-weight the offline samples. The offline stage involves local PDE solves which can be expensive; no PDE solve is required in the online stage, which is then considered to be cheap. 

In the offline stage, the first step is to specify a proposal PDF $p_{\tau_i}(\tau_i)$ for each interface parameter $\tau_i$ (for $i=1,\ldots,M$), and the proposal input PDF for each subdomain $D_i$ can be set as $p_{\xi_i,\tau_i}(\xi_i,\tau_i)=\pi_{\xi_i}(\xi_i)p_{\tau_i}(\tau_i)$, where $\pi_{\xi_i}(\xi_i)$ is a given PDF for the local system input parameter. In the proposal, $\xi_i$ and $\tau_i$ are set to be independent  for simplicity, and the proposal PDF $p_{\tau_i}(\tau_i)$ can be chosen to be any PDF as long as its support is large enough to cover the expected support of the target interface parameter $\tau_i^{\infty}(\xi)$ (introduced in Section \ref{dd}). After that, the Monte Carlo simulation is conducted for each local problem. That is, a large number $\Noff$ of samples $\{(\xi_i^{(s)},\tau_i^{(s)})\}_{s=1}^{\Noff}$  are drawn from the proposal input PDF $p_{\xi_i,\tau_i}(\xi_i,\tau_i)$ for each subdomain, where the superscript ($s$) denotes the $s$-th sample,  and the local solutions $\{u(x,\xi_i^{(s)},\tau_i^{(s)})\}_{s=1}^{\Noff}$ are computed by solving the local deterministic problems \eqref{lp1}--\eqref{lp3} for each input sample. 
Once the local solutions are computed, one can evaluate the local outputs of interest $\{y_i(u(x,\xi_i^{(s)},\tau_i^{(s)}))\}_{s=1}^{\Noff}$ and  the values of the coupling functions $\{h_{i,j}\}_{(i,j)\in\linter}$. Then, with the samples $\{((\xi_i^{(s)},\tau_{i}^{(s)}), h_{i,j}(\xi_i^{(s)},\tau_{i}^{(s)}))\}_{s=1}^{\Noff}$ of the coupling functions for $(i,j\,)\in \linter$, the coupling surrogates $\{\tilde{h}_{i,j}\}_{(i,j)\in\linter}$ can be constructed. Here, we focus on the deep neural network approximations for surrogate modeling \cite{rumelhart1986learning,he2016deep,zhu2019physics,wang2024deep}, as the size of the training data set (samples of the coupling functions) is typically large. 

In the online stage, the first step is to generate a large number $\Non$ of samples $\{\xi^{(s)}\}_{s=1}^{\Non}$ from the joint PDF of input parameters $\pi_{\xi}(\xi)$. Then the domain decomposition iterations for each sample $\xi^{(s)}$ are conducted using the coupling surrogates $\{\tilde{h}_{i,j}\}_{(i,j)\in\linter}$ to evaluate the corresponding target interface parameters $\tau_{i,j}^{\infty}(\xi^{(s)})$. As the domain decomposition iteration procedure here acts only with the coupling surrogates $\{\tilde{h}_{i,j}\}_{(i,j)\in\linter}$, no PDE solve is required. 
To determine the convergence of the domain decomposition iteration, the following error indicator is defined
\begin{eqnarray}
    \epsilon_{k}=\max_{(i,j\,)\in \linter}\left(\left\|\tau^{k}_{i,j}\left(\xi^{(s)}\right)-\tau^{k-1}_{i,j}\left(\xi^{(s)}\right)\right\|_{\infty}\right), \quad k=1,2,\dots,
    \label{error_indicator}
\end{eqnarray}
where $\|\cdot\|_{\infty}$ is the vector infinity norm and $k$ is the domain decomposition iteration step. When this error indicator is smaller than a given tolerance, the iteration terminates and we take $\tau^{\infty}_{i,j}(\xi^{(s)})=\tau^{k}_{i,j}(\xi^{(s)})$ {\lr for $(i,j\,)\in \linter$} and $\tau^{\infty}_i(\xi^{(s)})=\tau^{k}_i(\xi^{(s)})$ {\lr for $i=1,\dots,M$}. After that, the local target samples $\{(\xi_i^{(s)},\tau_i^{\infty}(\xi^{(s)}))\}_{s=1}^{\Non}$ for $i=1,\dots,M$ are obtained. With local target samples, each target conditional PDF $\pi_{\tau_i\mid\xi_i}(\tau_i\mid\xi_i)$ is approximated using our cKRnet (i.e., \eqref{estimated_PDF} and the procedures presented in Section \ref{section_ckrnet}), and the approximation (the estimated PDF) is denoted by $\hat{\pi}_{\tau_i\mid\xi_i}(\tau_i\mid\xi_i)$. The final step of the online stage is to re-weight the offline output samples $\{y_i(u(x,\xi_i^{(s)},\tau_i^{(s)}))\}_{s=1}^{\Noff}$ by the weights $\{w_i^{(s)}\}_{s=1}^{\Noff}$ defined in \eqref{conditional_weights}. 

The details of the online stage are summarized in Algorithm \ref{alg_on1}. Here, $\Non$ is the number of samples and $tol$ is a given tolerance for the domain decomposition iteration. In the density estimation step, for each subdomain $D_i$, a cKRnet is initialized and the optimization problem \eqref{optimal_parameter} is solved associated with the data set $\mathcal{S}_i=\{(\xi_i^{(s)},\tau_i^{\infty}(\xi^{(s)}))\}_{s=1}^{\Non}$ using the stochastic gradient descent (SGD) method. More specifically, the set of samples $\mathcal{S}_i$ is first divided into $N_m$ mini-batches $\{\mathcal{S}_i^{(n)}\}_{n=1}^{N_m}$ and each mini-batch contains $N_b$ samples such that $\Non=N_m \times N_b$. Let $\Theta_{i,1}^{(0)}$ denote the initial parameters of cKRnet and $p_{\mathsf{cKR}}\left(\tau_i\mid\xi_i;\Theta_{i,1}^{(0)}\right)$ denote the approximate conditional PDF induced by the initial cKRnet. Then, the model parameters at the $n$-th iteration of $e$-th epoch, which is denoted by $\Theta_{i,e}^{(n)}$, can be updated as
\begin{eqnarray}
    \Theta_{i,e}^{(n)}=\Theta_{i,e}^{(n-1)}
    - \eta \nabla_{\Theta} \mathcal{L}\left(\mathcal{S}_i^{(n)};\Theta_{i,e}^{(n-1)}\right),\quad n=1,\dots,N_m, \quad e=1,\ldots,N_e,\label{sgd_alg} 
\end{eqnarray}
where $\eta>0$ is the learning rate and $N_e$ is the maximum epoch number. With the samples in the $n$-th mini-batch $\mathcal{S}_i^{(n)}=\{(\xi_i^{(s)},\tau_i^{\infty}(\xi^{(s)}))\}_{s=1}^{N_b}$, $\mathcal{L}(\mathcal{S}_i^{(n)};\Theta_{i,e}^{(n-1)})$ is computed through 
\begin{eqnarray}
\mathcal{L}\left(\mathcal{S}_i^{(n)};\Theta_{i,e}^{(n-1)}\right)=-\frac{1}{N_b}\sum_{s=1}^{N_b} \log p_{\mathsf{cKR}}\left(\tau_i^{\infty}\left(\xi^{(s)}\right) \mid \xi_i^{(s)};\Theta_{i,e}^{(n-1)}\right).\label{sgd_loss}
\end{eqnarray}
In this work, the Adam \cite{kingma2014adam} optimizer version of SGD is employed to accelerate the training process for cKRnet.

\begin{algorithm}
    \renewcommand{\algorithmicrequire}{\textbf{Input:}}
	\renewcommand{\algorithmicensure}{\textbf{Output:}}
    \caption{CKR-DDUQ online}
    \label{alg_on1}
    \begin{algorithmic}[1]
    \Statex \textbf{Input:} 
    system input PDF $\pi_{\xi}(\xi)$, coupling surrogates $\{\tilde{h}_{i,j}\}_{(i,j)\in \linter}$ and offline output samples $\left\{y_i\left(u(x,\xi_i^{(s)},\tau_i^{(s)})\right)\right\}_{s=1}^{\Noff}$ for $i=1,\dots,M$.
    \State Generate samples $\left\{\xi^{(s)}\right\}^{\Non}_{s=1}$ of the PDF $\pi_{\xi}(\xi)$.
    \For{$s=1:\Non$}
        \State Initialize the interface parameters $\tau^0_{i,j}\left(\xi^{(s)}\right)$ for {\lr $(i,j\,) \in \linter$}.
        \State Set $k=0$ and $\epsilon_0>tol$.
            \While{$\epsilon_k>tol$}
            \For{{\lr$i=1:M$}}
            \For{{\lr $j\in \linteri$}}
            \State $\tau^{k+1}_{i,j}(\xi^{(s)})=\theta_{i,j} \tilde{h}_{i,j}\left(\xi^{(s)}_i,\tau^k_i(\xi^{(s)})\right)+(1-\theta_{i,j})\tau^{k}_{i,j}(\xi^{(s)})$.
            \EndFor
            \EndFor
        \State Update $k=k+1$.
        \State Update $\epsilon_{k}=\max_{\left(i,j\,\right)\in \linter}\left(\left\|\tau^{k}_{i,j}\left(\xi^{\left(s\right)}\right)-\tau^{k-1}_{i,j}\left(\xi^{\left(s\right)}\right)\right\|_{\infty}\right)$.
            \EndWhile
            \State Set $\tau^{\infty}_{i,j}\left(\xi^{(s)}\right)=\tau^{k}_{i,j}\left(\xi^{(s)}\right)$ for {\lr $(i,j\,) \in \linter$}. 
    \EndFor
    \For {$i=1:M$}
        \State Let $\tau_i^{\infty}\left(\xi^{(s)}\right)=\otimes_{j\in\linteri} \tau^{\infty}_{i,j}\left(\xi^{(s)}\right)$ and $\mathcal{S}_i=\left\{\left(\xi_i^{(s)},\tau_i^{\infty}\left(\xi^{(s)}\right)\right)\right\}_{s=1}^{\Non}$.
        \State Divide $\mathcal{S}_i$ into $N_m$ mini-batches $\left\{\mathcal{S}_i^{(n)}\right\}_{n=1}^{N_m}$.
        \State Initialize a cKRnet and the induced conditional PDF is $p_{\mathsf{cKR}}\left(\tau_i\mid\xi_i;\Theta_{i,1}^{(0)}\right)$.
        \For {$e=1:N_e$}
            \For {$n=1:N_m$}
                \State Compute the loss function (see \eqref{sgd_loss}) on the mini-batch $\mathcal{S}_i^{(n)}$.
                \State Update $\Theta_{i,e}^{(n)}$ using the stochastic gradient descent method (see \eqref{sgd_alg}).
            \EndFor
            \If {$e<N_e$}
            \State Let $\Theta_{i,e+1}^{(0)}=\Theta_{i,e}^{(N_m)}$.
            \State Shuffle the mini-batches $\left\{\mathcal{S}_i^{(n)}\right\}_{n=1}^{N_m}$ of $\mathcal{S}_i$.
            \EndIf
        \EndFor
        \State Return $\hat{\pi}_{\tau_i\mid\xi_i}\left(\tau_i\mid\xi_i\right)=p_{\mathsf{cKR}}\left(\tau_i\mid\xi_i;\Theta_{i,N_e}^{(N_m)}\right)$.
        \State Re-weight the pre-computed offline local outputs: $\left\{w^{(s)}_iy_i\left(u\left(x,\xi^{(s)}_i,\tau^{(s)}_i\right)\right)\right\}^{\Noff}_{s=1}$ with $w^{(s)}_i=\frac{\hat{\pi}_{\tau_i\mid\xi_i}\left(\tau^{(s)}_i\mid\xi^{(s)}_i\right)}{p_{\tau_i}\left(\tau^{(s)}_i\right)}$.
    \EndFor  
    \Statex \textbf{Output:} weighted samples $\left\{w^{(s)}_iy_i\left(u\left(x,\xi^{(s)}_i,\tau^{(s)}_i\right)\right)\right\}^{\Noff}_{s=1}$, for $i=1,\ldots,M$.
    \end{algorithmic}
\end{algorithm}

\subsection{Convergence analysis}\label{analysis}
For the following analysis, we consider the outputs of interest $\{y_i(u(x,\xi)|_{D_i})\}_{i=1}^M$ of \eqref{spdexi1}--\eqref{spdexi2} to be scalar. It is straightforward to extend the analysis to the case of multiple outputs. Following the notation in \cite{liao2015domain}, for a given number $a$, the actual probability of $y_i(u(x,\xi)|_{D_i})\leq a$ is denoted by
$P(y_i\leq a)$, i.e.,
\begin{eqnarray*}
P(y_i \leq a):=\int_{\Gamma} I^a_i(\xi) \pi_{\xi}(\xi)\, {\rm{d}}\xi,
\end{eqnarray*}
where
\begin{eqnarray}
I^a_i(\xi):=\left\{\begin{array}{cc}
    1 & \textrm{if $y_i\left( u\left(x,\xi\right)\big|_{D_i}\right)\leq a$}\\
    0 & \textrm{if $y_i\left( u\left(x,\xi\right)\big|_{D_i}\right)> a$}
    \end{array}\right. . \nonumber
\end{eqnarray}
Next, the support of the proposal input PDF $p_{\tau_i}(\tau_i)$ is denoted by $\Gamma_{P,i}$ and the support of the local target marginal PDF $\pi_{\tau_i}(\tau_i)=\int_{\Gamma_i}\pi_{\xi_i,\tau_i}(\xi_i,\tau_i)\, {\rm{d}} \xi_i$ is denoted by $\Gamma_{T, i}$. In the following, it is required that $p_{\tau_i}(\tau_i)$ is chosen so that $\Gamma_{P,i}\supseteq \Gamma_{T,i}$. 
The probability of the local output $y_i( u(x,\xi_i,\tau_i))$ associated with the local target input conditional PDF $\pi_{\tau_i \mid \xi_i}(\tau_i \mid \xi_i)$ is defined by
\begin{eqnarray}
\tilde{P}(y_i \leq a):=\int_{\Gamma_{T,i}}\int_{\Gamma_{i}} \tilde{I}^a_i(\xi_i,\tau_i) \pi_{\tau_i \mid \xi_i}(\tau_i \mid \xi_i) \pi_{\xi_i}(\xi_i)\, {\rm{d}}\xi_i \, {\rm{d}}\tau_i, \label{target_p}
\end{eqnarray}
where
\begin{eqnarray}
\tilde{I}^a_i(\xi_i,\tau_i):=\left\{\begin{array}{cc}
    1 & \textrm{if $y_i\left( u\left(x,\xi_i,\tau_i\right)\right)\leq a$}\\
    0 & \textrm{if $y_i\left( u\left(x,\xi_i,\tau_i\right)\right)> a$}
\end{array}\right. .  \label{indicator_local}
\end{eqnarray}
\begin{lemma}\label{lemma_gl}
    If the DD-convergence condition is satisfied for all $\xi\in \Gamma$, then
    for any given number $a$, $\tilde{P}(y_i \leq a)$=$P(y_i \leq a)$.
\end{lemma}
The proof of Lemma \ref{lemma_gl} can be seen in \cite{liao2015domain}. At the last step of the CKR-DDUQ algorithm (line 31 of Algorithm \ref{alg_on1}), weighted samples are obtained. We next define the probability associated with the weighted samples. To begin with, the  exact weights are considered:
\begin{eqnarray}
\mw^{(s)}_i=\frac{\pi_{\tau_i\mid\xi_i}\left(\tau^{(s)}_i\mid\xi^{(s)}_i\right)}{p_{\tau_i}\left(\tau^{(s)}_i\right)},\quad s=1,\dots,\Noff,\quad i=1,\dots,M,\label{w_exact}
\end{eqnarray}
which are associated with the exact local target input conditional PDF $\pi_{\tau_i \mid \xi_i}(\tau_i \mid \xi_i)$ rather than the estimated conditional PDF $\hat{\pi}_{\tau_i \mid \xi_i}(\tau_i \mid \xi_i)$ on line 31 of Algorithm \ref{alg_on1}.
With these exact weights, we define
\begin{eqnarray}
    P_{\mw_i}(y_i \leq a):= \frac{\sum^{\Noff}_{s=1}\mw^{(s)}_i\tilde{I}^a_i\left(\xi^{(s)}_i,\tau^{(s)}_i\right)}{\sum^{\Noff}_{s=1}\mw^{(s)}_i}.\label{phat_exact}
\end{eqnarray}
For the CKR-DDUQ outputs (see line 31 of Algorithm \ref{alg_on1}) with weights $\{w^{(s)}_i\}^{\Noff}_{s=1}$ obtained from the estimated conditional PDFs $\hat{\pi}_{\tau_i \mid \xi_i}(\tau_i \mid \xi_i)$ for  $i=1,\dots,M$, we define
\begin{eqnarray}
P_{w_i}(y_i \leq a):= \frac{\sum^{\Noff}_{s=1}w^{(s)}_i\tilde{I}^a_i\left(\xi^{(s)}_i,\tau^{(s)}_i\right)}{\sum^{\Noff}_{s=1}w^{(s)}_i},\label{phat}
\end{eqnarray}
where 
\begin{eqnarray}
    w^{(s)}_i=\frac{\hat{\pi}_{\tau_i\mid\xi_i}\left(\tau^{(s)}_i\mid\xi^{(s)}_i\right)}{p_{\tau_i}\left(\tau^{(s)}_i\right)},\quad s=1,\dots,\Noff,\quad i=1,\dots,M.\label{w_estimate}
\end{eqnarray}
In the following, $P_{w_i}(y_i \leq a)$ is shown to approach the actual probability $P(y_i\leq a)$ as the offline sample size $\Noff$  and the online sample size $\Non$ increase. The analysis takes two steps: the first step is to show that $P_{\mw_i}(y_i \leq a)\to P(y_i \leq a)$ as $\Noff$ goes to infinity and the second step is to show that $P_{w_i}(y_i \leq a)\to P_{\mw_i}(y_i \leq a)$ as $\Non$ goes to infinity.
\begin{lemma}\label{lemma_IS}
If the DD-convergence condition is satisfied for all $\xi\in \Gamma$, and the support of the local proposal input PDF $p_{\tau_i}(\tau_i)$ covers the support of the local target marginal PDF $\pi_{\tau_i}(\tau_i)=\int_{\Gamma_i}\pi_{\xi_i,\tau_i}(\xi_i,\tau_i)\, {\rm{d}} \xi_i$ i.e., $\Gamma_{P,i} \supseteq \Gamma_{T,i}$, then for any given number $a$, the following equality holds
    \begin{eqnarray}
    \lim_{\Noff \to \infty} P_{\mw_i}(y_i \leq a)=P(y_i \leq a).
    \label{lemma4.2_eq}
    \end{eqnarray} 
\end{lemma}
\begin{proof}
Following the standard convergence analysis of importance sampling in \cite{amaral2014decomposition,smith1992bayesian} and from equations \eqref{w_exact} and \eqref{phat_exact}, we have
\begin{eqnarray}
    \lim_{\Noff \to \infty}P_{\mw_i}(y_i \leq a)
    &=&\lim_{\Noff \to \infty}\frac{\sum^{\Noff}_{s=1}\mw^{(s)}_i\tilde{I}^a_i\left(\xi^{(s)}_i,\tau^{(s)}_i\right)}{\sum^{\Noff}_{s=1}\mw^{(s)}_i}\nonumber\\
    &=&\lim_{\Noff \to \infty}\frac{\frac{1}{\Noff}\sum^{\Noff}_{s=1}  \frac{\pi_{\tau_i \mid \xi_i}\left(\tau^{(s)}_i \mid \xi^{(s)}_i \right)}{p_{\tau_i}\left(\tau^{(s)}_i\right)} \tilde{I}^a_i\left(\xi^{(s)}_i,\tau^{(s)}_i\right)}{\frac{1}{\Noff}\sum^{\Noff}_{s=1} \frac{\pi_{\tau_i \mid \xi_i}\left(\tau^{(s)}_i \mid \xi^{(s)}_i \right)}{p_{\tau_i}\left(\tau^{(s)}_i\right)}}\nonumber,
    \end{eqnarray}
to which we apply the convergence property of Monte Carlo integration \cite{caflisch1998monte} and obtain
\begin{eqnarray}
    \lim_{\Noff \to \infty}P_{\mw_i}(y_i \leq a)
    &=&\frac{\int_{{\Gamma}_{P,i}} \int_{{\Gamma}_{i}} \frac{\pi_{\tau_i \mid \xi_i}\left(\tau_i \mid \xi_i\right)}{p_{\tau_i}\left(\tau_i\right)} \tilde{I}^a_i(\xi_i,\tau_i)p_{\tau_i}\left(\tau_i\right) \pi_{\xi_i}\left(\xi_i\right)\, {\rm{d}}\xi_i \, {\rm{d}}\tau_i}
    {\int_{{\Gamma}_{P,i}} \int_{{\Gamma}_{i}} \frac{\pi_{\tau_i\mid \xi_i}\left(\tau_i\mid\xi_i\right)}{p_{\tau_i}\left(\tau_i\right)}p_{\tau_i}\left(\tau_i\right) \pi_{\xi_i}\left(\xi_i\right)\, {\rm{d}}\xi_i \, {\rm{d}}\tau_i}                                      \nonumber \\
    &=&\frac{\int_{{\Gamma}_{P,i}} \int_{{\Gamma}_{i}} \pi_{\tau_i \mid \xi_i}\left(\tau_i \mid \xi_i\right) \tilde{I}^a_i(\xi_i,\tau_i) \pi_{\xi_i}\left(\xi_i\right)\, {\rm{d}}\xi_i \, {\rm{d}}\tau_i}
    {\int_{{\Gamma}_{P,i}} \int_{{\Gamma}_{i}} \pi_{\tau_i\mid \xi_i}\left(\tau_i\mid\xi_i\right)\pi_{\xi_i}\left(\xi_i\right)\, {\rm{d}}\xi_i \, {\rm{d}}\tau_i}       \nonumber .
    \end{eqnarray}
Combining $\Gamma_{P,i} \supseteq \Gamma_{T,i}$ and Lemma \ref{lemma_gl} gives
\begin{eqnarray}
    \lim_{\Noff \to \infty}P_{\mw_i}(y_i \leq a) 
    &=&\frac{\int_{{\Gamma}_{T,i}}\int_{\Gamma_{i}} \pi_{\tau_i \mid \xi_i}\left(\tau_i \mid \xi_i\right) \tilde{I}^a_i(\xi_i,\tau_i) \pi_{\xi_i}\left(\xi_i\right)\, {\rm{d}}\xi_i \, {\rm{d}}\tau_i}
    {\int_{{\Gamma}_{T,i}}\int_{\Gamma_{i}} \pi_{\tau_i\mid \xi_i}\left(\tau_i\mid\xi_i\right)\pi_{\xi_i}\left(\xi_i\right)\, {\rm{d}}\xi_i \, {\rm{d}}\tau_i}                                      \nonumber \\  
    &=&\int_{{\Gamma}_{T,i}} \int_{\Gamma_{i}}\tilde{I}^a_i(\xi_i,\tau_i) \pi_{\tau_i \mid \xi_i}\left(\tau_i \mid \xi_i\right)  \pi_{\xi_i}\left(\xi_i\right)\, {\rm{d}}\xi_i \, {\rm{d}}\tau_i   \nonumber \\  
    &=&\tilde{P}(y_i \leq a)\nonumber \\
    &=&P(y_i \leq a).
    \end{eqnarray}
\end{proof}

In the following analysis, two estimated PDFs
$p_{\mathsf{cKR}}(\tau_i\mid\xi_i;\Theta^*)$ and $p_{\mathsf{cKR}}(\tau_i\mid\xi_i;\Theta_{\Non}^*)$ are considered, where $\Theta^*$ and $\Theta_{\Non}^*$  are defined as 
\begin{eqnarray}
&&\Theta^{*}=\arg \min_{\Theta} \mathcal{L}\left(p_{\mathsf{cKR}}\left(\tau_i\mid\xi_i;\Theta\right)\right)  =\arg \min_{\Theta} - \mathbb{E}_{\pi_{ \xi_i,\tau_i}(\xi_i,\tau_i)} \log p_{\mathsf{cKR}}\left(\tau_i\mid\xi_i;\Theta\right),\label{pdf1}\\
&&\Theta_{\Non}^{*}=\arg \min_{\Theta} \hat{\mathcal{L}}\left(p_{\mathsf{cKR}}\left(\tau_i\mid\xi_i;\Theta\right)\right) = \arg \min_{\Theta} - \frac{1}{\Non}\sum_{s=1}^{\Non} \log p_{\mathsf{cKR}}\left(\tau^{\infty}_i\left(\xi^{(s)}\right)\mid\xi_i^{(s)};\Theta\right),\label{pdf2}
\end{eqnarray}
where $\{(\xi^{(s)}_i, \tau^{\infty}_i(\xi^{(s)}))\}^{\Non}_{s=1}$ are the target samples in the online stage (see line 17 of Algorithm \ref{alg_on1}).
\begin{lemma}\label{lemma_online}
    If the target samples $\{(\xi^{(s)}_i , \tau^{\infty}_i(\xi^{(s)}))\}^{\Non}_{s=1}$ are independent and identically distributed (i.i.d.) data following the distribution $\pi_{\xi_i,\tau_i}\left(\xi_i, \tau_i\right)$ for $i=1,\dots,M$ and the minimizer of $\mathcal{L}(p_{\mathsf{cKR}}(\tau_i\mid\xi_i;\Theta))$ is unique,
    then for any given number $a$, the following equality holds
    \begin{eqnarray}
    \lim_{\Non \to \infty} P_{w_i}(y_i \leq a)=P_{\mw_i}(y_i \leq a).
    \label{lemma4.3_eq}
    \end{eqnarray} 
\end{lemma}

\begin{proof}
We have
\begin{eqnarray*}
\begin{aligned}
p_{\mathsf{cKR}}\left(\tau_i\mid\xi_i;\Theta_{\Non}^*\right)-\pi_{\tau_i\mid\xi_i}\left(\tau_i\mid\xi_i \right)
=&p_{\mathsf{cKR}}\left(\tau_i\mid\xi_i;\Theta_{\Non}^*\right)-
p_{\mathsf{cKR}}\left(\tau_i\mid\xi_i;\Theta^*\right)\\
&+
p_{\mathsf{cKR}}\left(\tau_i\mid\xi_i;\Theta^*\right)-
\pi_{\tau_i\mid\xi_i}\left(\tau_i\mid\xi_i \right).
\end{aligned}
\end{eqnarray*}
Then $\forall (\xi_i^{(s)},\tau_i^{(s)}) \in \Gamma_{i}\times \Gamma_{P,i}$, the following inequality holds
\begin{eqnarray*}
\begin{aligned}
\left|p_{\mathsf{cKR}}\left(\tau_i^{(s)}\mid\xi_i^{(s)};\Theta_{\Non}^*\right)-\pi_{\tau_i\mid\xi_i}\left(\tau_i^{(s)}\mid\xi_i^{(s)} \right)\right|
\leq&
\underbrace{\left|p_{\mathsf{cKR}}\left(\tau_i^{(s)}\mid\xi_i^{(s)};\Theta_{\Non}^*\right)-
p_{\mathsf{cKR}}\left(\tau_i^{(s)}\mid\xi_i^{(s)};\Theta^*\right)\right|}_{e_1}\\
&+
\underbrace{\left|p_{\mathsf{cKR}}\left(\tau_i^{(s)}\mid\xi_i^{(s)};\Theta^*\right)-
\pi_{\tau_i\mid\xi_i}\left(\tau_i^{(s)}\mid\xi_i^{(s)} \right)\right|}_{e_2},
\end{aligned}
\end{eqnarray*}
where $e_1$ represents the statistical error and $e_2$ represents the approximation error. Since $\{(\xi^{(s)}_i , \tau^{\infty}_i(\xi^{(s)}))\}^{\Non}_{s=1}$ are i.i.d. data following the target PDF $\pi_{\xi_i,\tau_i}\left(\xi_i, \tau_i\right)$, according to the Law of Large Numbers, we have
\begin{eqnarray*}
&&\lim_{\Non \to \infty}  \hat{\mathcal{L}}\left(p_{\mathsf{cKR}}\left(\tau_i\mid\xi_i;\Theta\right)\right)=\mathcal{L}\left(p_{\mathsf{cKR}}\left(\tau_i\mid\xi_i;\Theta\right)\right).
\end{eqnarray*}
Then since the minimizer of $\mathcal{L}(p_{\mathsf{cKR}}(\tau_i\mid\xi_i;\Theta))$ is assumed to be unique, we obtain
\begin{eqnarray*}
\lim_{\Non \to \infty}  \Theta_{\Non}^*=\Theta^*.
\end{eqnarray*}
Thus
\begin{eqnarray*}
\lim_{\Non \to \infty}e_1=0,
\end{eqnarray*}
which yields that
\begin{eqnarray}
\lim_{\Non \to \infty}\left|p_{\mathsf{cKR}}\left(\tau_i^{(s)}\mid\xi_i^{(s)};\Theta_{\Non}^*\right)-\pi_{\tau_i\mid\xi_i}\left(\tau_i^{(s)}\mid\xi_i^{(s)} \right)\right|\leq e_2, \quad 
\forall (\xi_i^{(s)},\tau_i^{(s)}) \in \Gamma_{i}\times \Gamma_{P,i}.\label{e1}
\end{eqnarray}
The approximation error $e_2$ depends on the model capability of our cKRnet. Recent researches \cite{teshima2020coupling,ishikawa2023universal} show that normalizing flow models based on affine coupling can be universal distributional approximators. 
Given that our cKRnet is one of the affine coupling-based normalizing flow models, it is assumed that $e_2=0$ in this work, and \eqref{e1} then turns into 
\begin{eqnarray*}
\lim_{\Non \to \infty}p_{\mathsf{cKR}}\left(\tau_i^{(s)}\mid\xi_i^{(s)};\Theta_{\Non}^*\right)=\pi_{\tau_i\mid\xi_i}\left(\tau_i^{(s)}\mid\xi_i^{(s)} \right),\quad 
 \forall (\xi_i^{(s)},\tau_i^{(s)}) \in \Gamma_{i}\times \Gamma_{P,i}.
\end{eqnarray*}
With $\hat{\pi}_{\tau_i\mid\xi_i}(\tau_i\mid\xi_i)=p_{\mathsf{cKR}}(\tau_i\mid\xi_i;\Theta_{\Non}^*)$, we have
\begin{eqnarray}
&&\lim_{\Non \to \infty}\hat{\pi}_{\tau_i\mid\xi_i}\left(\tau_i^{(s)}\mid\xi_i^{(s)}\right)=\pi_{\tau_i \mid \xi_i}\left(\tau_i^{(s)} \mid \xi_i^{(s)}\right) ,\quad 
\forall (\xi_i^{(s)},\tau_i^{(s)}) \in \Gamma_{i}\times \Gamma_{P,i},\nonumber\\
\Rightarrow &&\lim_{\Non \to \infty} w_i^{(s)}=\mw_i^{(s)}, \quad s=1,\dots,\Noff,\nonumber\\
\Rightarrow &&\lim_{\Non \to \infty}P_{w_i}(y_i \leq a)=P_{\mw_i}(y_i \leq a).\label{same_weights}
\end{eqnarray}
\end{proof}
\begin{theorem}\label{theorem}
    Under the same conditions of Lemma \ref{lemma_IS} and Lemma \ref{lemma_online}. For any given number $a$, the following equality holds
    \begin{eqnarray}
    \lim_{\substack{\Non \to \infty \\ \Noff \to \infty}} P_{w_i}(y_i \leq a)=P(y_i \leq a).
    \label{thorem_overall}
    \end{eqnarray}
\end{theorem}
\begin{proof}
Combing \eqref{lemma4.2_eq} in Lemma \ref{lemma_IS} and \eqref{lemma4.3_eq} in Lemma \ref{lemma_online}, this theorem is proven. 
\end{proof}

%% file: sections/section5.tex
In this section, numerical experiments are conducted to illustrate the effectiveness of our CKR-DDUQ (cKRnet based domain decomposed uncertainty quantification) approach. Four test problems are considered: a two-component diffusion problem, a three-component diffusion problem, a random domain decomposition problem, and a Stokes problem. The local Monte Carlo simulation in the offline stage is conducted on an Intel Xeon Gold 6130 2.1 GHz processor with 96 GB RAM using Matlab 2018a. The coupling surrogates modeling and the whole online stage are conducted on an NVIDIA GTX 1080 Ti GPU with 64 GB RAM using Python 3.7 and TensorFlow 2.0.

\subsection{Two-component diffusion problem with \texorpdfstring{$28$}{28} parameters}\label{2domain_test}
\begin{figure}[!htp]
	\setlength{\unitlength}{1.5cm}
	\centerline{
	\begin{picture}(4,2.3)(-0.2,-0.2)
        \put(0,0){\line(1,0){4}}
        \put(0,0){\line(0,1){2}}
        \put(0,2){\line(1,0){4}}
        \put(4,0){\line(0,1){2}}
        \put(0.6,1){$D_1$}
        \put(3.1,1){$D_2$}
        \put(0,0){\circle*{0.1}}
        \put(-0.3,-0.25){\scriptsize (0,0)}
        \put(4,0){\circle*{0.1}}
        \put(3.8,-0.25){\scriptsize (2,0)}
        \put(4,2){\circle*{0.1}}
        \put(4.05,1.95){\scriptsize  (2,1)}
        \put(2,0){\circle*{0.1}}
        \put(1.8,-0.25){\scriptsize (1,0)}
        \put(-1.5,0.5){\vector(1,0){0.5}}
        \put(-1.5,0.5){\vector(0,1){0.5}}
        \put(-1.4,0.2){$x_1$}
        \put(-1.85,0.6){$x_2$}
        \linethickness{0.7mm}
        \put(2,0){\line(0,1){2}}
	\end{picture}
	}
	\caption{Illustration of the spatial domain with two components.}
	\label{fig_2d}
\end{figure}
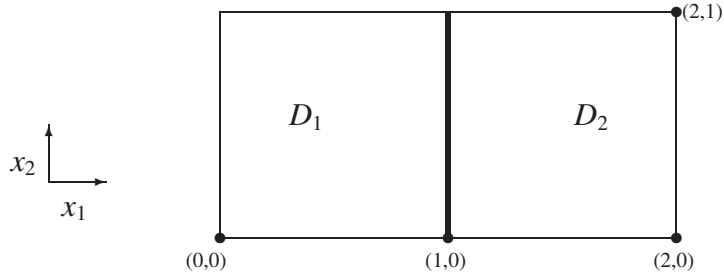

We start with the following diffusion problem, while the governing equation is
\begin{eqnarray}
    {\lr -\div \left(a\left(x,\xi\,\right)\nabla u\left(x,\xi\,\right)\right)}=f(x,\xi\,)
    & \quad \textrm{in}& D\times\Gamma, \label{ht1} \\
   u\left(x,\xi\,\right)=0
    & \quad\textrm{on}& \partial D \times\Gamma,\label{ht2}
\end{eqnarray}
where $a(x,\xi)$ is the permeability coefficient and $f$ is the source function. For this test problem, the source term is specified as $f(x,\xi)=100$ and the physical domain $D=(0,2)\times(0,1)\subset \dsR^{2}$ is decomposed into $M=2$ subdomains as shown in Figure \ref{fig_2d}.

On each local subdomain $D_i$ (for $i=1,\ldots,M$), $a(x,\xi)|_{D_i}$ is assumed to be a random field with mean function $a_{i,0}(x)$, constant standard deviation $\sigma$ and covariance function $C(x,x')$
\begin{eqnarray}
    C(x,x')=\sigma^2 \exp\left(-\frac{|x_1-x_1'|}{L_c}-\frac{|x_2-x_2'|}{L_c}
    \right),\label{covariance}
\end{eqnarray}
where $x=[x_1,x_2]^{\top}$, $x'=[x_1',x_2']^{\top}$ and $L_c$ is the correlation length. 
The random fields are assumed to be independent between different subdomains. 
Each random field is approximated with the truncated Karhunen--Lo\`eve (KL) expansion \cite{babuska1,elmmil11,ghaspa03} as follows
\begin{eqnarray}
    a(x,\xi\,)|_{D_i}\approx
    a_{i,0}(x)+\sum_{m=1}^{\mathfrak{N}_i}\sqrt{\lambda_{i,m}}a_{i,m}(x)\xi_{i,m}, \quad i=1,\ldots,M,
    \label{kl_expan}
\end{eqnarray}
where $\mathfrak{N}_i$ is the number of KL modes retained for subdomain $D_i$, $\{a_{i,m}(x)\}^{\mathfrak{N}_i}_{m=1}$ and $\{\lambda_{i,m}\}^{\mathfrak{N}_i}_{m=1}$ are the eigenfunctions and eigenvalues of the covariance function \eqref{covariance} defined on $D_i$, and $\{\xi_{i,m}: i=1,\ldots,M\textrm{ and } m=1,\ldots,\mathfrak{N}_i\}$ are  uncorrelated random variables. 

In this test problem, we set $a_{1,0}(x)=a_{2,0}(x)=2$, $\sigma=0.5$ and $L_c=1$. The number of retained KL modes is chosen as $\mathfrak{N}_1=\mathfrak{N}_2=14$. Thus, the input vector $\xi$ comprises $28$ parameters with $\xi_1:=[\xi_{1,1},\xi_{1,2},\ldots,
\xi_{1,\mathfrak{N}_1}]^{\top}\in\dsR^{14}$ and $\xi_2:=[\xi_{2,1},\xi_{2,2},\ldots,\xi_{2,\mathfrak{N}_2}]^{\top}\in\dsR^{14}$ , where the random variables $\xi_{i,m}$ are set to be independent truncated Gaussian distributions with mean $0$, standard deviation $0.5$, and range $[-1,1]$. 
The outputs of interest are defined as 
\begin{eqnarray}
y_1(\xi)=\int_{Z_1} u(x,\xi) \,{\rm d} x_2,\label{o1}\\
y_2(\xi)=\int_{Z_2} u(x,\xi) \,{\rm d} x_2, \label{o2}
\end{eqnarray}
where $x=[x_1,x_2]^{\top}$, $Z_1:=\{x|\>x_1=0.5,\>0\leq x_2\leq 1\}$ and  $Z_2:=\{x|\>x_1=1.5,\>0\leq x_2 \leq 1\}$.

To decompose the global problem \eqref{ht1}--\eqref{ht2}, the parallel Dirichlet-Neumann domain decomposition method \cite{quavalbook} is employed. Specifically, on subdomain $D_1$, a local problem with a Dirichlet condition posed on $\partial_2 D_1$ is solved, while on subdomain $D_2$, a local problem with a Neumann condition imposed on $\partial_1 D_2$ is solved. The coupling functions are
\begin{eqnarray}
h_{2,1}:=u\left(x,\xi_2,\tau_2\right)\big|_{\partial_1 D_2}, \label{h21}\\
h_{1,2}:=\frac{\partial u\left(x,\xi_1,\tau_1\right)}{\partial n}\bigg|_{\partial_2 D_1}, \label{h12}
\end{eqnarray}
where the interface parameters $\tau_1$ and $\tau_2$ are defined in \eqref{lam_in}. In this two-component test problem, we have $\tau_1=\tau_{2,1}$ and $\tau_2=\tau_{1,2}$. The acceleration parameters in \eqref{ddupdate} are set to $\theta_{2,1} = 0.1$ and $\theta_{1,2}=0$ in this section. During the local Monte Carlo simulation procedure in the offline stage, each local solution is computed using a bilinear finite element approximation \cite{braess2001finite,elman2014finite} with a mesh size $h=1/16$. The interface functions $\{g_{i,j}\}_{(i,j)\in\linter}$ defined on the interface $\pjoi$ then have discrete forms associated with the grids on the interface. To obtain a lower-dimensional representation of the interface parameters $\tau_{i,j}$, the proper orthogonal decomposition (POD) method \cite{gunzburger2007reduced,holmes2012turbulence} is used here. That is, before the offline step, we generate a limited number of samples of $\xi$ and employ the domain decomposition method to solve the problem \eqref{ht1}--\eqref{ht2} for each sample. This process yields snapshots of the interface functions. With the snapshots, the POD bases for each interface function $g_{i,j}$ can be constructed, and the coefficients of the POD bases are treated as the interface parameters. In this test problem, the number of samples generated to construct the POD bases is set to one hundred. To obtain sufficiently accurate representations of interface functions, two POD modes for $g_{2,1}$ and six POD modes for $g_{1,2}$ are retained, i.e., $\tau_1\in\dsR^{2}$ and $\tau_2\in\dsR^{6}$. The local inputs for subdomain $D_1$ and $D_2$ are then $[\xi_1^{\top},\tau_1^{\top}]^{\top}\in\dsR^{16}$ and $[\xi_2^{\top},\tau_2^{\top}]^{\top}\in\dsR^{20}$ respectively. 

\begin{figure}[!htp]
    \centerline{
    \begin{tabular}{cc}
    \includegraphics[width=0.37\textwidth]{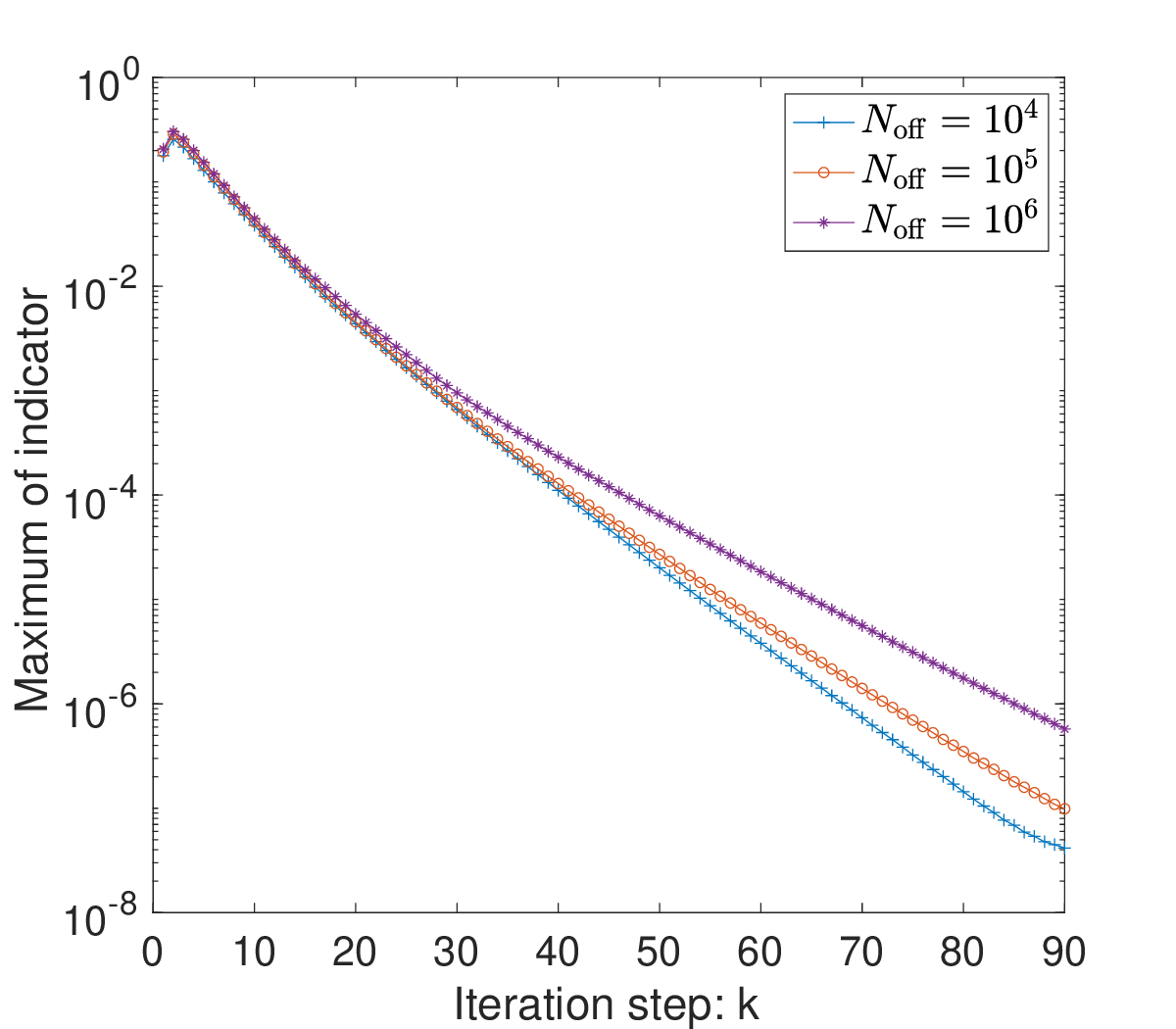}&
    \includegraphics[width=0.37\textwidth]{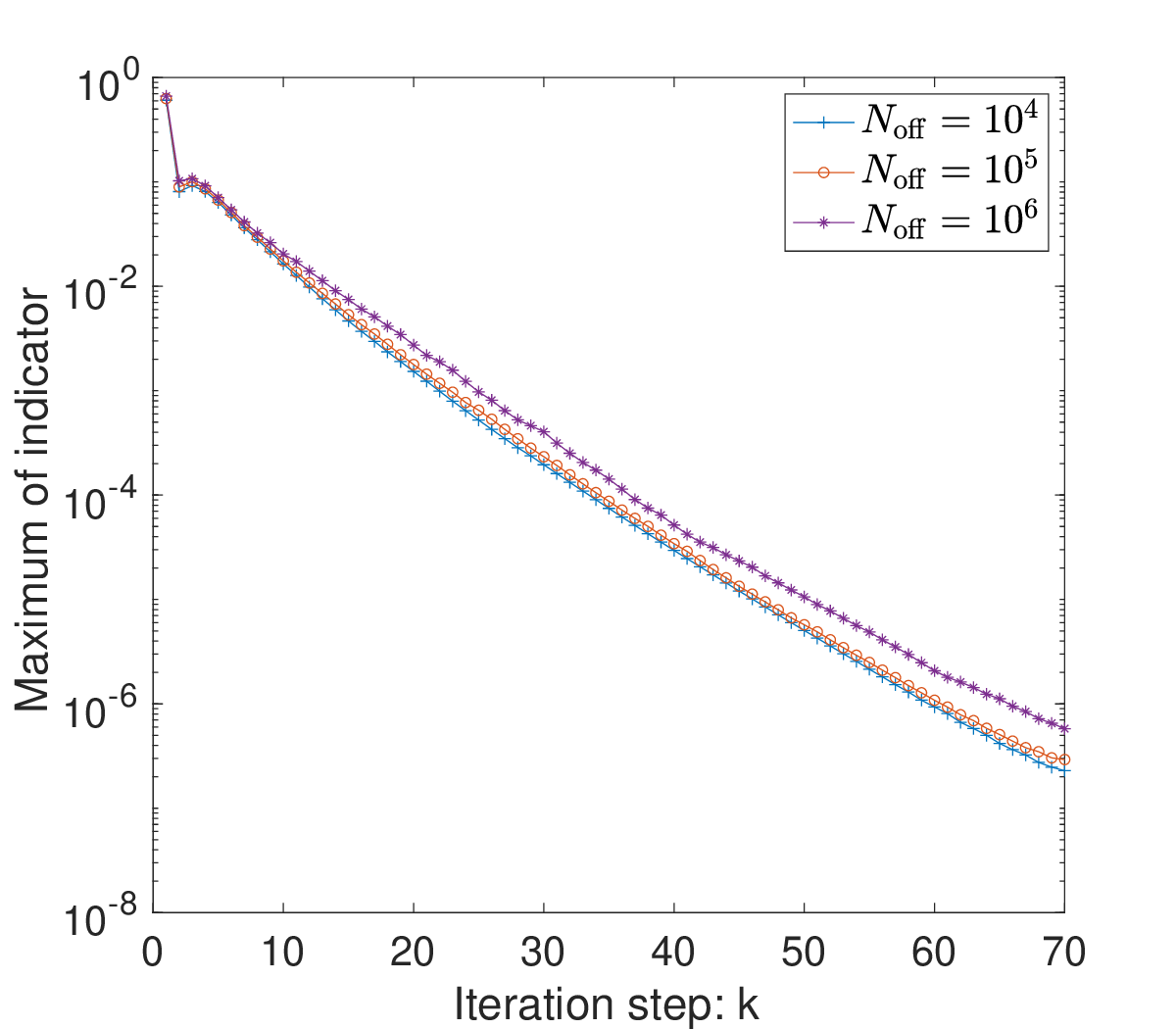}\\
    (a) Error indicator on $D_1$ & (b) Error indicator on $D_2$
    \end{tabular}}
    \caption{\lr
    Maximum of the error indicator on subdomain $D_1$ ($\max_{s=1:{\Noff}}\|\tau^{k+1}_1(\xi^{(s)})-\tau^k_1(\xi^{(s)})\|_{\infty}$) and that on subdomain $D_2$ ($\max_{s=1:{\Noff}}\|\tau^{k+1}_2(\xi^{(s)})-\tau^k_2(\xi^{(s)})\|_{\infty}$) for the coupling surrogates, two-component diffusion test problem.}
    \label{fig_2d_surrogate_convergence}
\end{figure}

In the offline stage, the proposal PDF $p_{\tau_1}(\tau_1)$ and $p_{\tau_2}(\tau_2)$ are chosen as two multivariate Gaussian distributions, of which the mean vectors and the covariance matrices are estimated with the snapshots generated in the POD step. The coupling surrogates are constructed with ResNets \cite{he2016deep}, where each ResNet comprising ten hidden layers with an equal width of sixty four neurons and the LeakyReLU activation function is utilized. In the online stage, we choose $\Non=\Noff$ for each subdomain to simplify the illustration and set the tolerance of the domain decomposition iteration to $tol=10^{-6}$. To show the domain decomposition iteration convergence property, the maximum of the error indicator on each subdomain $D_i$ (for $i=1,2$) is computed, i.e., $\max_{s=1:\Noff}\|\tau^{k+1}_i(\xi^{(s)})-\tau^k_i(\xi^{(s)})\|_{\infty}$ where $\tau^{k+1}_i(\xi^{(s)})$ is updated as line 8 of Algorithm \ref{alg_on1}. Figure \ref{fig_2d_surrogate_convergence} shows the maximum of the error indicator with respect to domain decomposition iterations with $\Noff=10^4,10^5$ and $10^6$, where the iterations are conducted using the coupling surrogates $ \{\tilde{h}_{i,j}(\xi_i,\tau_i)\}_{(i,j\,)\in \linter}$. It is clear that the error indicators reduce exponentially as iteration step $k$ increases.  

When estimating the target conditional PDF in the online stage (line $16$  to line $30$ of Algorithm \ref{alg_on1}),  cKRnet is set to contain $R=2$ invertible mappings (i.e., $f_{\bc,[r]}$ defined in \eqref{inner_loop} or the nonlinear invertible layer $L_N$) for subdomain $D_1$ and $R=3$ invertible mappings for subdomain $D_2$. Each invertible mapping except $L_N$ has $L=4$ conditional affine coupling layers, where $\bs_{\bc}$ and $\bt_{\bc}$ are implemented through a fully connected network containing two hidden layers with thirty two neurons in each hidden layer. In addition, the rectified linear unit function (ReLU) function is employed as the activation function. For comparison, we directly apply KRnet in the standard DDUQ online algorithm (see \cite{liao2015domain}), where the target input PDF $\pi_{\xi_i,\tau_i}(\xi_i,\tau_i)$ is estimated using KRnet (see \cite{tang2020deep}). This direct combination of KRnet and DDUQ is referred to as KR-DDUQ in the following, and the setting of KRnet is the same as that for cKRnet as above. Both KRnet and cKRnet are trained by the Adam \cite{kingma2014adam} optimizer with default settings and a learning rate of $0.001$. 

With the trained cKRnet, the offline output samples $\{y_i(u(x,\xi^{(s)}_i,\tau^{(s)}_i))\}^{\Noff}_{s=1}$ are re-weighted by weights $\{w^{(s)}_i\}_{s=1}^{\Noff}$ for each output (see line $31$ of Algorithm \ref{alg_on1}). For simplicity, the offline output samples are rewritten as $\{y_i(\xi^{(s)}_i,\tau^{(s)}_i)\}^{\Noff}_{s=1}$ in the following. Then the mean and the variance of each output estimated using CKR-DDUQ (the weighted samples obtained using Algorithm \ref{alg_on1}) are computed as 
\begin{eqnarray}
    \pE_{\Noff}\left(y_i\right)&:=&\frac{\sum^{\Noff}_{s=1}w^{(s)}_iy_i\left(\xi^{(s)}_i,\tau^{(s)}_i\right)}{\sum^{\Noff}_{s=1}w^{(s)}_i},\label{dduq_mean}\\
    \pV_{\Noff}\left(y_i\right)&:=&\frac{\sum^{\Noff}_{s=1}w^{(s)}_i\bigg(y_i\left(\xi^{(s)}_i,\tau^{(s)}_i\right)-\pE_{\Noff}\left(y_i\right)\bigg)^2}{\sum^{\Noff}_{s=1}w^{(s)}_i}.\label{dduq_var}
\end{eqnarray}
Furthermore, the mean and the variance of each output estimated using KR-DDUQ are defined in the same way except that the weights $\{w^{(s)}_i\}_{s=1}^{\Noff}$ are computed by \eqref{dduq_weights}.

To compare the mean and variance estimates of CKR-DDUQ and KR-DDUQ, a set of reference output samples are obtained by solving the global problem \eqref{ht1}--\eqref{ht2} using the Monte Carlo method with $N_{\rm ref}=10^7$ samples. Denoting the output samples as $\{y^{\rm{ref}}_i(\xi^{(s)})\}^{N_{\rm {ref}}}_{s=1}$ for $i=1,2$, the reference mean and variance of each output are computed as
\begin{eqnarray}
\pE_{N_{\textrm{ref}}}\left(y^{\rm ref}_i\right)&:=&\sum^{N_{\textrm{ref}}}_{s=1}\frac{1}{N_{\textrm{ref}}}y^{\rm ref}_i\left(\xi^{(s)}\right),\label{ref_mean}\\
\pV_{N_{\textrm{ref}}}\left(y^{\rm ref}_i\right)&:=&\sum^{N_{\textrm{ref}}}_{s=1}\frac{1}{N_{\textrm{ref}}}\bigg(y^{\rm ref}_i\left(\xi^{(s)}\right)-\pE_{N_{\textrm{ref}}}\left(y^{\rm ref}_i\right)\bigg)^2.\label{ref_var}
\end{eqnarray}
Next, we assess the errors of CKR-DDUQ and KR-DDUQ in estimating the mean and variance with the following quantities
\begin{eqnarray}
    \epsilon_i&:=&\left|\frac{\pE_{\Noff}\left(y_i\right)-\pE_{N_{\textrm{ref}}}\left(y^{\rm{ref}}_i\right)} {\pE_{N_{\textrm{ref}}}\left(y^{\rm{ref}}_i\right)}\right|,\label{dduq_mean_error} \\
    \eta_i&:=&\left|\frac{\pV_{\Noff}\left(y_i\right)-\pV_{N_{\textrm{ref}}}\left(y^{\rm{ref}}_i\right)} {\pV_{N_{\textrm{ref}}}\left(y^{\rm{ref}}_i\right)} \right|. \label{dduq_var_error}
\end{eqnarray}

\begin{figure}[!htp]
    \centerline{
    \begin{tabular}{cc}
    \includegraphics[width=0.37\textwidth]{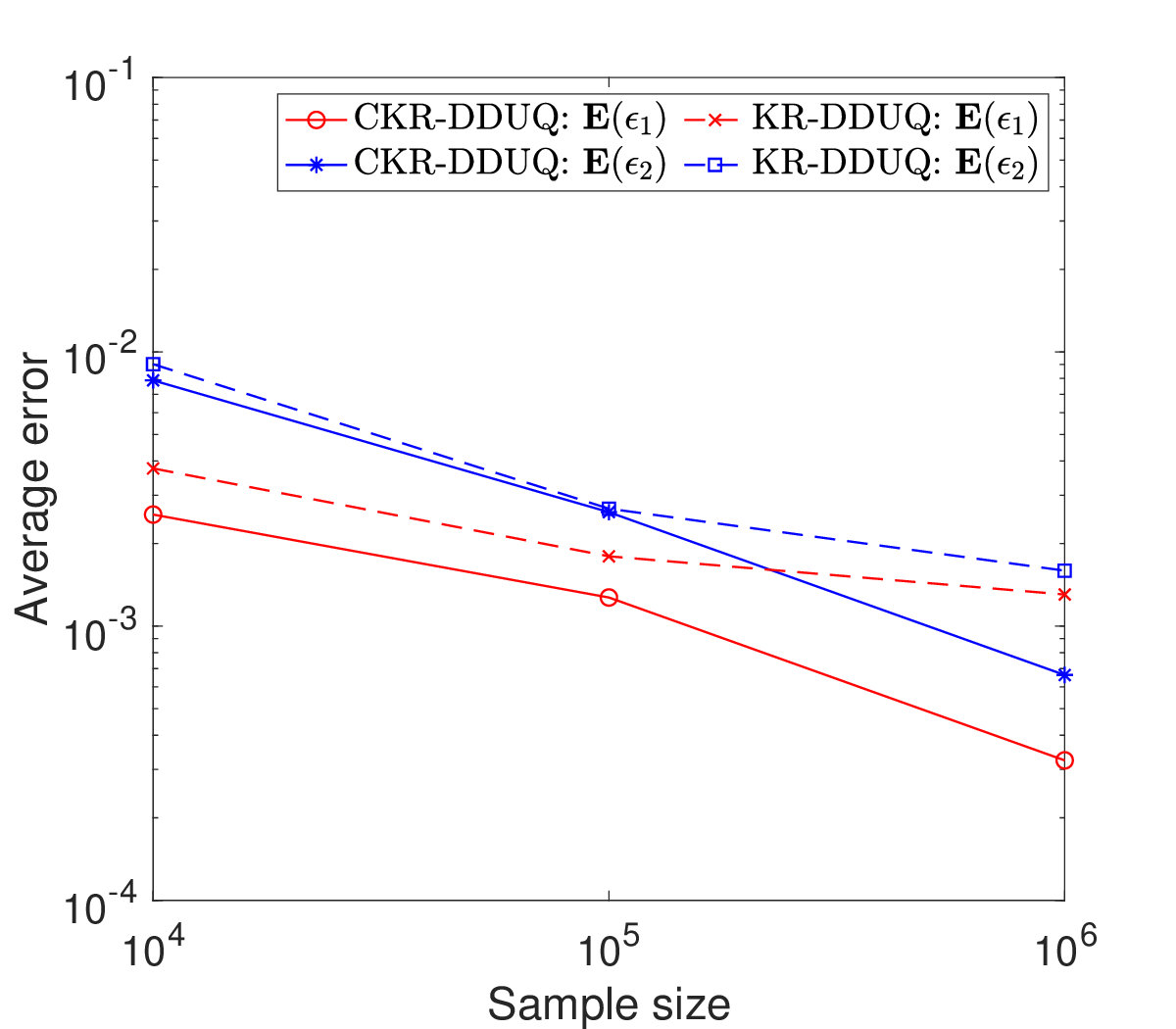}&
    \includegraphics[width=0.37\textwidth]{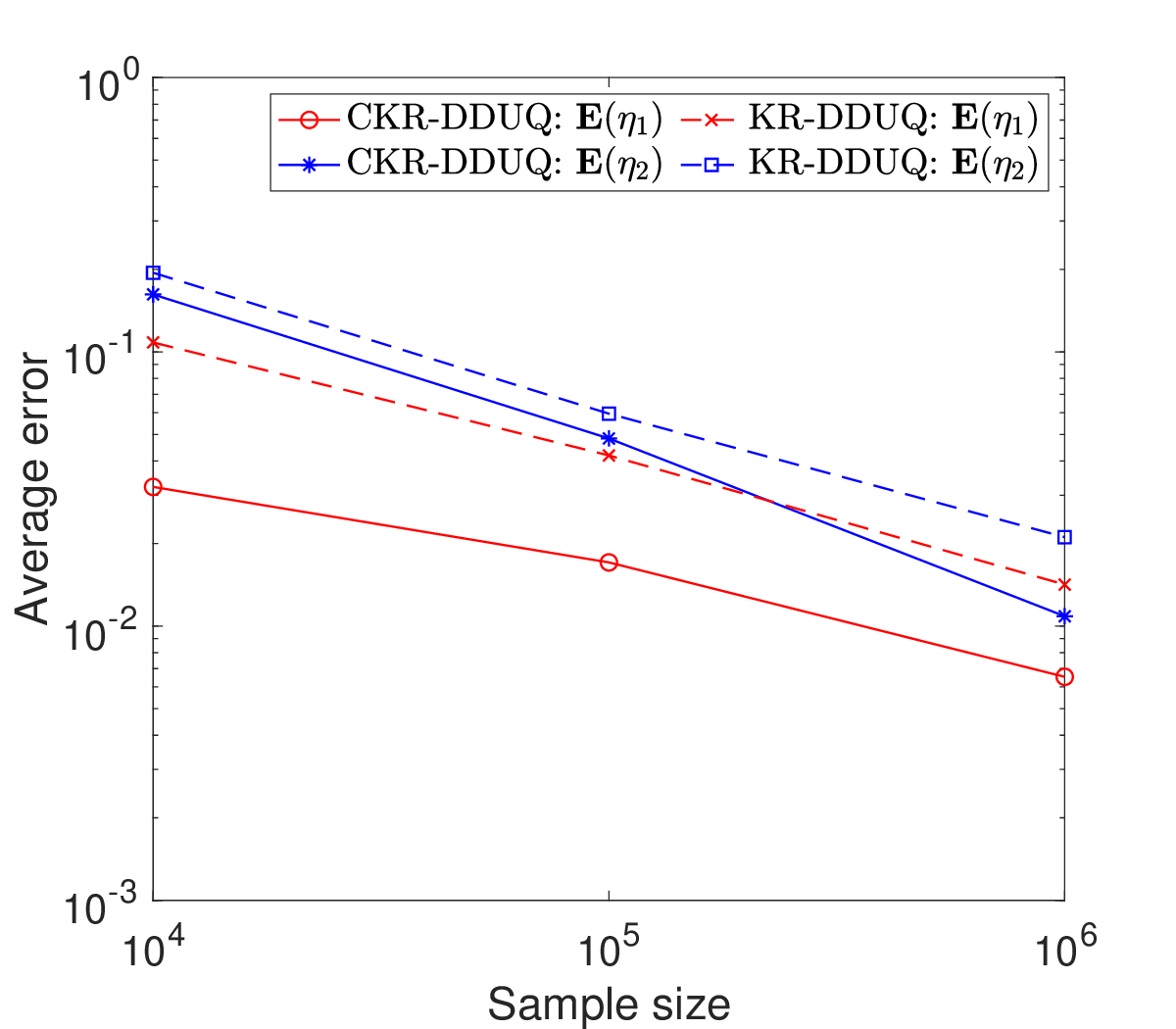}\\
    (a) Average mean errors & (b) Average variance errors
    \end{tabular}}
    \caption{\lr Average CKR-DDUQ and KR-DDUQ errors in mean and variance estimates for each output $y_i$ ($\pE(\epsilon_i)$ and $\pE(\eta_i))$, $i=1,2$, two-component diffusion test problem.}
    \label{fig_2d_error}
\end{figure}

Fixing the reference mean $\pE_{N_{\textrm{ref}}}(y^{\textrm{ref}}_i)$ and variance $\pV_{N_{\textrm{ref}}}(y^{\rm{ref}}_i)$ for each output, we repeat the CKR-DDUQ and KR-DDUQ processes thirty times for each $\Noff$ and compute the averages of $\epsilon_i$ and $\eta_i$, which are denoted by $\pE(\epsilon_i)$ and $\pE(\eta_i)$, respectively. Figure \ref{fig_2d_error} shows the average mean and variance errors of CKR-DDUQ and KR-DDUQ for this test problem. It can be seen that all errors decrease as the sample size increases for both methods, while the errors of CKR-DDUQ are always smaller than those of KR-DDUQ with the same sample size. 

To further compare the results obtained by CKR-DDUQ and KR-DDUQ, the effective sample size  \cite{kong1994sequential,doucet2000sequential,liu1996metropolized} is considered. The effective sample size is theoretically defined as the equivalent number of independent samples generated directly from the target distribution, and it is an important measure for the efficiency of importance sampling procedures. Following \cite{martino2017effective,doucet2001sequential}, we compute the effective sample size through
\begin{eqnarray}
    N_i^{\text {eff}}=\frac{\left(\sum_{s=1}^{N_{\text {off }}} w_i^{(s)}\right)^2}{\sum_{s=1}^{N_{\text {off }}}\left(w_i^{(s)}\right)^2}, \quad i=1, \ldots, M, 
\end{eqnarray}
where $w_i^{(s)}$ are the weights obtained in Algorithm \ref{alg_on1} (or \eqref{dduq_weights} for KR-DDUQ). 
Similarly to the averages of mean and variance errors discussed above, the average effective sample sizes are denoted by $\pE(N_i^{\text {eff }})$, where the CKR-DDUQ and KR-DDUQ processes are repeated thirty times and the averages  are rounded to integers. Table \ref{table_2d_ess} shows the average effective sample sizes of the two methods with respect to $\Noff=10^4,10^5$ and $10^6$. It is clear that CKR-DDUQ gives more effective samples than KR-DDUQ, which is consistent with the errors shown in Figure \ref{fig_2d_error}---CKR-DDUQ has smaller errors than KR-DDUQ. 
\begin{table}[!ht]
    \centering
    \caption{Average effective sample sizes, two-component diffusion test problem.}
    \begin{tblr}
    {
      cells = {c},
      cell{2}{1}     = {r=2}{},
      cell{4}{1}     = {r=2}{},
      hline{1-2,4,6} = {-}  {},
      hline{3,5}     = {2-5}{},
    }
               & Method   & $\Noff=10^4$  & $\Noff=10^5$ & $\Noff=10^6$ \\
    $\pE(N_1^{\text {eff }})$ & CKR-DDUQ  & 360          & 3958         & 40659        \\
                              & KR-DDUQ   & 294          & 3407         & 35682       \\
    $\pE(N_2^{\text {eff }})$ & CKR-DDUQ  & 46           & 577          & 6156         \\
                              & KR-DDUQ   & 44           & 532          & 5475         
    \end{tblr}\label{table_2d_ess}
\end{table}

\begin{figure}[!htp]
    \centerline{
    \begin{tabular}{cc}
    \includegraphics[width=0.37\textwidth]{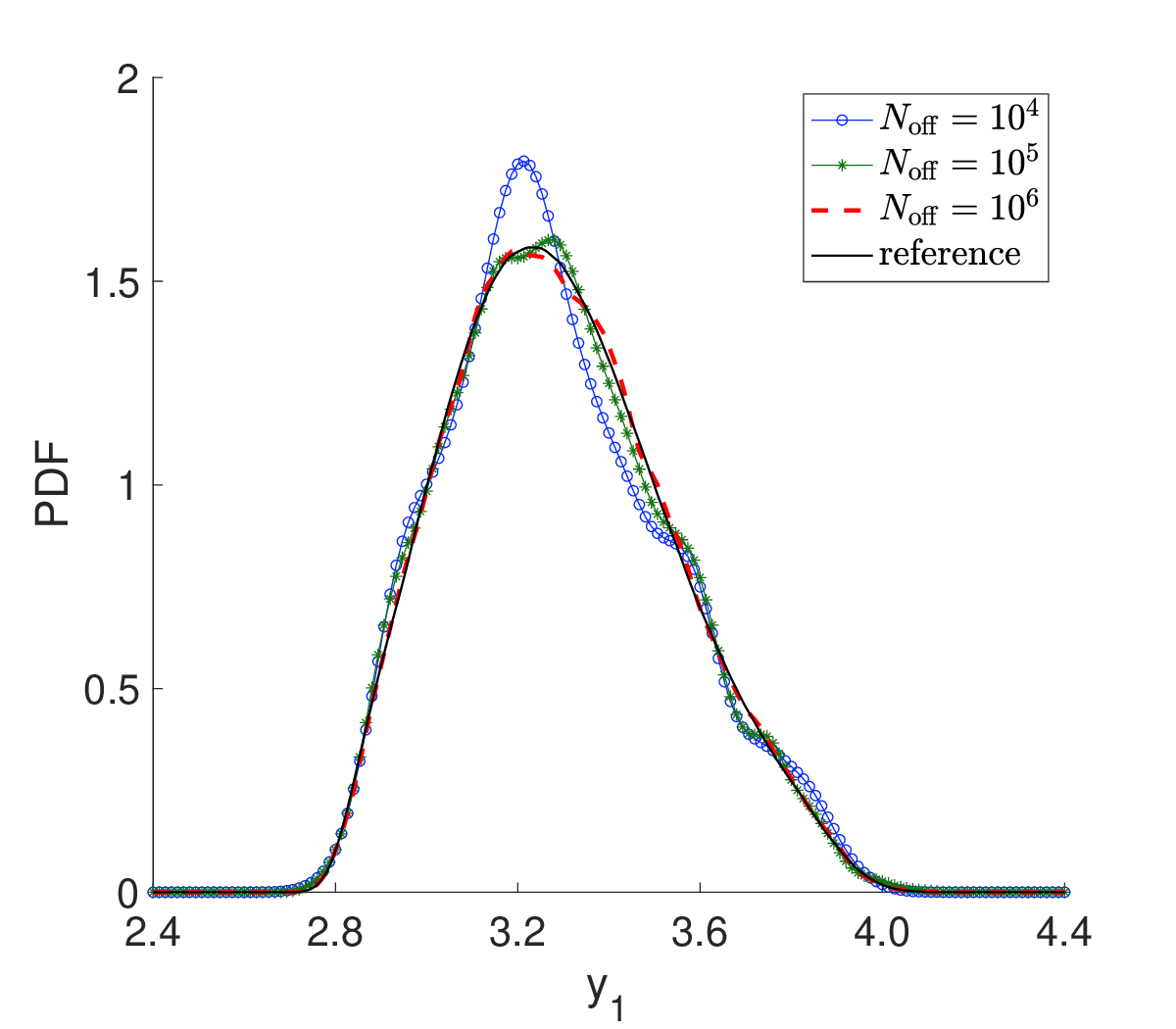}&
    \includegraphics[width=0.37\textwidth]{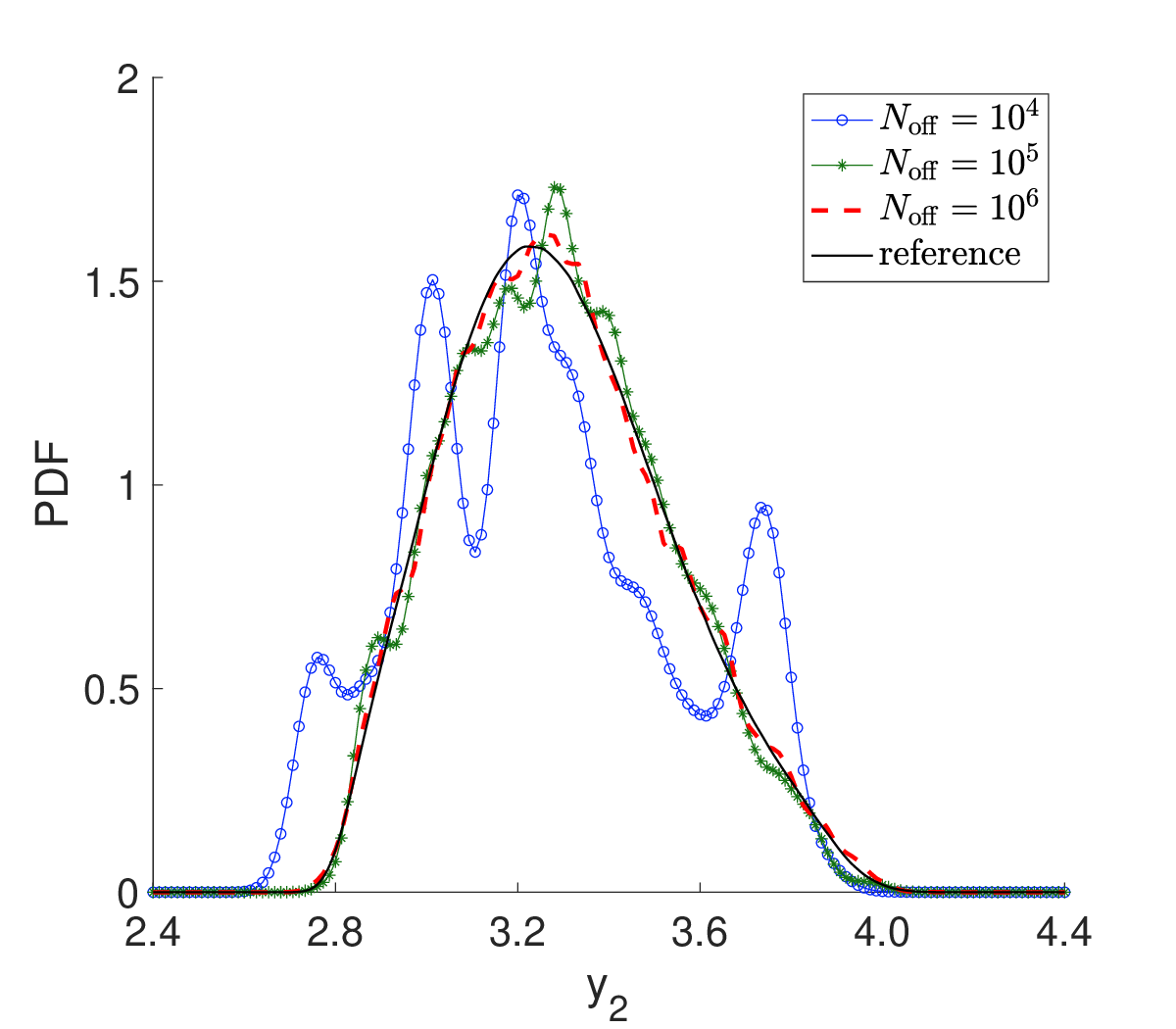}
    \end{tabular}}
    \caption{PDFs of the outputs of interest, two-component diffusion test problem.}
    \label{fig_2d_output_pdf}
\end{figure}

The PDFs of the outputs are also assessed. As each  output in our test problems is set to a scalar, kernel density estimation (KDE) \cite{hansen2008uniform} is applied to estimate the output PDFs for simplicity. For CKR-DDUQ, KDE is applied to the weighted samples $\{w^{(s)}_iy_i(u(x,\xi^{(s)}_i,\tau^{(s)}_i))\}^{\Noff}_{s=1}$ (for $i=1,\ldots,M$), and to compute the reference PDF, KDE is applied to the samples $\{y^{\rm{ref}}_i(\xi^{(s)})\}^{N_{\rm {ref}}}_{s=1}$. Figure \ref{fig_2d_output_pdf} shows these PDFs, where it is clear that as the number of samples ($\Noff$) increases, our CKR-DDUQ estimates of the PDFs approach the reference PDFs. 

\subsection{Three-component diffusion problem with \texorpdfstring{$42$}{42} parameters}\label{3domain_test}
\begin{figure}[!htp]
	\setlength{\unitlength}{1.5cm}
	\centerline{
	\begin{picture}(6,2.3)(0,0)
        \put(0,0){\line(1,0){6}}
        \put(0,0){\line(0,1){2}}
        \put(0,2){\line(1,0){6}}
        \put(4,0){\line(1,0){2}}
        \put(6,0){\line(0,1){2}}
        \put(0.8,1){$D_1$}
        \put(2.8,1){$D_2$}
        \put(4.8,1){$D_3$}
        \put(0,0){\circle*{0.1}}
        \put(-0.3,-0.25){\scriptsize (0,0)}
        \put(2,0){\circle*{0.1}}
        \put(1.8,-0.25){\scriptsize (1,0)}
        \put(4,0){\circle*{0.1}}
        \put(3.8,-0.25){\scriptsize (2,0)}
        \put(6,0){\circle*{0.1}}
        \put(5.8,-0.25){\scriptsize (3,0)}
        \put(6,2){\circle*{0.1}}
        \put(6.05,1.95){\scriptsize (3,1)}
        \put(-1.5,0.5){\vector(1,0){0.5}}
        \put(-1.5,0.5){\vector(0,1){0.5}}
        \put(-1.4,0.2){$x_1$}
        \put(-1.85,0.6){$x_2$}
        \linethickness{0.7mm}
        \put(2,0){\line(0,1){2}}
        \put(4,0){\line(0,1){2}}
	\end{picture}
	}
	\caption{Illustration of the spatial domain with three components.}
	\label{fig_3d}
\end{figure}
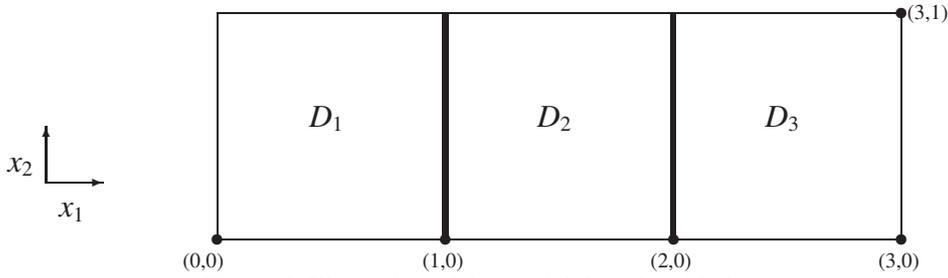

In this test problem, the diffusion equation \eqref{ht1}--\eqref{ht2} is posed on the physical domain $D=(0,3)\times(0,1)\subset \dsR^{2}$, which consists of three components as shown in Figure \ref{fig_3d}.
The homogeneous Dirichlet boundary condition is specified on the boundary $\pD$ and the source function is set to $f=100$. The permeability coefficient $a(x,\xi)$ is again set to be a random field with covariance function \eqref{covariance} and KL expansion \eqref{kl_expan} on each local subdomain. In this test problem, we set $a_{1,0}(x)=a_{3,0}(x)=3$, $a_{2,0}(x)=2$, $\sigma=0.5$ and $L_c=0.5$, and choose $\mathfrak{N}_1=\mathfrak{N}_2=\mathfrak{N}_3=14$, which results in  a total of $42$ parameters in the input vector $\xi$ with $\xi_1\in\dsR^{14}$, $\xi_2\in\dsR^{14}$ and $\xi_3\in\dsR^{14}$. The distribution of $\xi$ is set to the same as that in Section~\ref{2domain_test}.
The outputs of this test problem are defined as 
\begin{eqnarray}
    y_1(\xi)=\int_{Z_1} u(x,\xi) \,{\rm d} x_2,\\
    y_2(\xi)=\int_{Z_2} u(x,\xi) \,{\rm d} x_1,\\
    y_3(\xi)=\int_{Z_3} u(x,\xi) \,{\rm d} x_2,
\end{eqnarray}
where $x=[x_1,x_2]^{\top}$, $Z_1:=\{x|\> x_1=0.5,\>0\leq x_2 \leq 1\}$, $Z_2:=\{x|\> 1\leq x_1 \leq 2,\> x_2=0.5\}$ and $Z_3:=\{x|\> x_1=2.5,\> 0\leq x_2 \leq 1\}$.

The parallel Dirichlet-Neumann domain decomposition method and the bilinear finite element method with a mesh size $h=1/16$ are also utilized in this test problem. One hundred samples are used to generate the POD bases of the interface functions and the interface input parameters in the reduced dimensional representation are $\tau_1 \in\dsR^{6}$, $\tau_2 \in\dsR^{4}$, and $\tau_3 \in\dsR^{6}$. Then the local inputs for the three subdomains are $[\xi_1^{\top},\tau_1^{\top}]^{\top}\in\dsR^{20}$, $[\xi_2^{\top},\tau_2^{\top}]^{\top}\in\dsR^{18}$ and $[\xi_3^{\top},\tau_3^{\top}]^{\top}\in\dsR^{20}$, respectively. 

\begin{figure}[!ht]
    \centerline{
    \begin{tabular}{ccc}
    \includegraphics[width=0.32\textwidth]{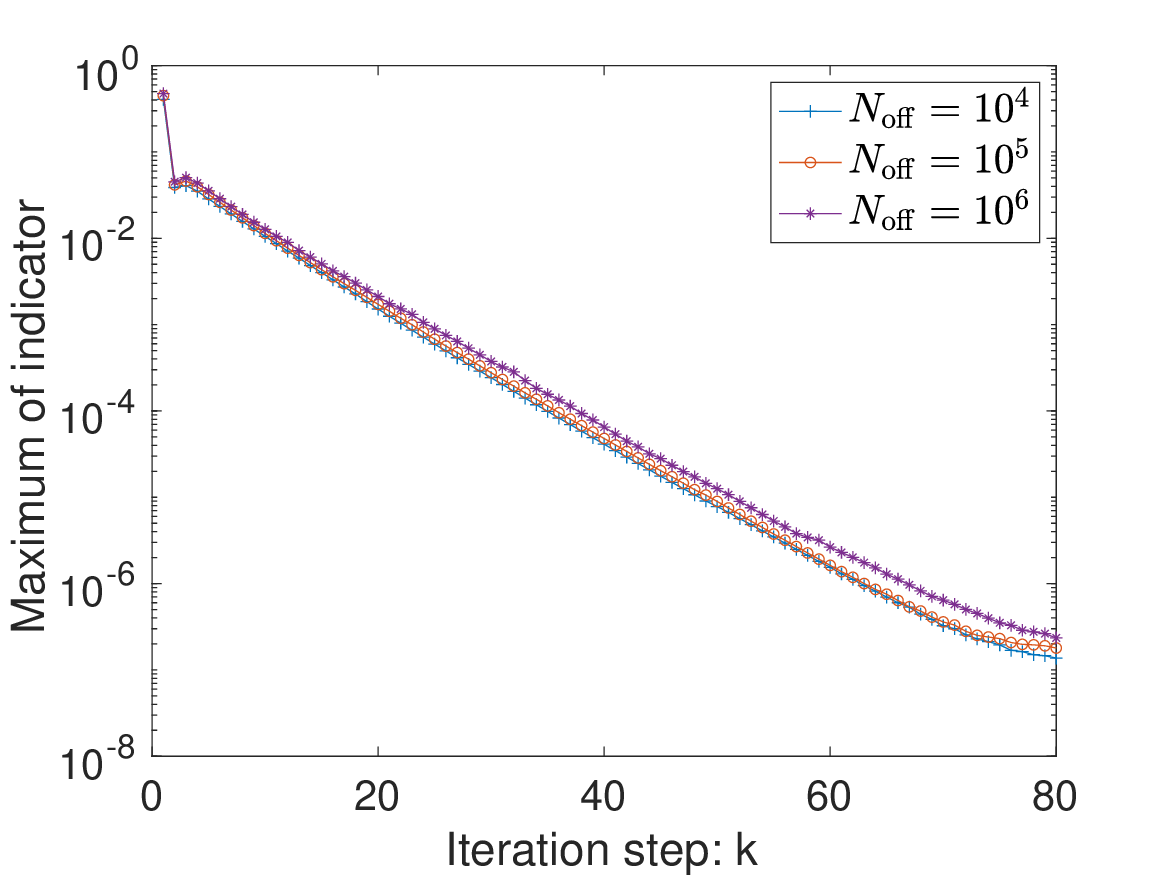}&
    \includegraphics[width=0.32\textwidth]{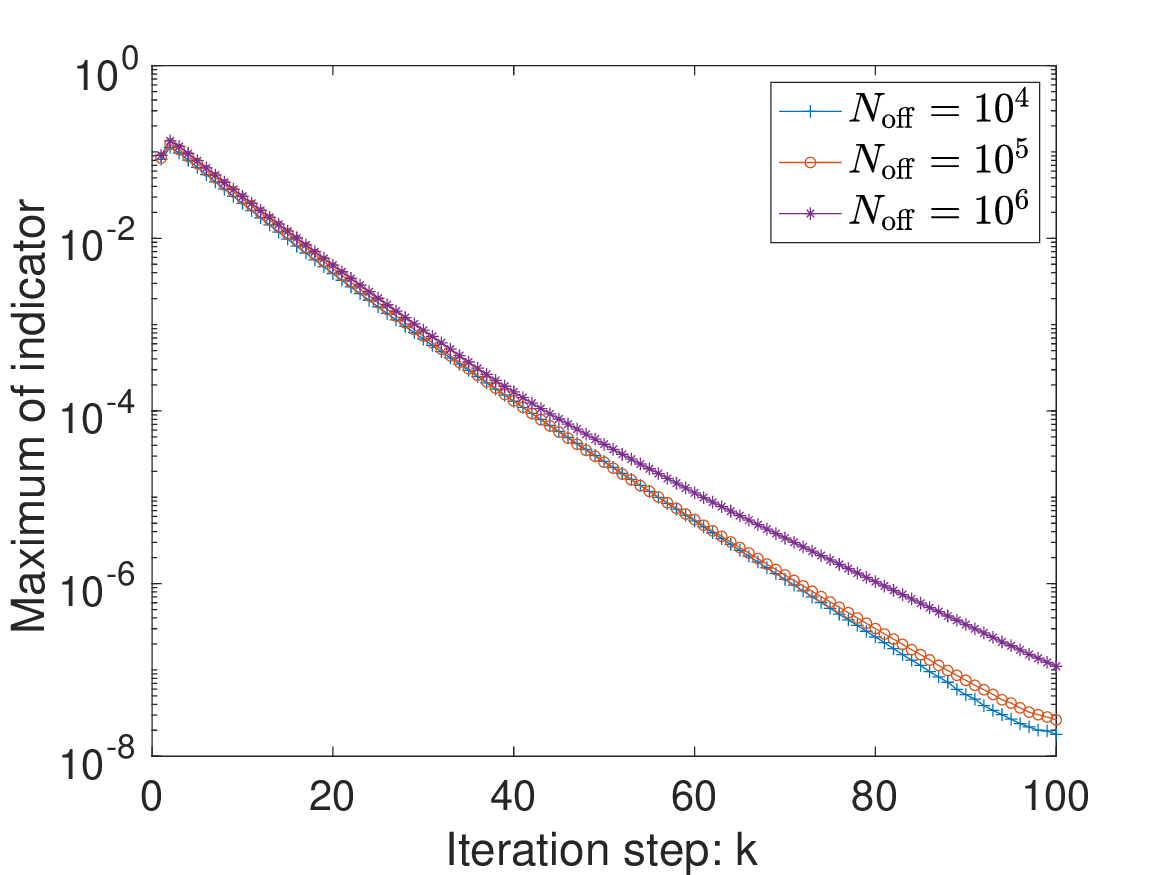}&
    \includegraphics[width=0.32\textwidth]{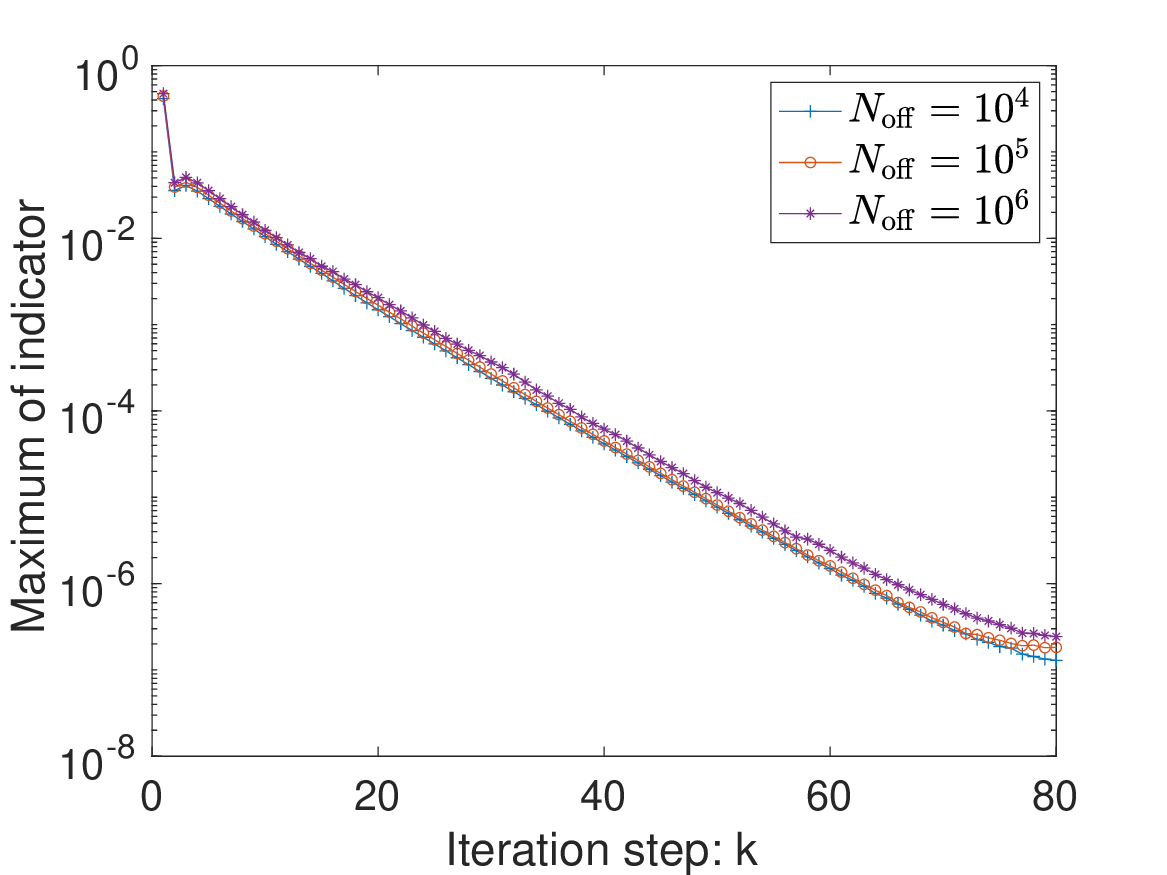}\\
    (a) Error indicator on $D_1$ & (b) Error indicator on $D_2$&
    (c) Error indicator on $D_3$ \\
    \end{tabular}}
    \caption{\lr
    Maximum of the error indicator on subdomain $D_1$ ($\max_{s=1:{\Noff}}\|\tau^{k+1}_1(\xi^{(s)})-\tau^k_1(\xi^{(s)})\|_{\infty}$), subdomain $D_2$ ($\max_{s=1:{\Noff}}\|\tau^{k+1}_2(\xi^{(s)})-\tau^k_2(\xi^{(s)})\|_{\infty}$) and that on subdomain $D_3$ ($\max_{s=1:{\Noff}}\|\tau^{k+1}_3(\xi^{(s)})-\tau^k_3(\xi^{(s)})\|_{\infty}$) for the coupling surrogates, three-component diffusion test problem.}
    \label{fig_3d_surrogate_convergence}
\end{figure}

In the offline stage, each proposal PDF $p_{\tau_i}(\tau_i)$ is still a multivariate Gaussian distribution, where the mean vectors and the covariance matrices are estimated with the snapshots obtained in the POD procedure. For constructing the coupling surrogates, we employ ResNets consisting of ten hidden layers, where each hidden layer contains sixty four neurons and the activation function is the LeakyReLU function. In the online stage, $\Non$ is set to equal $\Noff$ and the threshold of domain decomposition iteration is $tol=10^{-6}$. Figure \ref{fig_3d_surrogate_convergence} shows the maximum of the error indicator (introduced in Section \ref{2domain_test}) with respect to domain decomposition iterations using these coupling surrogates for $\Noff=10^4,10^5$ and $10^6$, where it is clear that the maximum of the error indicator reduces exponentially as the iteration step $k$ increases for each subdomain.

\begin{figure}[!htp]
    \centerline{
    \begin{tabular}{cc}
    \includegraphics[width=0.37\textwidth]{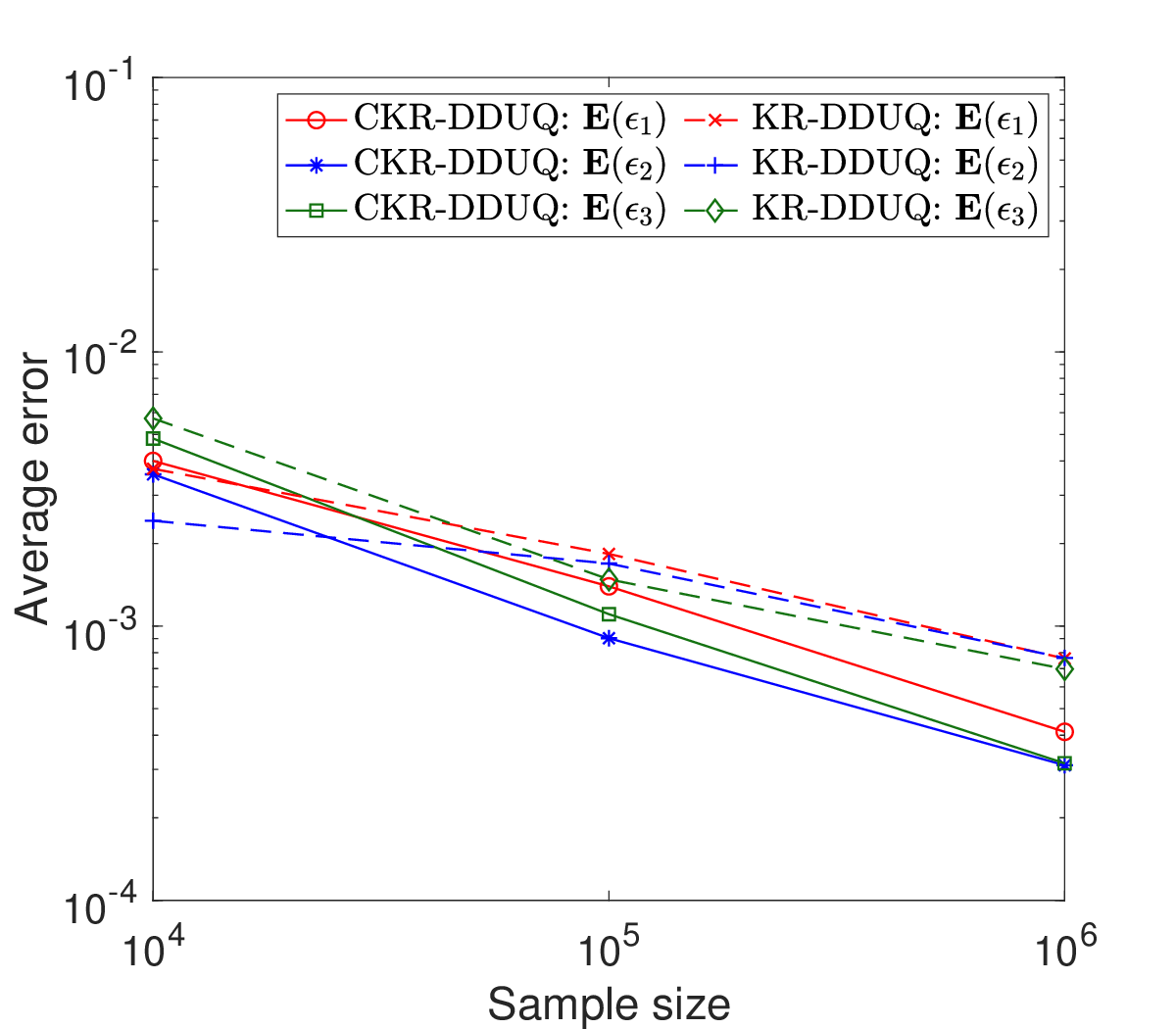}&
    \includegraphics[width=0.37\textwidth]{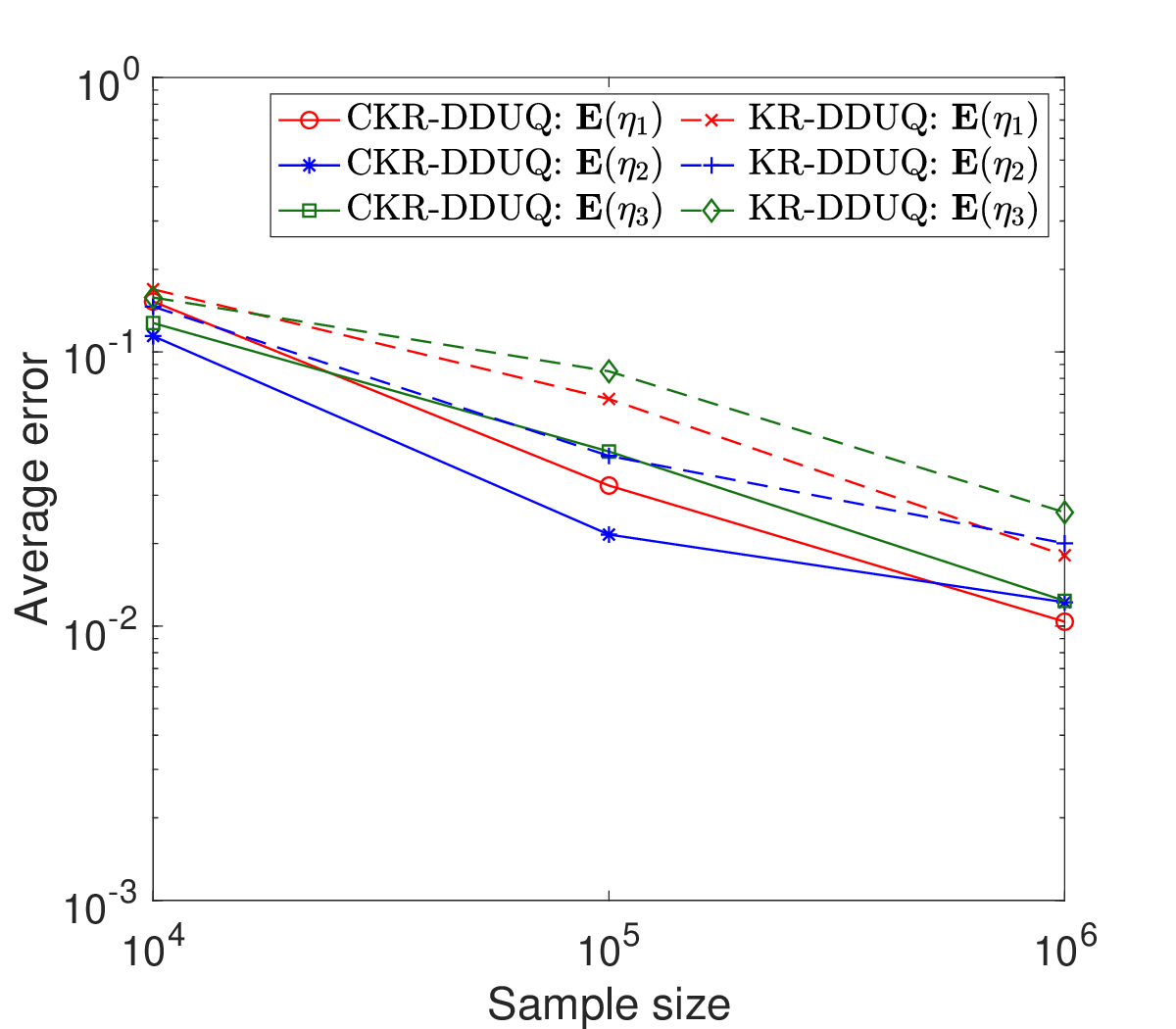}\\
    (a) Average mean errors & (b) Average variance errors
    \end{tabular}}
    \caption{\lr Average CKR-DDUQ and KR-DDUQ errors in mean and variance estimates for each output $y_i$ ($\pE(\epsilon_i)$ and $\pE(\eta_i))$, $i=1,2,3$, three-component diffusion test problem.}
    \label{fig_3d_error}
\end{figure}

\begin{figure}[!ht]
    \centerline{
    \begin{tabular}{ccc}
    \includegraphics[width=0.32\textwidth]{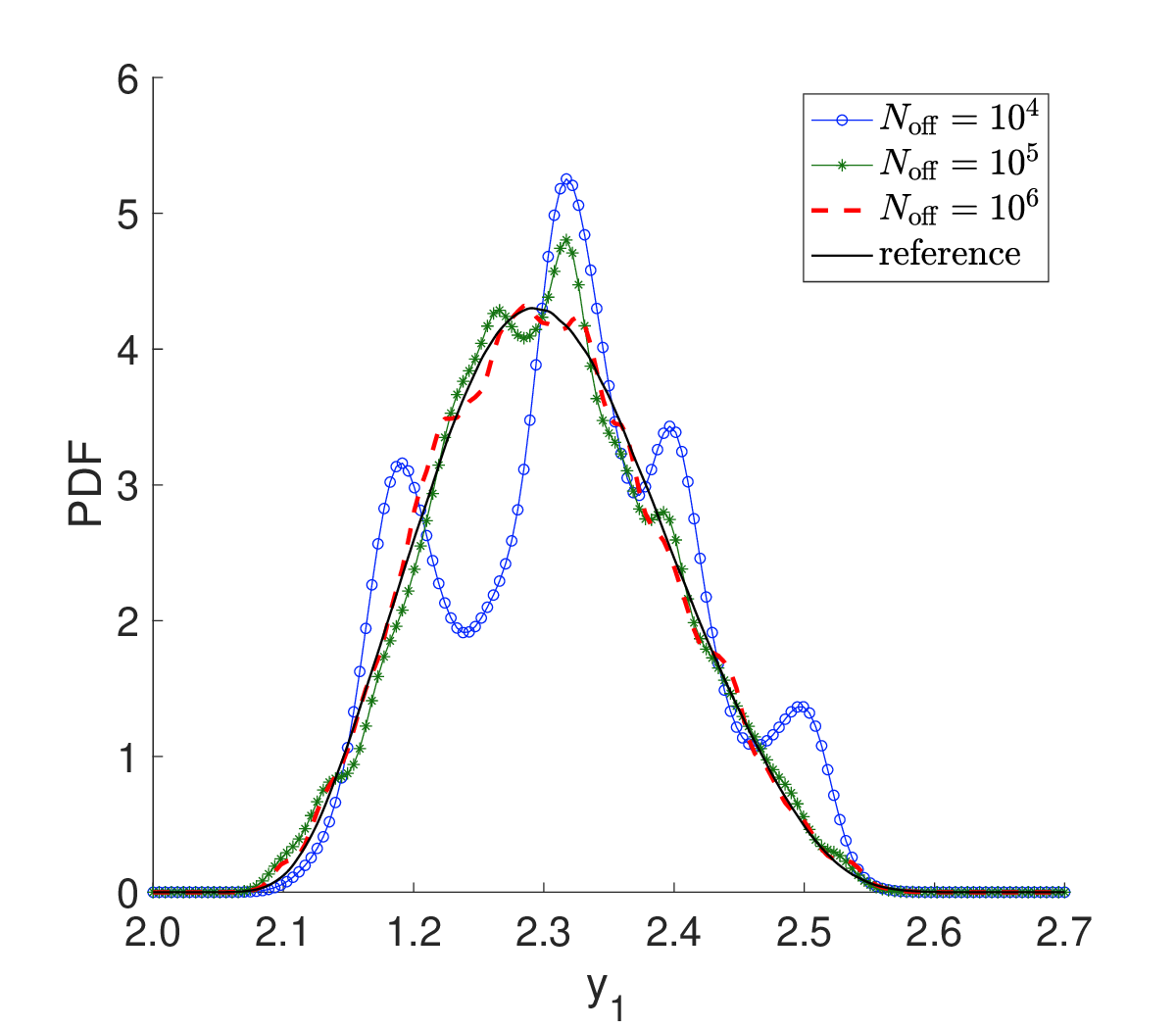}&
    \includegraphics[width=0.32\textwidth]{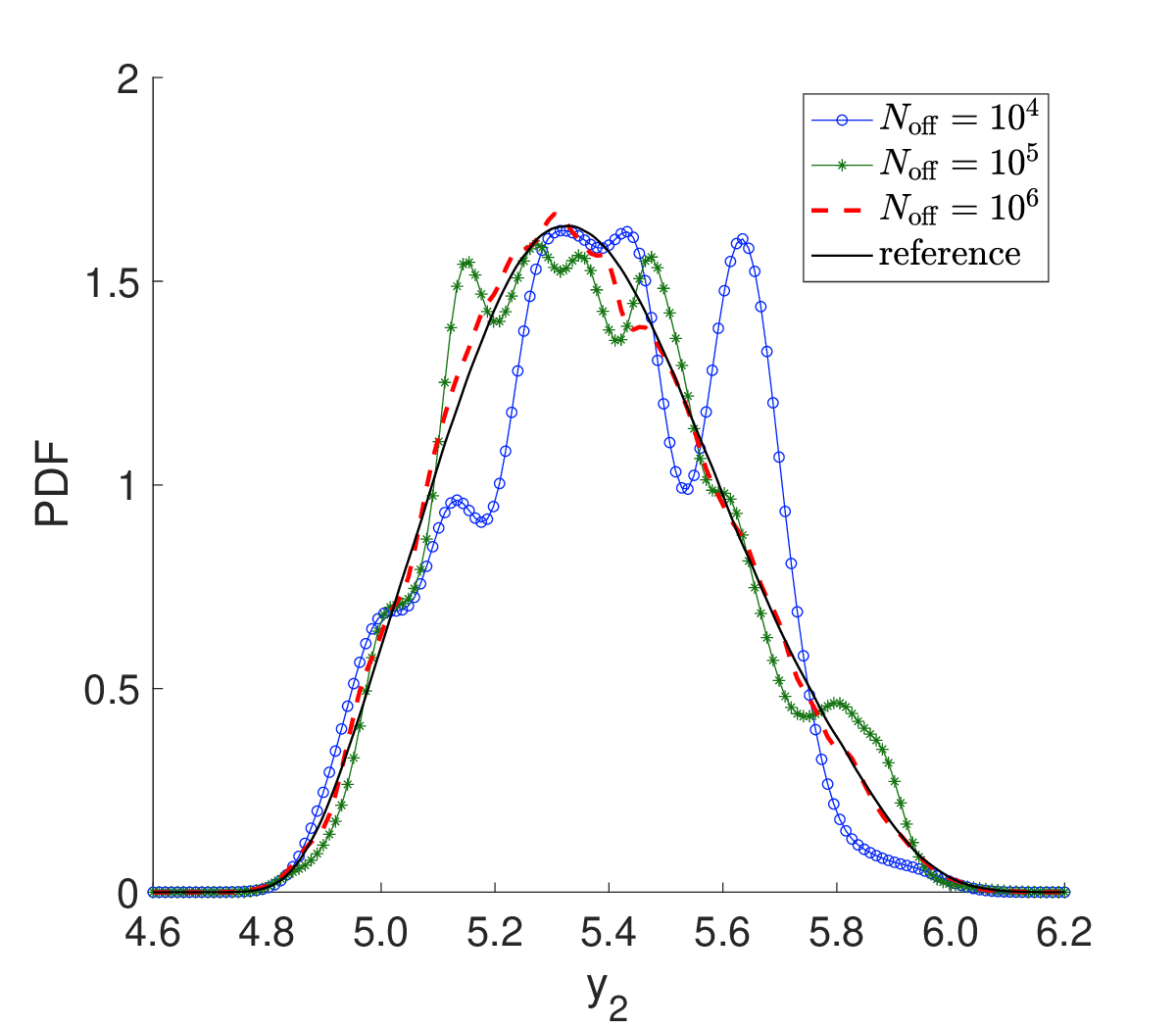}&
    \includegraphics[width=0.32\textwidth]{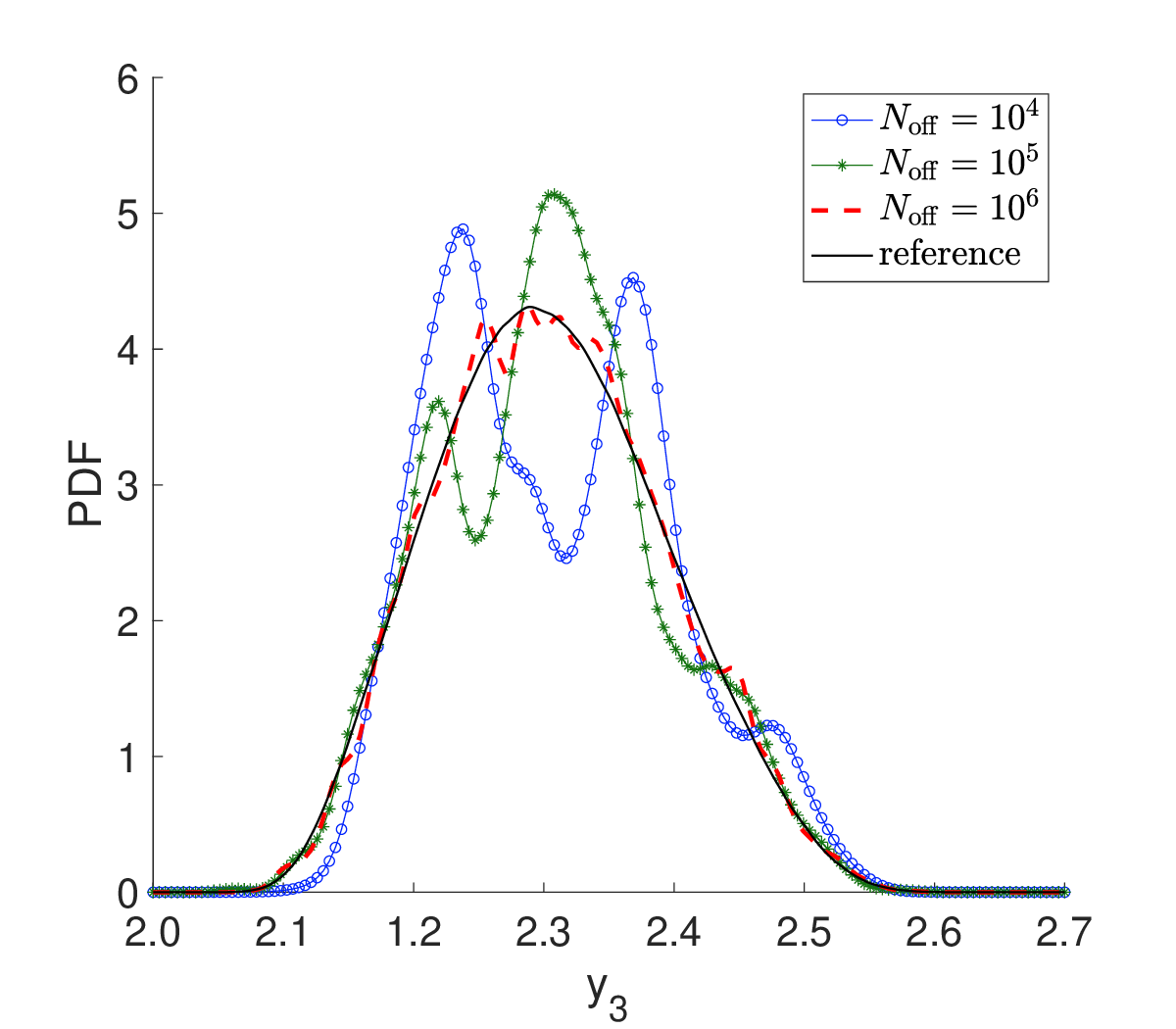}
    \end{tabular}}
    \caption{PDFs of the outputs of interest, three-component diffusion test problem.}
    \label{fig_3d_output_pdf}
\end{figure}

In the density estimation step, both cKRnet and KRnet are set to $R=3$ for subdomains $D_1$ and $D_3$, and $R=4$ for subdomain $D_2$. In addition, we employ $L=4$ affine coupling layers for both cKRnet and KRnet, and each affine coupling layer incorporates a fully connected network with two hidden layers, each consisting of thirty two neurons, and employs the ReLU activation function. Moreover, the Adam optimizer with default settings is employed to train cKRnet and KRnet, and the learning rate for Adam optimizer is set to $0.001$. A set of reference output samples are also generated by solving the problem globally with a sample size of $N_{\rm ref}=10^7$. Then, the average mean and variance errors of CKR-DDUQ and KR-DDUQ computed by thirty repeats are shown in Figure \ref{fig_3d_error}(a) and Figure \ref{fig_3d_error}(b), which demonstrate that the errors decrease as the sample size increases and CKR-DDUQ gives more accurate results than KR-DDUQ with the same sample size. Table \ref{table_3d_ess} presents the average effective sample sizes of the two methods for different $\Noff$. Again, the average effective sample sizes of CKR-DDUQ are larger than those of KR-DDUQ. The  PDFs of outputs in this test problem are shown in Figure \ref{fig_3d_output_pdf}, where it is clear that as $\Noff$ increases, the PDFs generated by CKR-DDUQ approach the reference PDFs.

\begin{table}[!ht]
    \centering
    \caption{Average effective sample sizes, three-component diffusion test problem.}
    \begin{tblr}
    {
      cells = {c},
      cell{2}{1} = {r=2}{},
      cell{4}{1} = {r=2}{},
      cell{6}{1} = {r=2}{},
      hline{1-2,4,6,8} = {-}{},
      hline{3,5,7} = {2-5}{},
    }
             & Method   & $\Noff=10^4$ & $\Noff=10^5$ & $\Noff=10^6$ \\
    $\pE(N_1^{\text {eff }})$ & CKR-DDUQ & 53           & 672          & 7184  \\
                              & KR-DDUQ  & 50           & 597          & 6401  \\
    $\pE(N_2^{\text {eff }})$ & CKR-DDUQ & 88           & 979          & 10934 \\
                              & KR-DDUQ  & 83           & 907          & 10159 \\
    $\pE(N_3^{\text {eff }})$ & CKR-DDUQ & 33           & 451          & 4834  \\
                              & KR-DDUQ  & 31           & 409          & 4359         
               
    \end{tblr}\label{table_3d_ess}
\end{table}

\subsection{Random domain decomposition diffusion problem with \texorpdfstring{$25$}{25} parameters}\label{rdd_test} 
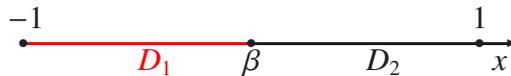
\begin{figure}[!htp]
    \setlength{\unitlength}{1cm}
    \centerline{
    \begin{picture}(6,1)(-0.2,-0.2)
    \put(-1,0){\color{red}\line(1,0){3}}
    \put(2,0){\line(1,0){3}}
    \put(2,0){\circle*{0.1}}
    \put(1.9,-0.35){$\beta$}
    \put(0.5,-0.35){$\color{red} D_1$}
    \put(3.5,-0.35){$D_2$}
    \put(-1,0){\circle*{0.1}}
    \put(-1.2,0.2){$-1$}
    \put(5,0){\circle*{0.1}}
    \put(4.9,0.2){$1$}
    \put(5,0){\vector(1,0){.5}}
    \put(5.15,-0.35){$x$}
    \end{picture}
    }
    \caption{Illustration of the spatial domain with a random interface $\beta$. }
    \label{rdd_domain}
\end{figure}

Following \cite{wintar02}, this test problem focuses on random domain decomposition. A one-dimensional diffusion problem is considered, of which the governing equation is 
\begin{eqnarray}    
    -\frac{\mathrm{d}}{\mathrm{d}x}\left(a\left(x,\xi\,\right)\frac{\mathrm{d} u\left(x,\xi\,\right)}{ \mathrm{d}x}\right)=f\left(x,\xi\,\right)& \quad \textrm{in}& D\times\Gamma, \label{ht3} \\    
   u\left(x,\xi\,\right)=0
    & \quad\textrm{on}& \partial D \times\Gamma.\label{ht4}
\end{eqnarray}
and the physical domain is shown in Figure \ref{rdd_domain}, where the location of the interface is described by a random variable $\beta$. Here, we set $\beta$ to be uniformly distributed over $[-0.2, 0.2]$. A homogeneous Dirichlet boundary condition $u=0$ is specified on $x=-1$ and $x=1$, and the source function is set to $f=1$.
Following the test problems in \cite{xiuhes05}, the permeability coefficient is set to be
\begin{eqnarray}
    a(x,\xi\,)|_{D_i}=a_{i,0}(x)+\varsigma_i\sum_{m=1}^{\mathfrak{N}_i}\frac{1}{m^2\pi^2} \cos(2\pi m x)\xi_{i,m}, \quad i=1,2,
    \label{sim_kl_expan}
\end{eqnarray}
where $a_{1,0}(x)=5$, $a_{2,0}(x)=1$, $\varsigma_1=17$, $\varsigma_2=3$, and 
$\xi_{1,m}$ and $\xi_{2,m}$ are independent uniformly distributed random variables with ranges $\xi_{1,m}\in [2, 3]$ for $m=1,\ldots,\mathfrak{N}_1$ and $\xi_{2,m}\in [0.4, 0.6]$ for $m=1,\ldots,\mathfrak{N}_2$. The dimensionality of the random parameters are set to $\mathfrak{N}_1=\mathfrak{N}_2=12$. The global solution now depends on $\xi=[\xi_{1,1},\ldots,\xi_{1,12},\xi_{2,1},\ldots,\xi_{2,12}]$ and $\beta$. The outputs of interest in this test problem are defined by
\begin{eqnarray}
    y_1(\xi,\beta\,):=u(-0.1,\xi,\beta\,), \quad y_2(\xi,\beta\,):=u(0.1,\xi,\beta\,).
    \label{ground-water-out}
\end{eqnarray}

For each subdomain $D_i\, (i=1,2)$, the local solution is denoted by $u(x, \xi_i,\beta,\tau_i)$, where $\tau_i$ is the interface parameter defined in Section~\ref{section_pre}. Since the range of $\beta$ covers the points $x=-0.1$ and $x=0.1$
where the outputs are evaluated, each local solution has the possibility to provide both outputs, and the outputs of local systems are denoted by
\begin{eqnarray}
y_{1,1}&=&u(-0.1,\xi_1,\beta,\tau_1),\\
y_{1,2}&=&u(0.1,\xi_1,\beta,\tau_1),\\
y_{2,1}&=&u(-0.1,\xi_2,\beta,\tau_2),\\
y_{2,2}&=&u(0.1,\xi_2,\beta,\tau_2).
\end{eqnarray}
For each subdomain $D_i$, the offline samples of each $y_{i,j}$, $j=1,2$ are obtained based on the conditions of $\beta$ (details are discussed in \cite{liao2015domain}) and are re-weighted by weights from Algorithm \ref{alg_on1}. Then, the  mean, the variance and the PDF of each output $y_{i,j}$, $i,j=1,2$ are estimated using the weighted samples. After that, following \cite{wintar00} (also discussed in \cite{liao2015domain}), the estimated mean, variance and PDF of each $y_i$ for $i=1,2$, can be obtained through the combination of contributions from each $y_{i,j}$. 

\begin{figure}[!htp]
    \centerline{           
    \begin{tabular}{cc}
    \includegraphics[width=0.37\textwidth]{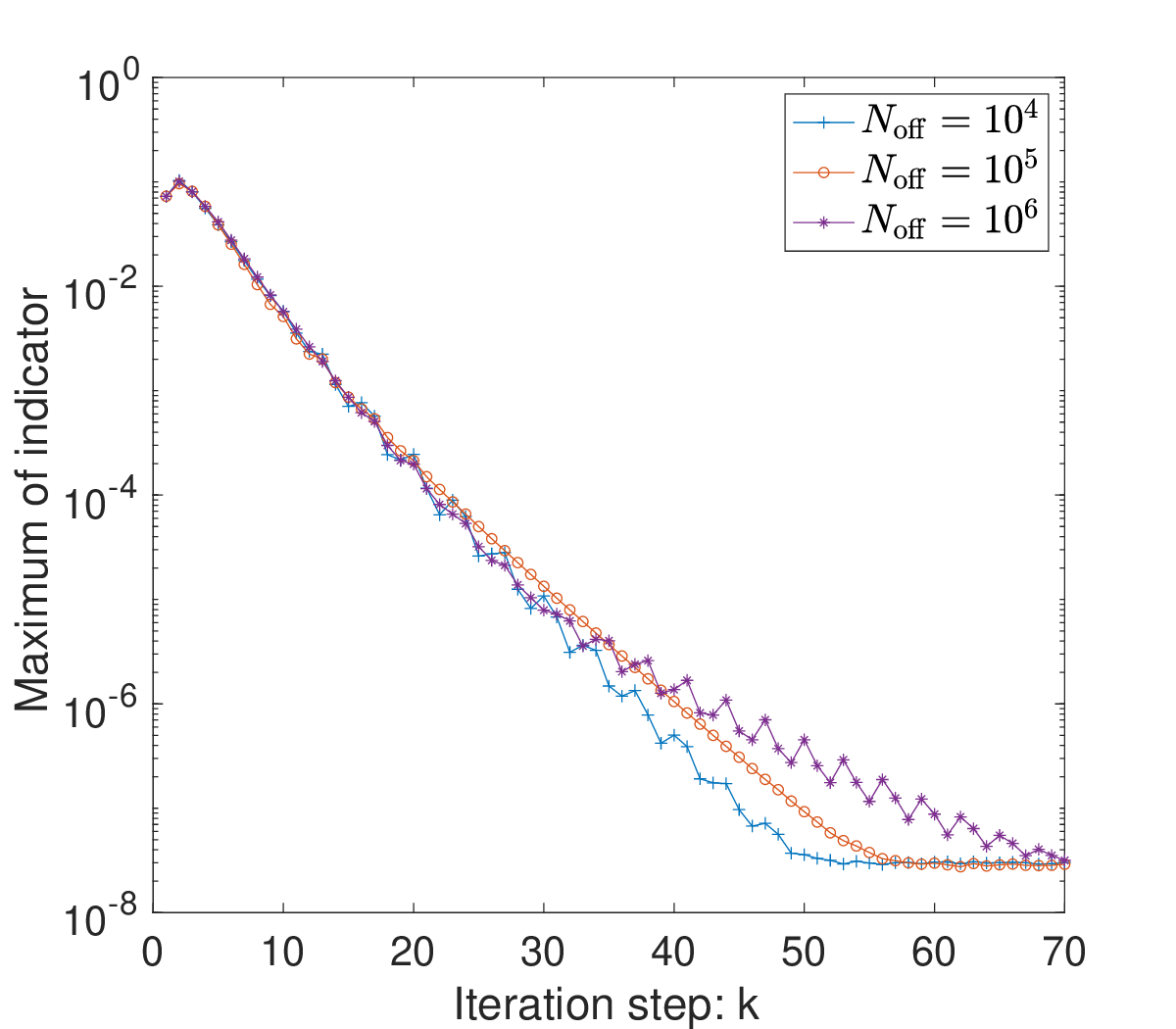}&
    \includegraphics[width=0.37\textwidth]{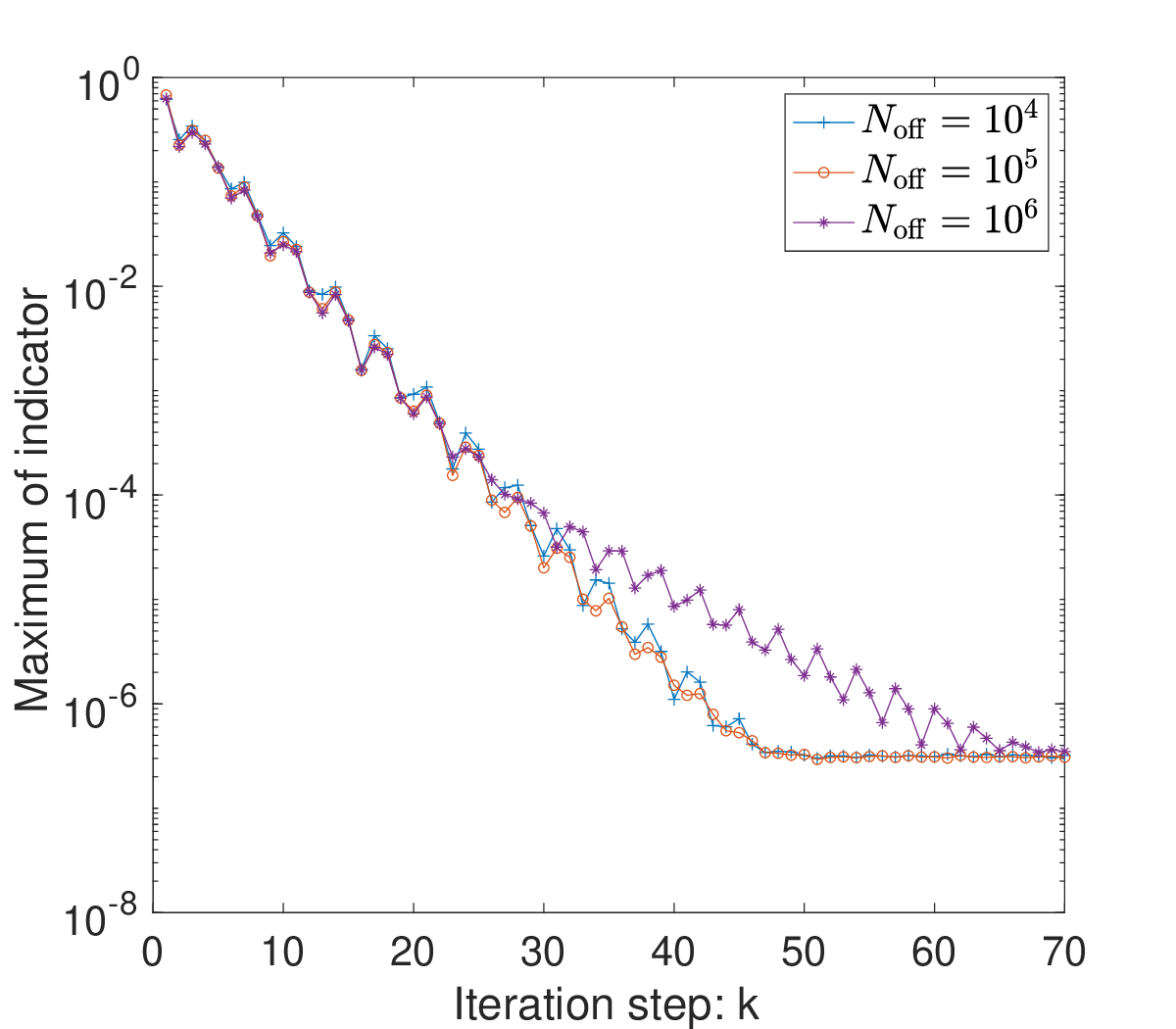}\\
    (a) Error indicator on $D_1$ & (b) Error indicator on $D_2$
    \end{tabular}}
    \caption{\lr    
    Maximum of the error indicator on subdomain $D_1$ ($\max_{s=1:{\Noff}}\|\tau^{k+1}_1(\xi^{(s)})-\tau^k_1(\xi^{(s)})\|_{\infty}$) and that on subdomain $D_2$ ($\max_{s=1:{\Noff}}\|\tau^{k+1}_2(\xi^{(s)})-\tau^k_2(\xi^{(s)})\|_{\infty}$) for the coupling surrogates, random domain decomposition test problem.}
    \label{fig_rdd_surrogate_convergence}
\end{figure}

The parallel Dirichlet-Neumann domain decomposition method is used to decompose the global problem, and a linear finite element approximation with 101 grid points is used to solve each local problem. In the offline stage, each proposal PDF $p_{\tau_i}(\tau_i)$ is a Gaussian distribution, where the mean and variance are estimated with ten samples generated by directly solving the deterministic version of \eqref{ht3}--\eqref{ht4} using the domain decomposition method. We continue to employ ResNets to construct the coupling surrogates. These ResNets consist of ten hidden layers, where each layer contains sixty four neurons and the LeakyReLU activation function is used. In the online stage, $\Non$ is set to equal $\Noff$ and the threshold of domain decomposition iteration is $tol=10^{-6}$. From Figure \ref{fig_rdd_surrogate_convergence}, the maximum of the error indicator of the domain decomposition iteration using these surrogates decreases exponentially as the number of iteration increases. 

\begin{figure}[!htp]
    \centerline{
    \begin{tabular}{cc}
    \includegraphics[width=0.37\textwidth]{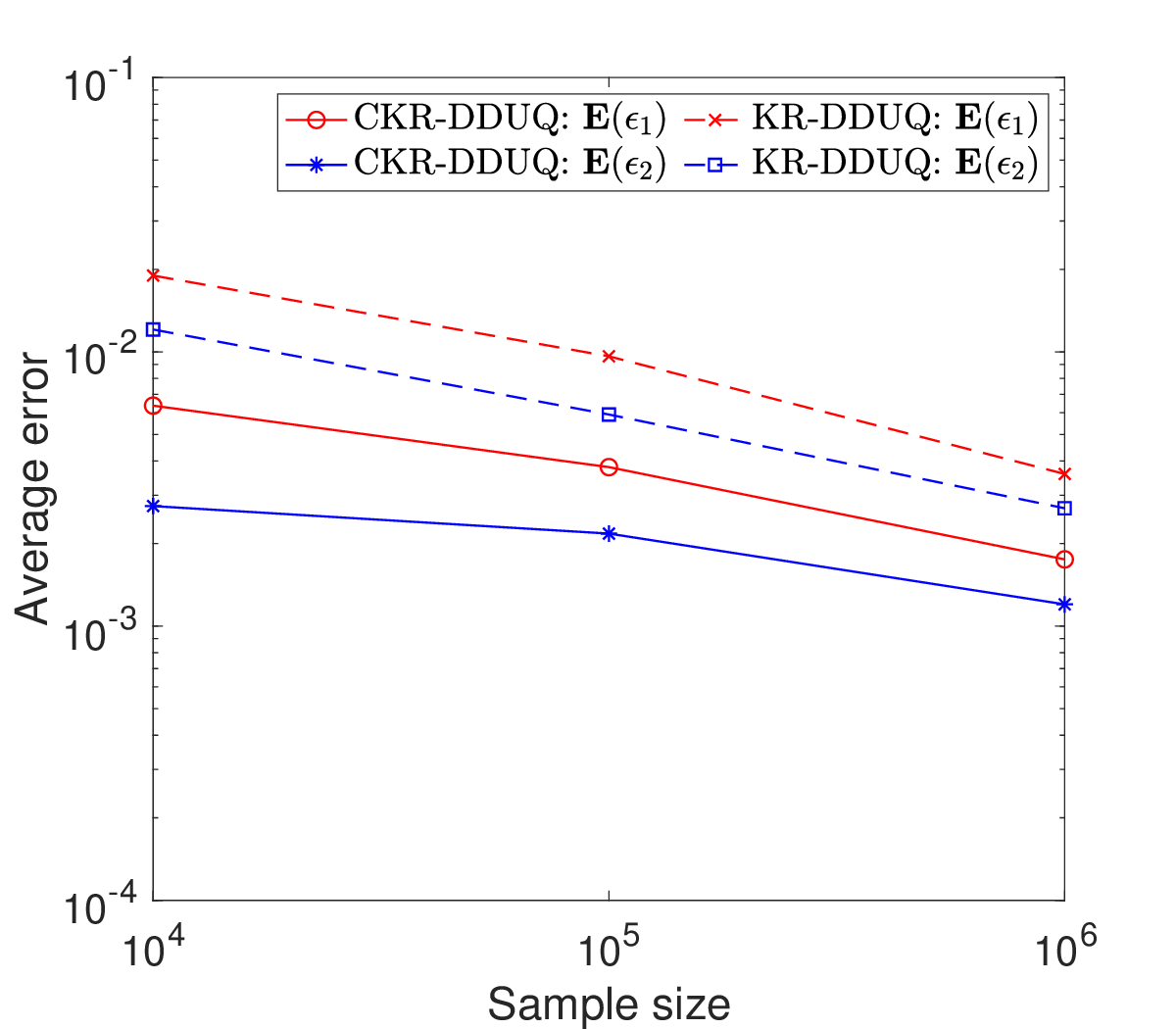}&
    \includegraphics[width=0.37\textwidth]{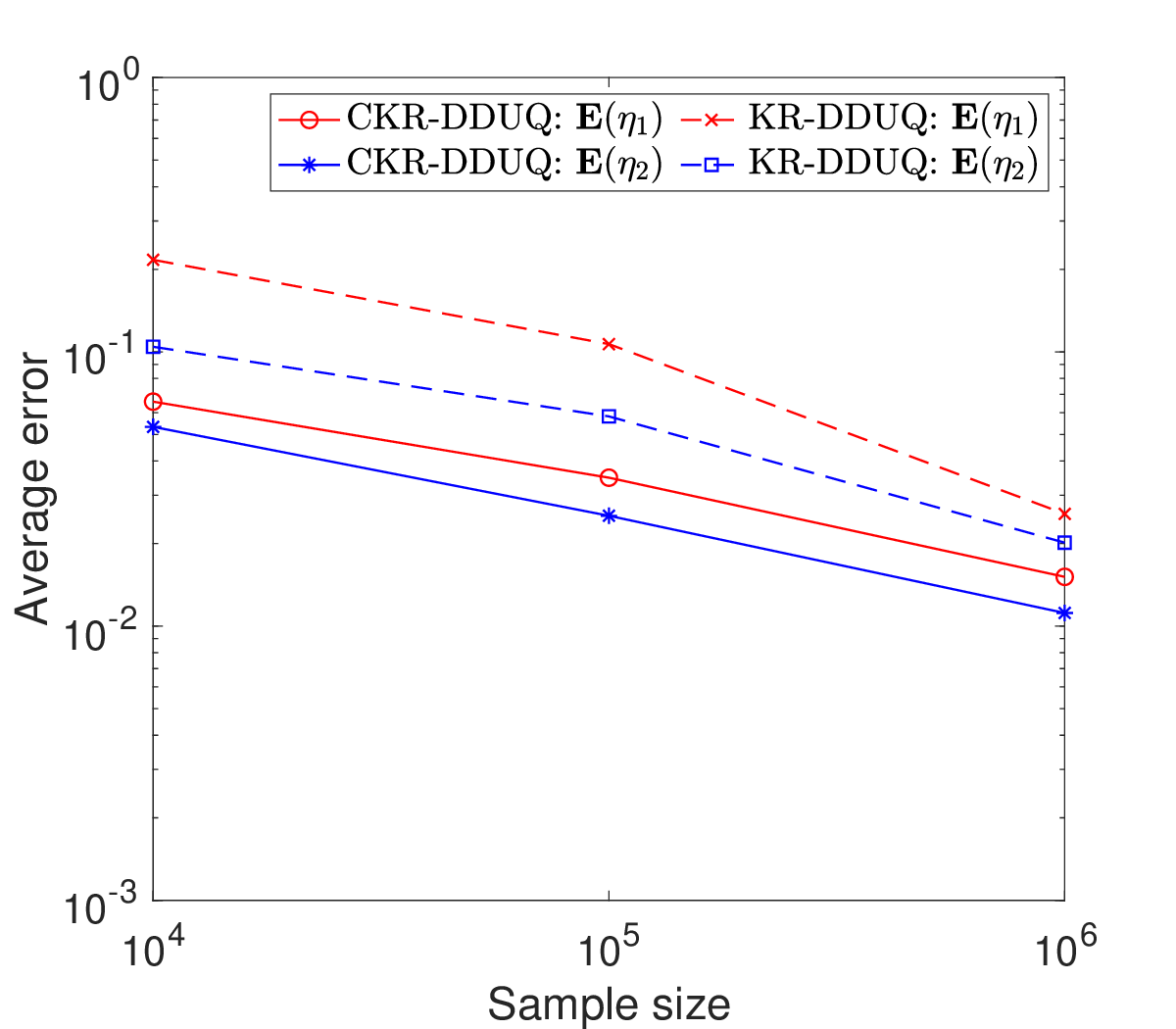}\\
    (a) Average mean errors & (b) Average variance errors
    \end{tabular}}
    \caption{\lr Average CKR-DDUQ and KR-DDUQ errors in mean and variance estimates for each output $y_i$ ($\pE(\epsilon_i)$ and $\pE(\eta_i))$, $i=1,2$, random domain decomposition test problem.}
    \label{fig_rdd_error}
\end{figure}

\begin{figure}[!ht]
    \centerline{
    \begin{tabular}{cc}
    \includegraphics[width=0.37\textwidth]{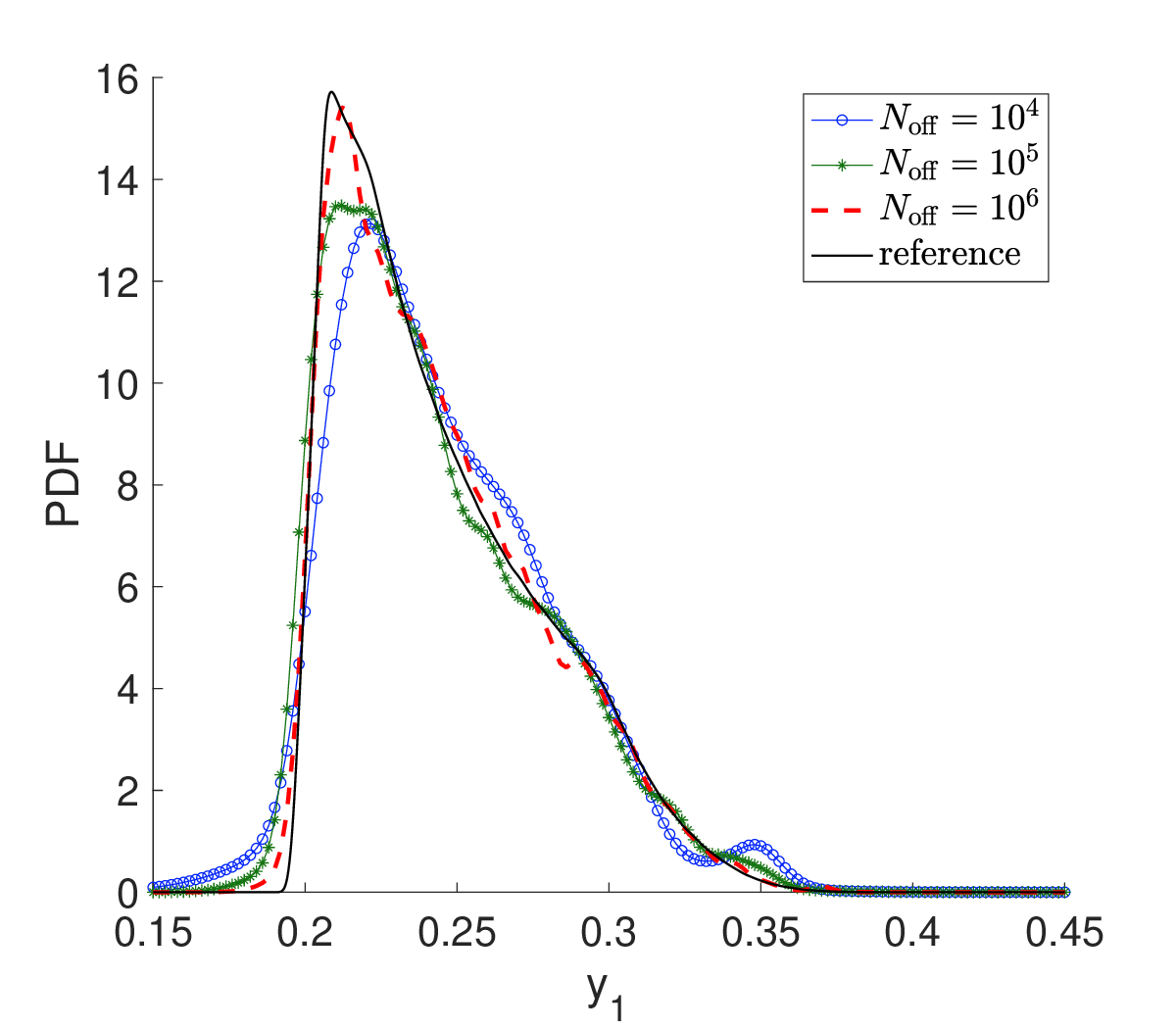}&
    \includegraphics[width=0.37\textwidth]{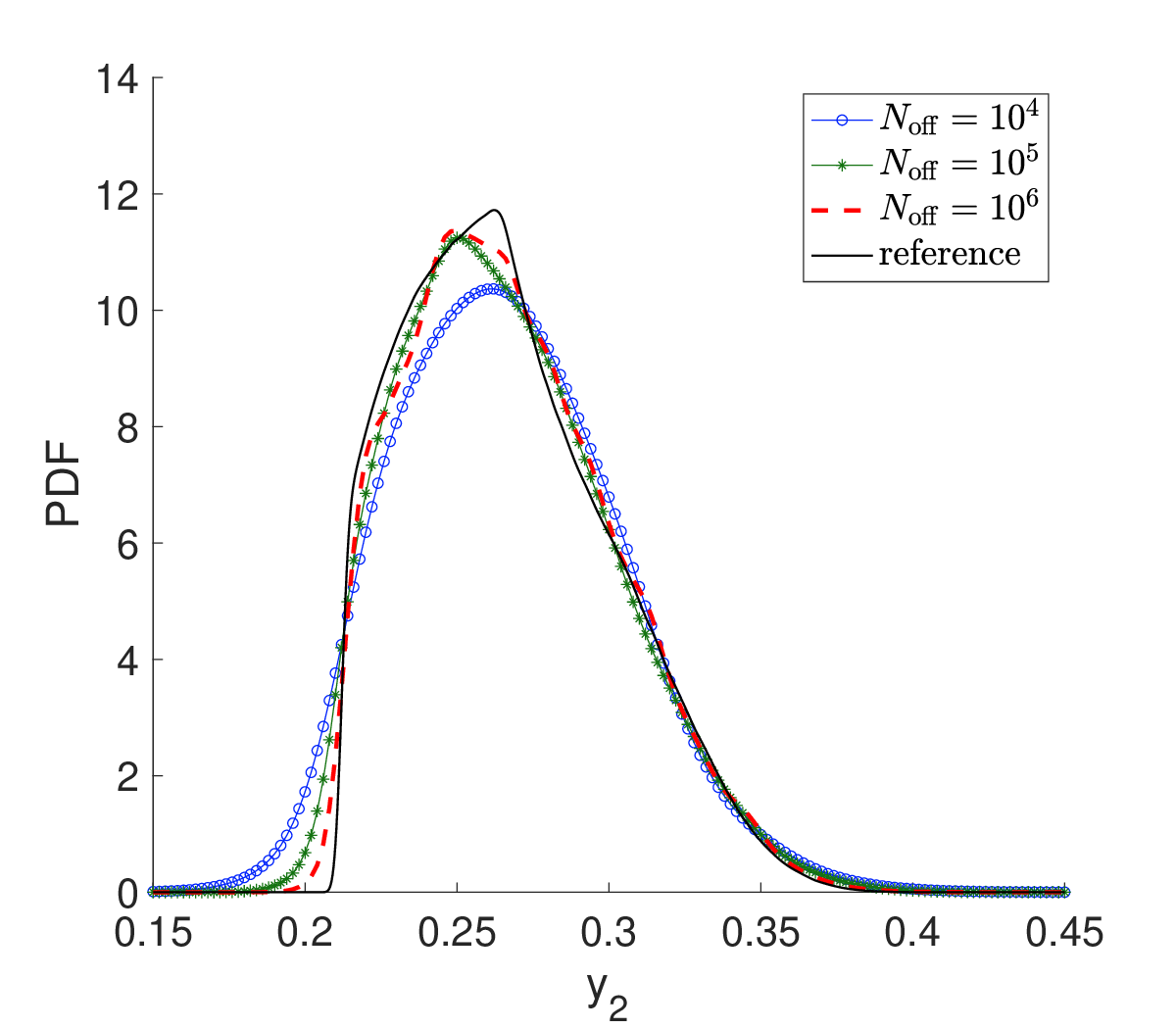}
    \end{tabular}}
    \caption{PDFs of the outputs of interest, random domain decomposition test problem.}
    \label{fig_rdd_output_pdf}
\end{figure}

For density estimation in the online stage,  both cKRnet and KRnet in this test problem are set to have $R=2$ and $L=12$ affine coupling layers. Each affine coupling layer contains a fully connected network that has two hidden layers with sixty four neurons in each hidden layer, and the ReLU function is employed as the activation function. All the models are trained using the Adam optimizer with default settings and a learning rate of 0.001. Again, the sample size of the reference results is set to $N_{\rm ref}=10^7$, and both CKR-DDUQ and KR-DDUQ processes are repeated thirty times to compute the average errors and effective sample sizes. Figure \ref{fig_rdd_error} presents the average mean and variance errors of CKR-DDUQ and KR-DDUQ. It is clear that each error reduces as the sample size increases, while all errors of CKR-DDUQ are smaller than those of KR-DDUQ with the same sample size. In this test problem, the average effective sample size is denoted by $\pE(N_{i,j}^{\text {eff}})$ for output $y_{i,j}$ (for $i,j=1,2$). Table \ref{table_rdd_ess} shows that for each output and $\Noff$, the average effective sample size of CKR-DDUQ is more than that of KR-DDUQ. From Figure \ref{fig_rdd_output_pdf}, it can be seen that the PDFs of $y_1(\xi,\beta)$ and $y_2(\xi,\beta)$ generated by CKR-DDUQ approach the reference PDFs as the sample size increases. 

\begin{table}[!ht]
    \centering
    \caption{Average effective sample sizes, random domain decomposition test problem.}
    \begin{tblr}
    {
      cells = {c},
      cell{2}{1} = {r=2}{},
      cell{4}{1} = {r=2}{},
      cell{6}{1} = {r=2}{},
      cell{8}{1} = {r=2}{},
      hline{1-2,4,6,8,10} = {-}{},
      hline{3,5,7,9} = {2-5}{},
    }
                         & Method   & $\Noff=10^4$ & $\Noff=10^5$ & $\Noff=10^6$  \\
    $\pE(N_{1,1}^{\text {eff}})$ & CKR-DDUQ & 308          & 391          & 2021  \\
                                 & KR-DDUQ  & 14           & 86           & 1123  \\               
    $\pE(N_{1,2}^{\text {eff}})$ & CKR-DDUQ & 232          & 246          & 505    \\
                                 & KR-DDUQ  & 7            & 21           & 180     \\                    
    $\pE(N_{2,1}^{\text {eff}})$ & CKR-DDUQ & 298          & 2414         & 23750   \\
                                 & KR-DDUQ  & 49           & 1102         & 18461    \\                   
    $\pE(N_{2,2}^{\text {eff}})$ & CKR-DDUQ & 1271         & 12908        & 128375   \\
                                 & KR-DDUQ  & 273          & 8010         & 109548                        
    \end{tblr}\label{table_rdd_ess}
\end{table}

\subsection{Two-component Stokes problem with \texorpdfstring{$19$}{19} parameters}\label{stokes_test}
The Stokes equations for this test problem are
\begin{eqnarray}
    \div \left(a\left(x,\xi\,\right)\nabla \vec{u}\left(x,\xi\,\right)\right)+\nabla p\left(x,\xi\,\right)=\vec{f}\left(x,\xi\,\right) & \quad \textrm{in}& D\times\Gamma, \label{st1} \\
    \div \vec{u}\left(x,\xi\,\right)=0 & \quad \textrm{in}& D\times\Gamma, \label{st2} \\
    \vec{u}\left(x,\xi\,\right)=\vec{g}\left(x,\xi\,\right) & \quad \textrm{on}& \partial D_D\times\Gamma, \label{st3} \\
    a\left(x,\xi\,\right)\frac{\partial \vec{u}\left(x,\xi\,\right)}{\partial n}- \vec{n}p=\vec{0} & \quad \textrm{on}& \partial D_N\times\Gamma, \label{st4} 
\end{eqnarray}
where $D\in \dsR^2$, $\vec{u}=[u_1,u_2]^{\top}$ and  $p$ are the flow velocity and the scalar pressure respectively, and $\vec{n}$ is the outward-pointing normal to the boundary. The boundary $\partial D$ is split into two parts: the Dirichlet boundary $\partial D_D$ and the Neumann boundary $\partial D_N$, such that $\partial D=\partial D_D \cup \partial D_N$ and $\partial D_D \cap \partial D_N=\emptyset$. $\vec{f}$ is the source function and $\vec{g}$ specifies the Dirichlet boundary condition. 
\begin{figure}[!htp]
    \setlength{\unitlength}{1.5cm}
    \centerline{
    \begin{picture}(6,2.3)(-0.2,-0.2)
    \put(0,0){\line(1,0){1}}
    \put(0,0){\line(0,1){2}}
    \put(0,2){\line(1,0){1}}
    \put(1,0){\line(0,1){0.75}}
    \put(1,2){\line(0,-1){0.75}}
    \put(5,0){\line(1,0){1}}
    \put(5,2){\line(1,0){1}}
    \put(5,0){\line(0,1){0.75}}
    \put(5,2){\line(0,-1){0.75}}
    \put(6,0){\line(0,1){2}}
    \put(1,0.75){\line(1,0){4}}
    \put(1,1.25){\line(1,0){4}}
    \put(0.3,1){$D_1$}
    \put(5.3,1){$D_2$}
    \put(0,0){\circle*{0.1}}
    \put(-0.3,-0.25){\scriptsize (0,0)}
    \put(1,0){\circle*{0.1}}
    \put(0.8,-0.25){\scriptsize   (1,0)}
    \put(1,0.75){\circle*{0.1}}
    \put(1.03,0.55){\scriptsize   (1,0.75)}
    \put(1,1.25){\circle*{0.1}}
    \put(1.03,1.3){\scriptsize   (1,1.25)}
    \put(1,2){\circle*{0.1}}
    \put(1.04,1.9){\scriptsize   (1,2)}
    \put(3,1.25){\circle*{0.1}}
    \put(3.03,1.3){\scriptsize   (3,1.25)}
    \put(5,1.25){\circle*{0.1}}
    \put(4.3,1.3){\scriptsize   (5,1.25)}
    \put(6,0){\circle*{0.1}}
    \put(5.8,-0.25){\scriptsize (6,0)}
    \put(-1.5,0.5){\vector(1,0){0.5}}
    \put(-1.5,0.5){\vector(0,1){0.5}}
    \put(-1.4,0.2){$x_1$}
    \put(-1.85,0.6){$x_2$}
    \linethickness{0.7mm}
    \put(3,0.75){\line(0,1){0.5}}
    \end{picture}
    }
    \caption{Illustration of the spatial domain with two components for the Stokes test problem. }
    \label{fig_2d_stokes}
\end{figure}
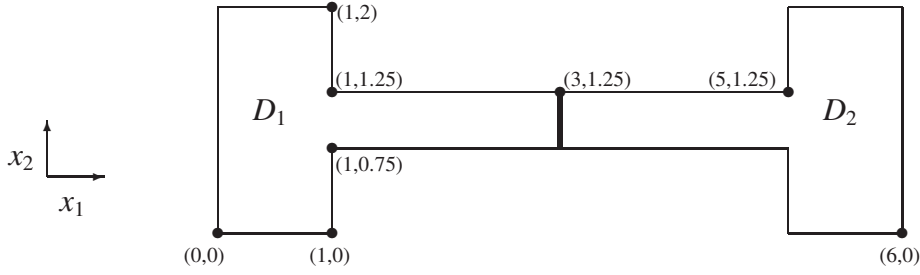
    
The physical domain considered in this test problem is shown in Figure \ref{fig_2d_stokes}, where the domain is decomposed into two components. The source function $\vec{f}$ is set to zero, and the Neumann boundary condition \eqref{st4} is applied on the outflow boundary $\partial D_N:= \{\{6\} \times {(0,2)}\}$, while other boundaries are the Dirichlet boundaries. The velocity profile $\vec{u}=\left[2x_2(2-x_2),0\right]^{\top}$ is imposed on the inflow boundary $\{\{0\}\times (0,2)\}$, and $\vec{u}=[0,0]^{\top}$ is imposed on all other Dirichlet boundaries. In \eqref{st1}, we focus on the situation that there exists uncertainties in the flow viscosity $a(x,\xi)$, which is again set to be a random field with covariance function \eqref{covariance} and KL expansion \eqref{kl_expan} posed on each local subdomain. Here, we set $a_{1,0}(x)=a_{2,0}(x)=1$ and $\sigma=0.5$. The covariance length of $a(x,\xi)|_{D_1}$ is set to $L_{c1}=2$ and that of $a(x,\xi)|_{D_2}$ is $L_{c2}=3$. The number of KL modes retained is $\mathfrak{N}_1=12$ and $\mathfrak{N}_2=7$, which results in a total of $19$ parameters in the input vector $\xi$ with $\xi_1\in\dsR^{12}$ and $\xi_2\in\dsR^{7}$. Each random variable in $\xi$ is set to follow the same distribution as that in Section~\ref{2domain_test}. The outputs of interest in this test problem are set to 
\begin{eqnarray}
    y_1(\xi)=\frac{1}{100}\int_{Z_1} p(x,\xi) \,{\rm d} x_2,\\
    y_2(\xi)=\frac{1}{10}\int_{Z_2} p(x,\xi) \,{\rm d} x_2,
\end{eqnarray}
where $x=[x_1,x_2]^{\top}$, $Z_1:=\{x|\> x_1=1.03125,\>0.75\leq x_2 \leq 1.25\}$ and $Z_2:=\{x|\> x_1=4.96875,\> 0.75\leq x_2 \leq 1.25\}$.

\begin{figure}[!htp]
    \centerline{
    \begin{tabular}{cc}
    \includegraphics[width=0.37\textwidth]{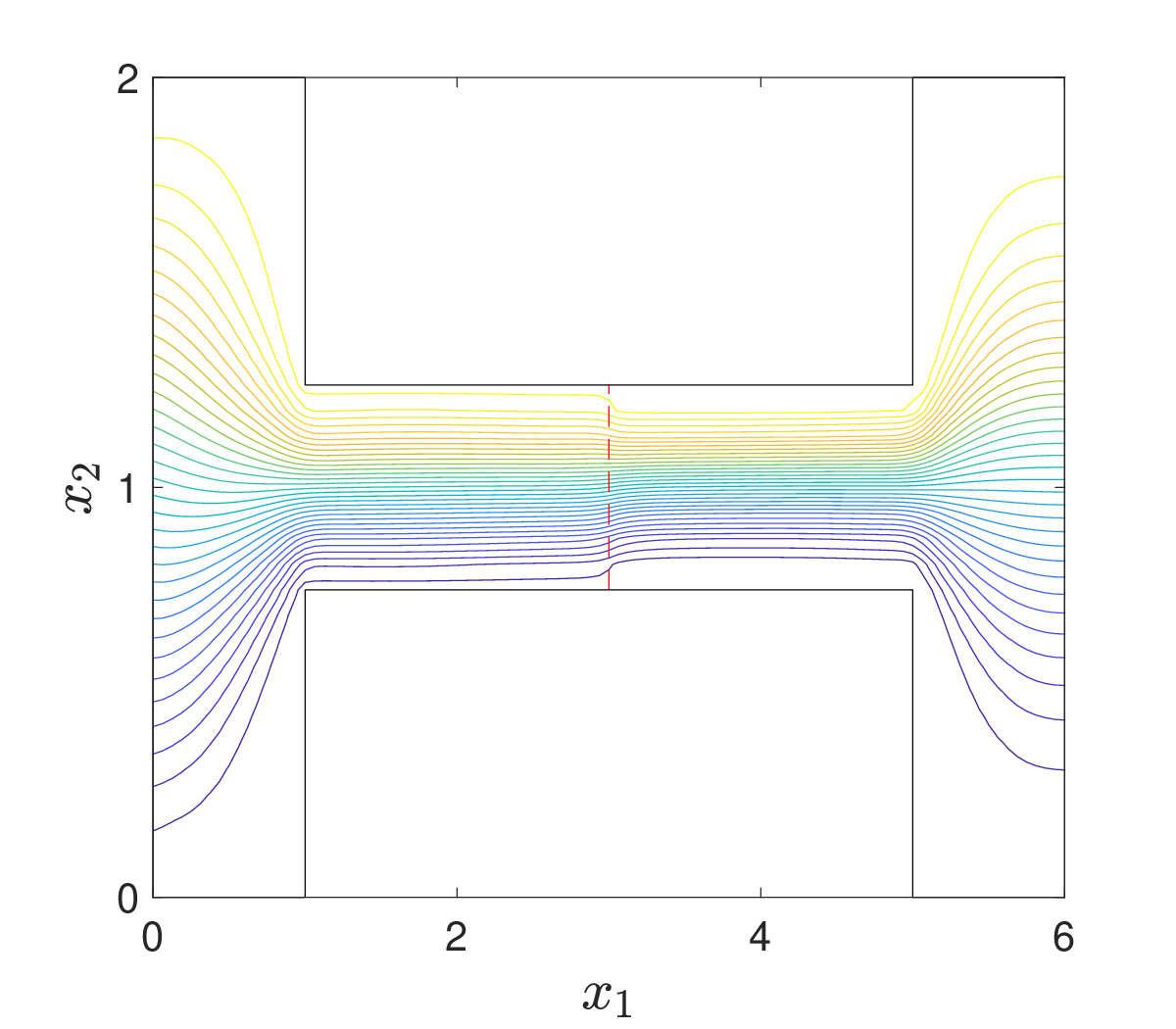}&
    \includegraphics[width=0.37\textwidth]{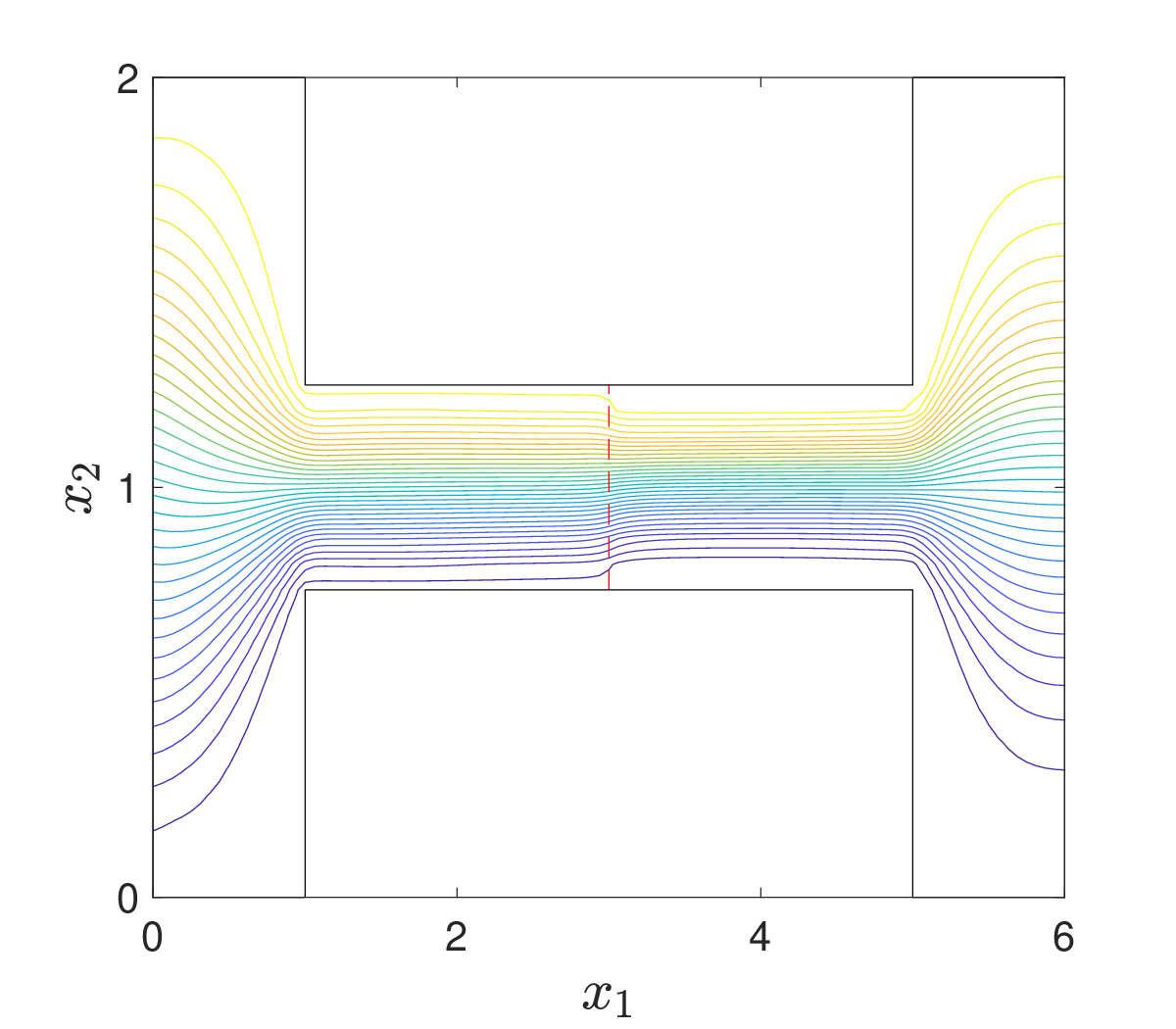}\\
    (a) Streamline: global finite element method& (b) Streamline: domain decomposition method\\
    \includegraphics[width=0.37\textwidth]{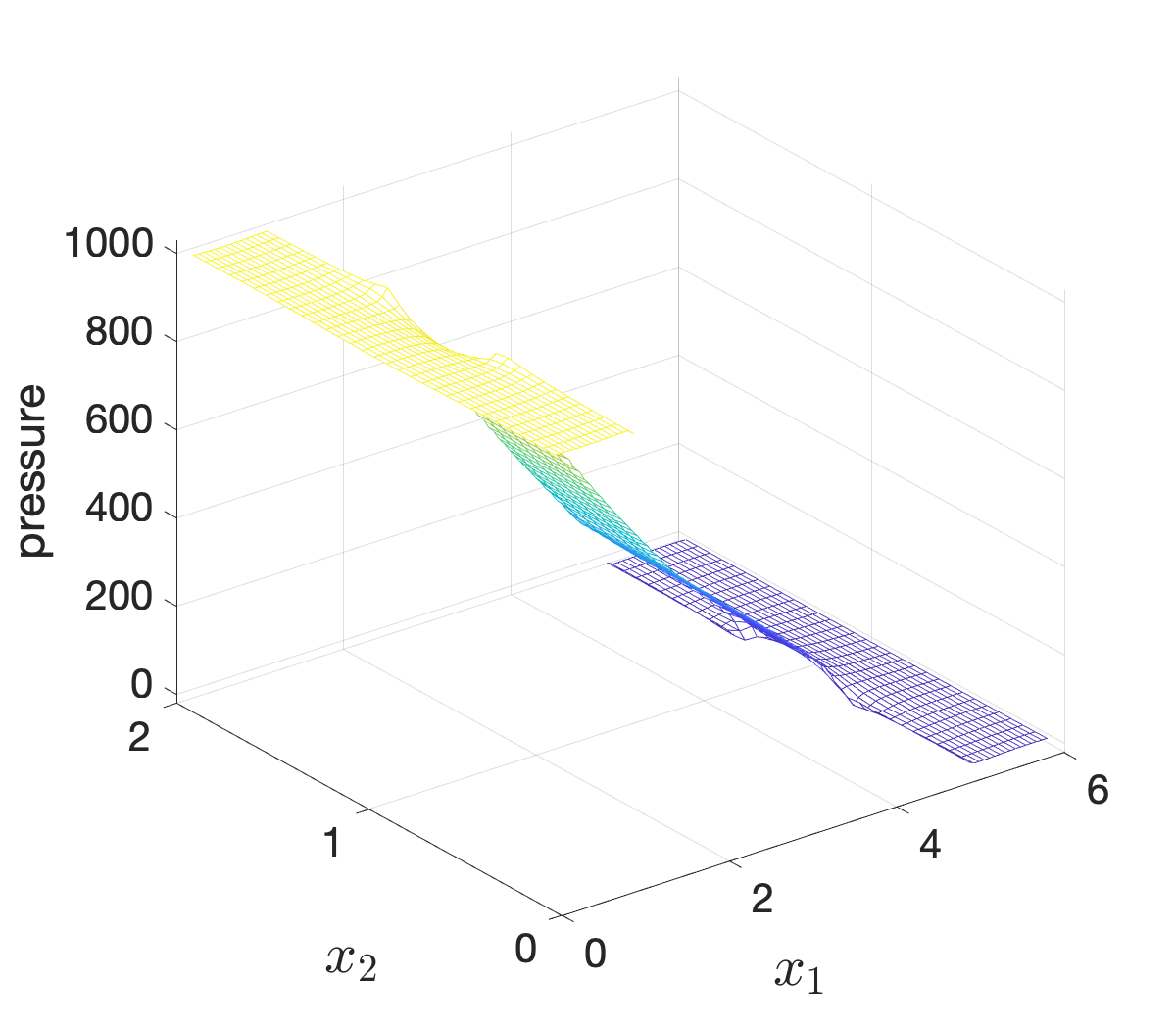}&
    \includegraphics[width=0.37\textwidth]{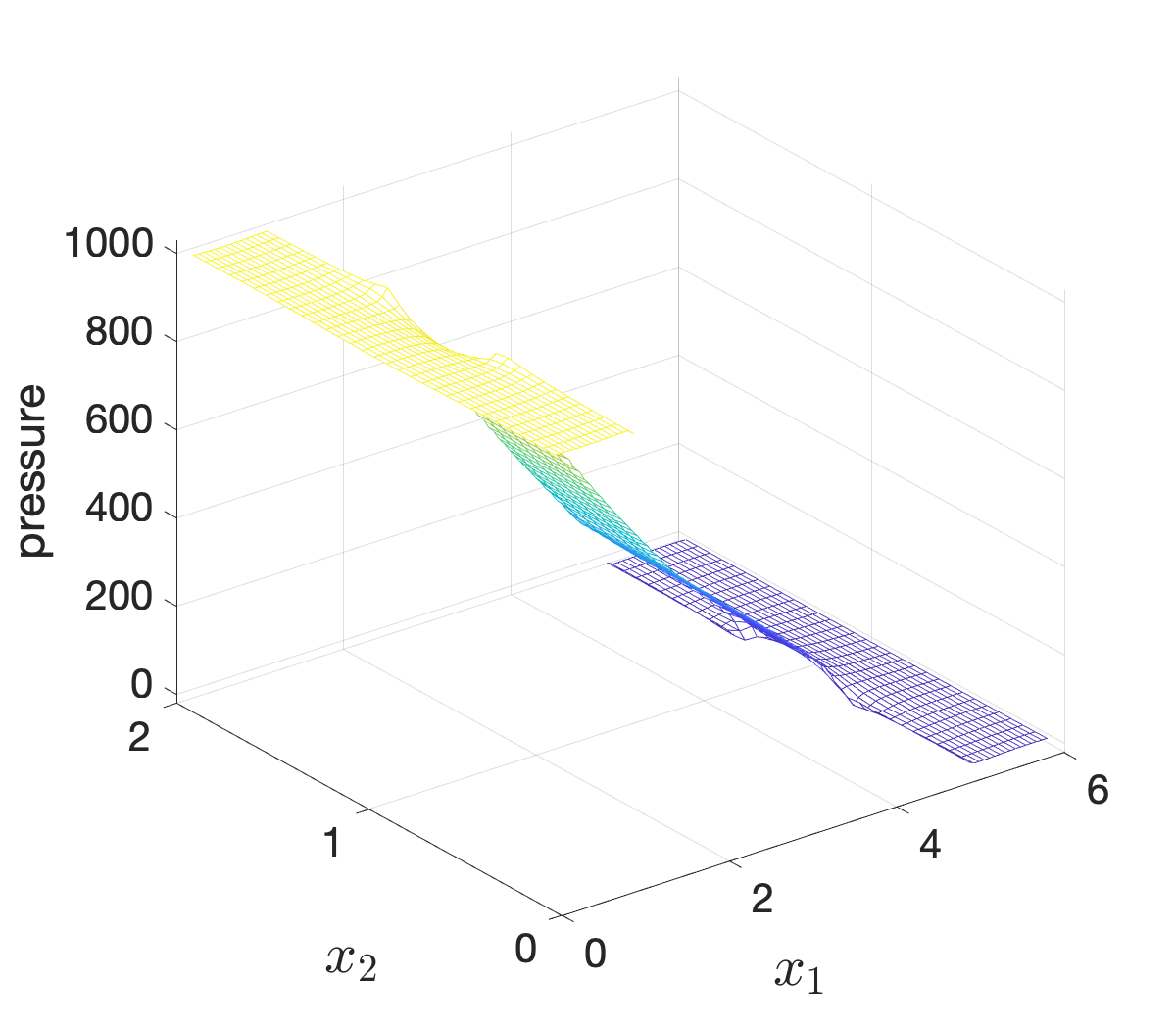}\\
    (c) Pressure: global finite element method& (d) Pressure: domain decomposition method
    \end{tabular}}
    \caption{\lr The global finite element solution and the domain decomposition solution responding to a given realization of $\xi$, Stokes test problem.}
    \label{fig_stokes_result_example}
\end{figure}

The Robin-type domain decomposition method based on the stabilized $\boldsymbol{Q}_1$--$\boldsymbol{P}_0$ (continuous bilinear velocity and discontinuous constant pressure on rectangles) mixed finite element \cite{feng2023robin} is used to decompose the global problem \eqref{st1}--\eqref{st4}, where a uniform square mesh with size $h=1/16$ is used. Figure \ref{fig_stokes_result_example} shows the flow streamlines and the pressure fields generated by both the global stabilized $\boldsymbol{Q}_1$--$\boldsymbol{P}_0$ method \cite{kechkar1992analysis,elman2014finite} and the domain decomposition version \cite{feng2023robin} responding to a given realization of $\xi$. It can be seen that there is no visual difference between the results of the two methods. Before the offline stage, one hundred samples of the global input parameter are generated, and the deterministic version of \eqref{st1}--\eqref{st4} are solved with the domain decomposition procedure to provide snapshots, which are used to construct POD bases of the interface functions. Three POD modes for $g_{2,1}$ and two POD modes for $g_{1,2}$ are retained, i.e., $\tau_1\in\dsR^{3}$ and $\tau_2\in\dsR^{2}$. Thus, the local inputs for subdomain $D_1$ and subdomain $D_2$ are $[\xi_1^{\top},\tau_1^{\top}]^{\top}\in\dsR^{15}$ and $[\xi_2^{\top},\tau_2^{\top}]^{\top}\in\dsR^{9}$ respectively. 

\begin{figure}[!htp]
    \centerline{
    \begin{tabular}{cc}
    \includegraphics[width=0.37\textwidth]{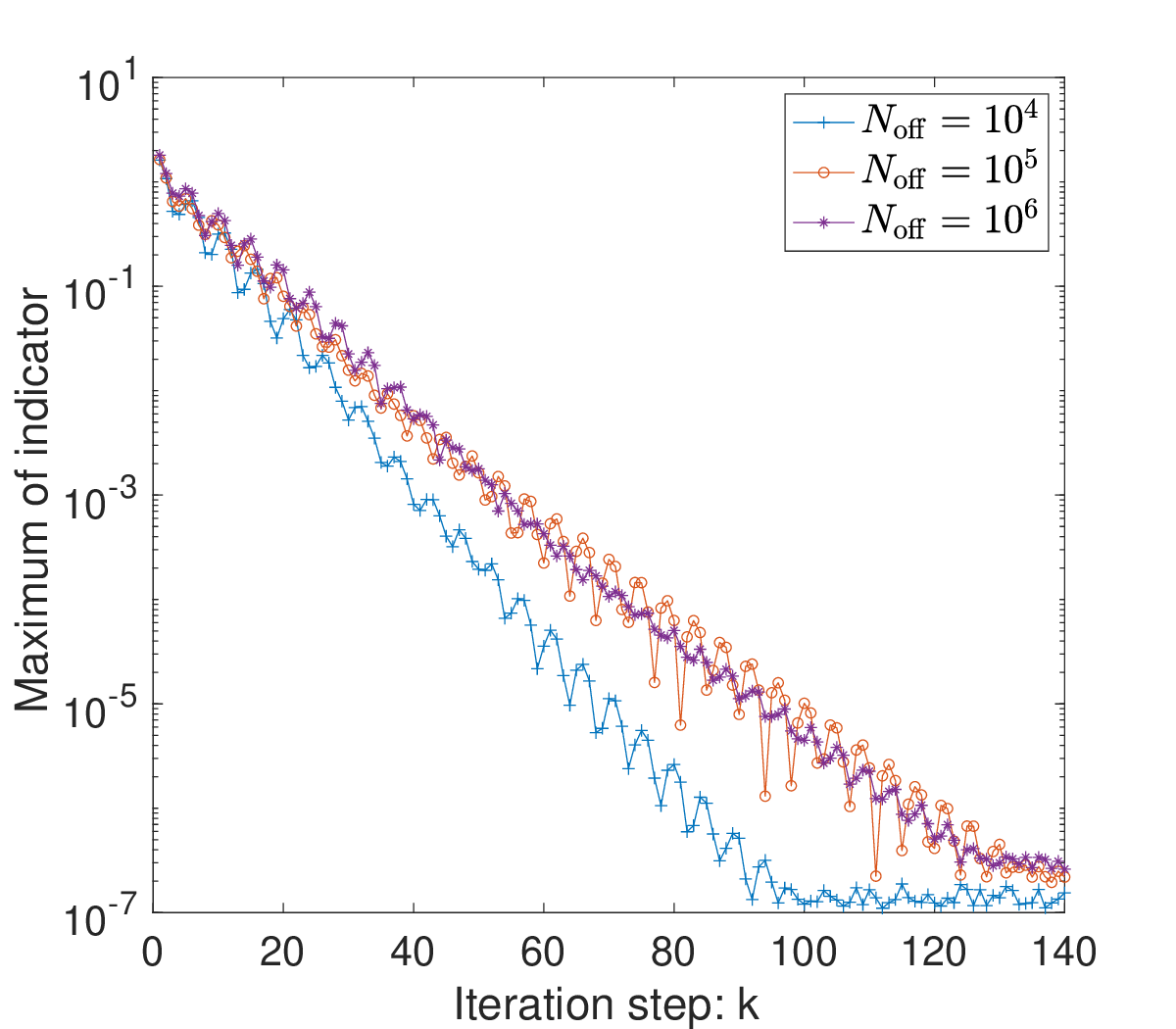}&
    \includegraphics[width=0.37\textwidth]{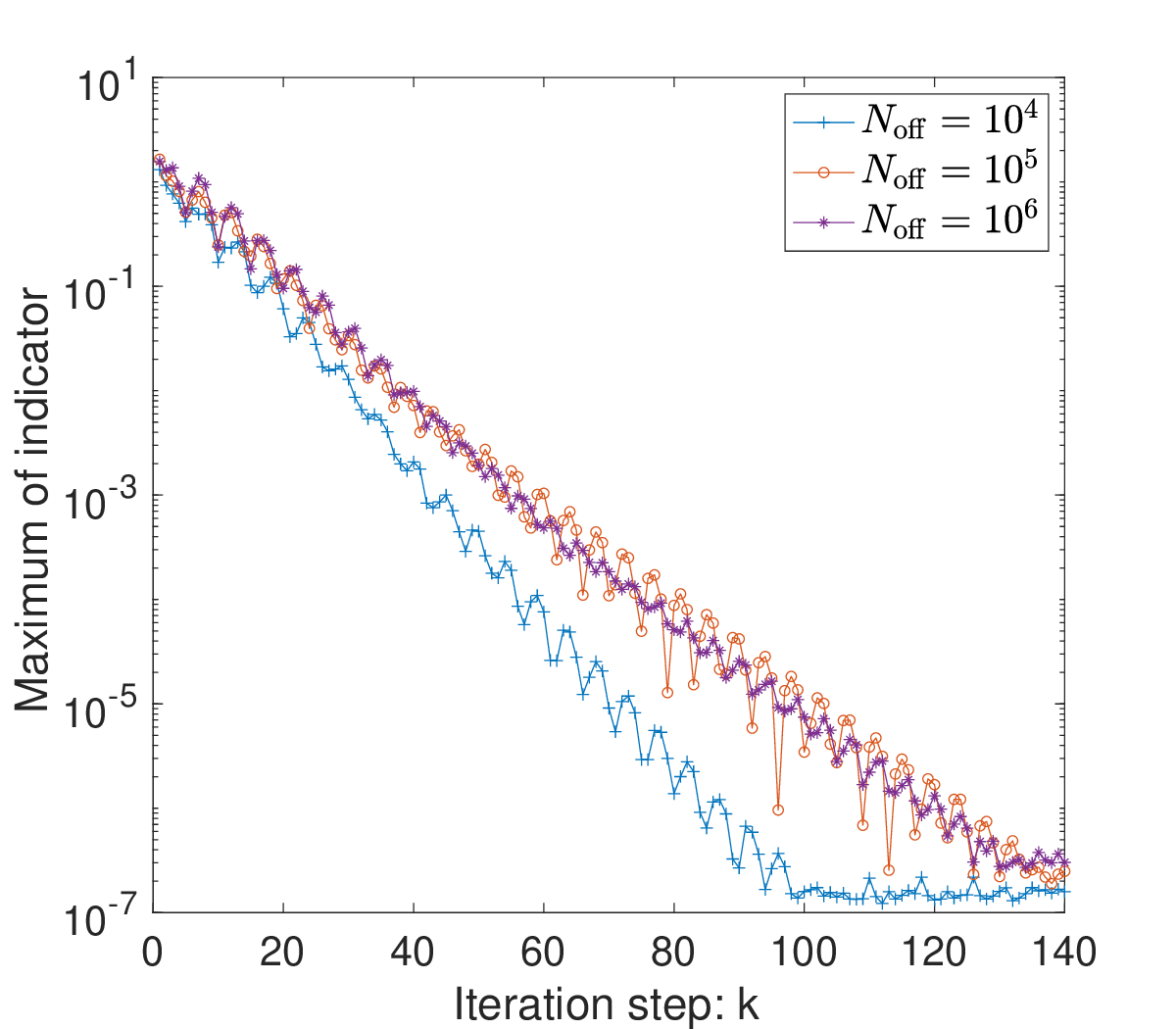}\\
    (a) Error indicator on $D_1$ & (b) Error indicator on $D_2$
    \end{tabular}}
    \caption{
    Maximum of the error indicator on subdomain $D_1$ ($\max_{s=1:{\Noff}}\|\tau^{k+1}_1(\xi^{(s)})-\tau^k_1(\xi^{(s)})\|_{\infty}$) and that on subdomain $D_2$ ($\max_{s=1:{\Noff}}\|\tau^{k+1}_2(\xi^{(s)})-\tau^k_2(\xi^{(s)})\|_{\infty}$) for the coupling surrogates, Stokes test problem.}
    \label{fig_stokes_surrogate_convergence}
\end{figure}

\begin{figure}[!htp]
    \centerline{
    \begin{tabular}{cc}
    \includegraphics[width=0.37\textwidth]{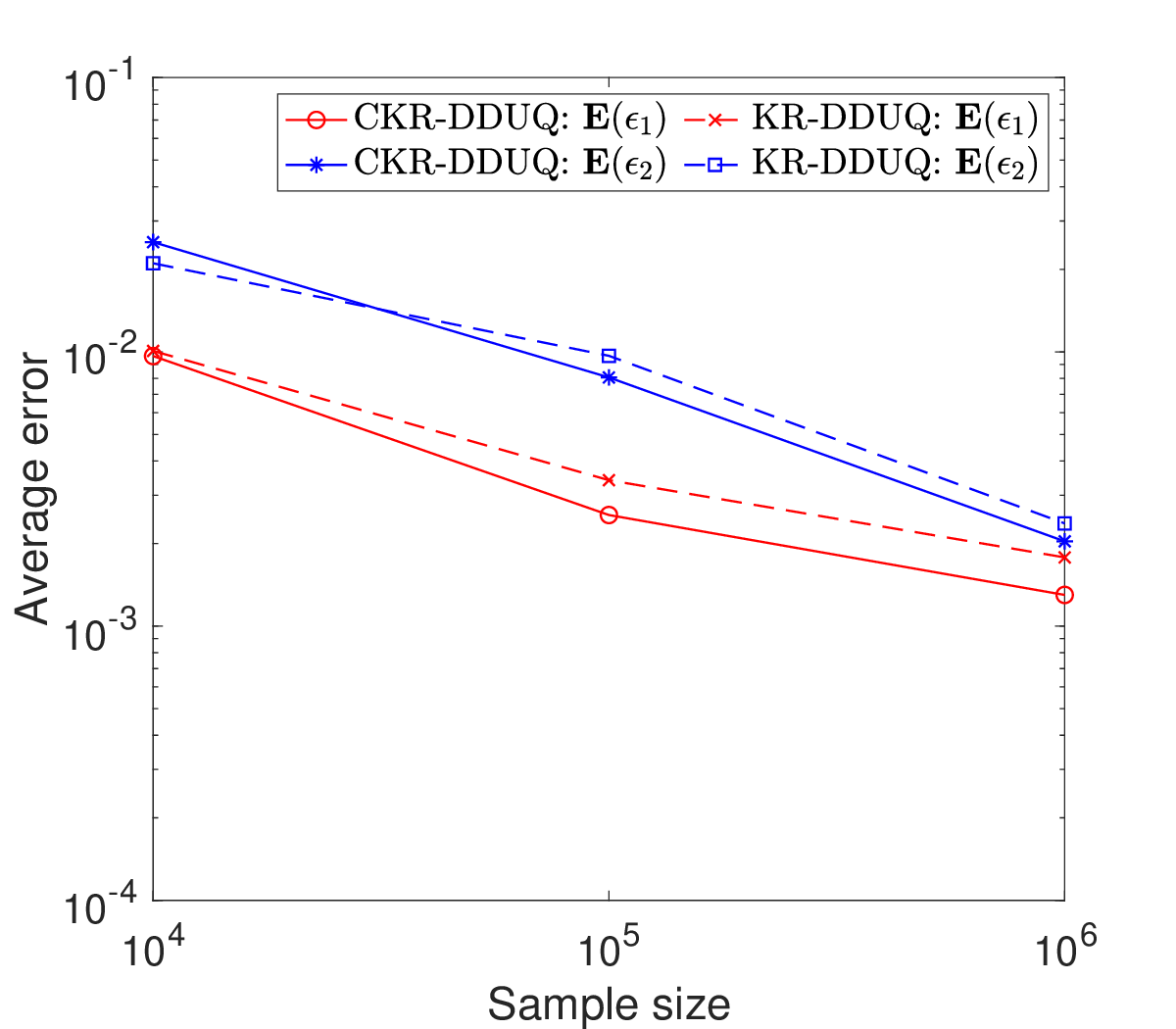}&
    \includegraphics[width=0.37\textwidth]{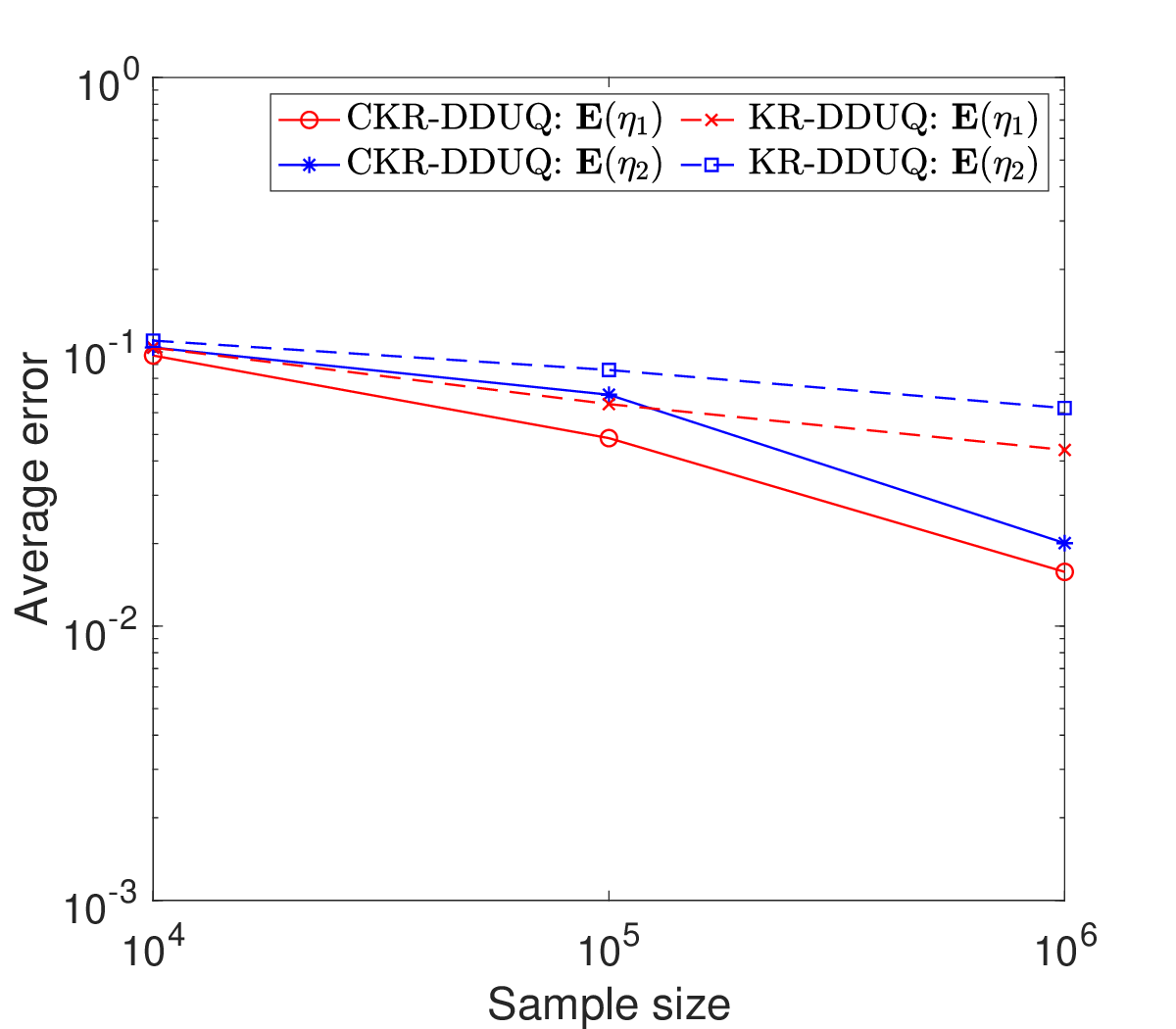}\\
    (a) Average mean errors & (b) Average variance errors
    \end{tabular}}
    \caption{\lr Average CKR-DDUQ and KR-DDUQ errors in mean and variance estimates for each output $y_i$ ($\pE(\epsilon_i)$ and $\pE(\eta_i))$, $i=1,2$, Stokes test problem.}
    \label{fig_stokes_error}
\end{figure}

\begin{figure}[!ht]
    \centerline{
    \begin{tabular}{cc}
    \includegraphics[width=0.37\textwidth]{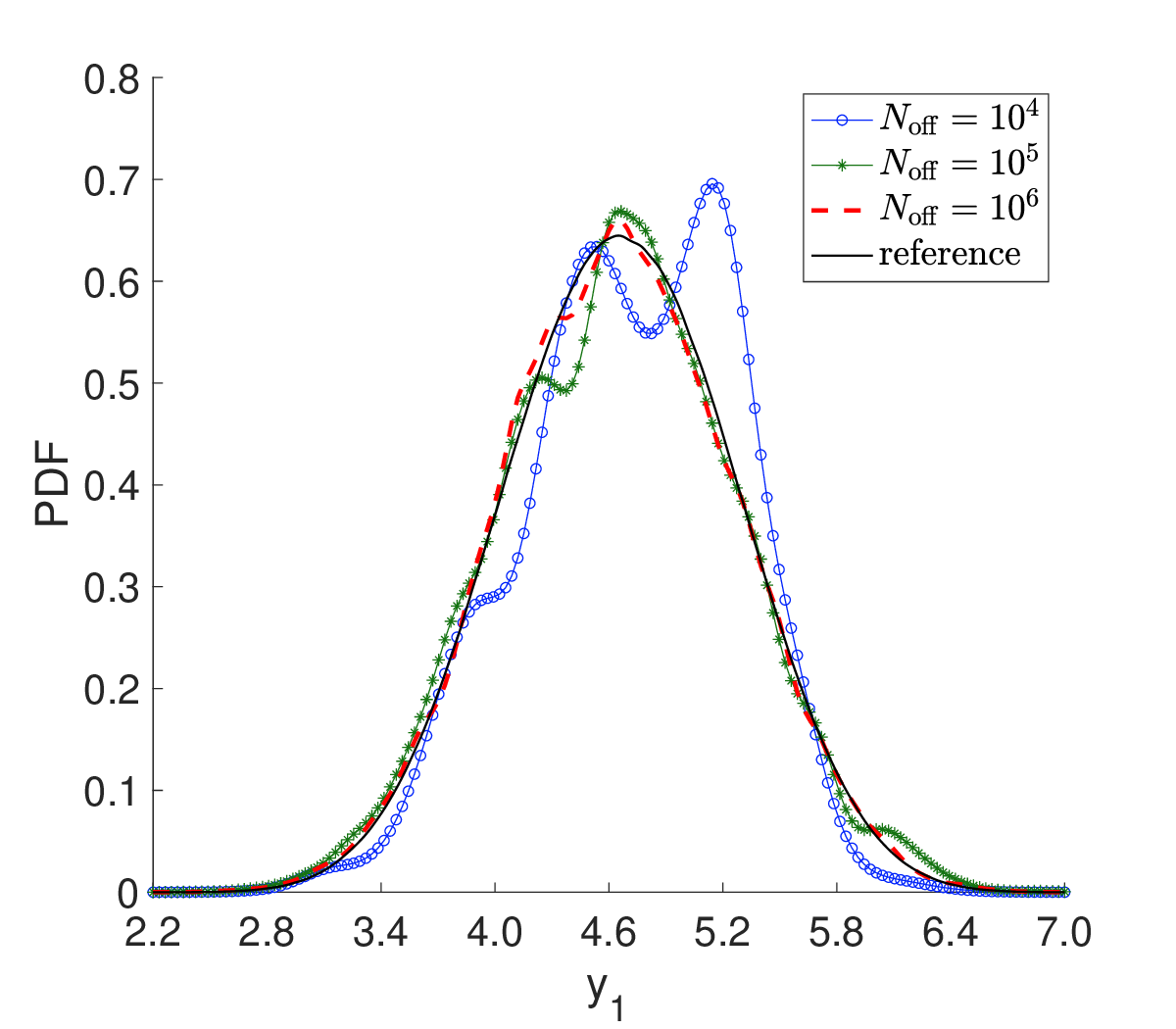}&
    \includegraphics[width=0.37\textwidth]{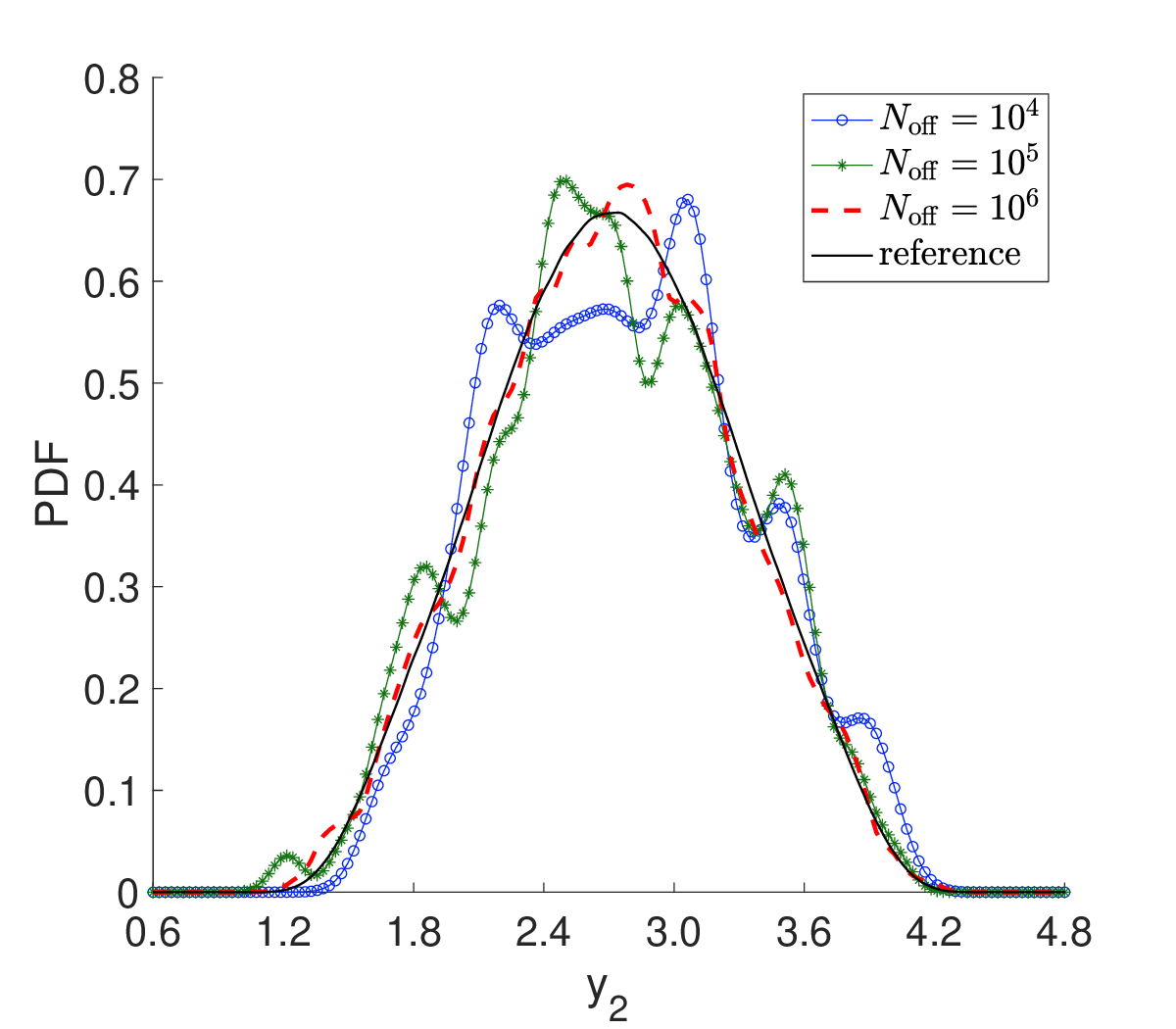}
    \end{tabular}}
    \caption{PDFs of the outputs of interest, Stokes test problem.}
    \label{fig_stokes_output_pdf}
\end{figure}

In the offline stage, two multivariate Gaussian distributions are used as the proposal PDFs $p_{\tau_1}(\tau_1)$ and $p_{\tau_2}(\tau_2)$, of which the mean vectors and the covariance matrices are estimated by the snapshots obtained in the POD procedure. We use ResNets with ten hidden layers and the LeakyReLU activation function to build the coupling surrogates, and each hidden layer contains sixty four neurons. In the online stage, $\Non$ is set to equal $\Noff$ and the threshold of domain decomposition iteration is specified as $tol=10^{-6}$. Figure \ref{fig_stokes_surrogate_convergence} shows the maximum of the error indicator (introduced in Section \ref{2domain_test}) with respect to domain decomposition iterations using these coupling surrogates for $\Noff=10^4,10^5$ and $10^6$, where it is clear that the maximum of the error indicator reduces exponentially as the iteration step $k$ increases for each subdomain.
In the density estimation step, both cKRnet and KRnet are set to $R = 3$ for subdomain $D_1$ and $R=2$ for subdomain $D_2$. We take $L=4$ affine coupling layers for both cKRnet and KRnet, and each affine coupling layer is set to contain a fully connected network with two hidden layers and thirty two neurons in each hidden layer. The activation function is set as the ReLU function here. Both cKRnet and KRnet are trained by the Adam optimizer with default settings and the learning rate is $0.001$. After generating a set of reference output samples with $N_{\rm ref}=10^7$ through solving the global problem with stabilized $\boldsymbol{Q}_1$--$\boldsymbol{P}_0$, we compute the average mean and variance errors of CKR-DDUQ and KR-DDUQ with thirty repeats, and compute the average errors and the average effective sample sizes. The errors of the outputs are shown in Figure \ref{fig_stokes_error}, where it is clear that all errors decrease as the sample size increases, and CKR-DDUQ always gives smaller errors than KR-DDUQ. From Table \ref{table_stokes_ess}, it can be seen that the average effective sample sizes of CKR-DDUQ are larger than those of KR-DDUQ with the same $\Noff$ value. The  PDFs of outputs in this test problem are shown in Figure \ref{fig_stokes_output_pdf}, where it is clear that as $\Noff$ increases, the PDFs generated by CKR-DDUQ approach the reference PDFs.

\begin{table}[!ht]
    \centering
    \caption{Average effective sample sizes, Stokes test problem.}
    \begin{tblr}
    {
      cells = {c},
      cell{2}{1} = {r=2}{},
      cell{4}{1} = {r=2}{},
      hline{1-2,4,6} = {-}{},
      hline{3,5} = {2-5}{},
    }
               & Method   & $\Noff=10^4$ & $\Noff=10^5$ & $\Noff=10^6$ \\
    $\pE(N_1^{\text {eff }})$  & CKR-DDUQ & 92          & 884         & 8462        \\
                               & KR-DDUQ  & 90          & 797         & 8132        \\            
    $\pE(N_2^{\text {eff }})$  & CKR-DDUQ & 58          & 365         & 3911         \\
                               & KR-DDUQ  & 55          & 352         & 3651        
    \end{tblr}\label{table_stokes_ess}
\end{table}

%% file: sections/section6.tex
In this paper we present a CKR-DDUQ approach to propagate uncertainties across different physical subdomains. It targets improving the original DDUQ approach through developing new density estimation techniques. In the DDUQ online stage, the target local input PDFs need to be estimated using target samples. While the target samples are obtained through domain decomposition iterations with coupling surrogates (no PDE solve is required), it is cheap to obtain a large number of target samples. However, these local inputs can be high-dimensional, and estimating their joint PDFs in DDUQ is a bottleneck. In this work, we utilize the structure of the target inputs $[\xi_i^\top,\tau_i^\top]^{\top}$ (see Section \ref{dd}), where elements of $\xi_i$  are typically assumed to be independent with each other but each target input interface parameter $\tau^{\infty}_i$ (see Section \ref{dd}) is dependent on $\xi_i$, and then develop cKRnet to estimate the conditional PDFs $\pi_{\tau_i\mid\xi_i}(\tau_i\mid\xi_i)$ for $i=1,\dots,M$.
Estimating conditional PDFs reduces the complexity compared to estimating joint PDFs---cKRnet only treats $\tau_i$ as the unknown random vector and considers $\xi_i$ as a conditional input, and then the dimension of random variables whose PDFs need to be estimated is reduced. The convergence analysis of CKR-DDUQ is conducted, and numerical results show the effectiveness of both cKRnet for conditional density estimation and the overall CKR-DDUQ approach. As POD is still applied for dimensional reduction of the interface parameters, this may not be optimal when the coupling functions have highly nonlinear structures. A possible solution is to utilize new deep learning based dimension reduction methods for the coupling surrogates, e.g., the VAE methods, which will be the focus of our future work.